%% file: pure-metric.tex
\documentclass[twoside]{book}

\usepackage{lectures}
\usepackage[colorlinks=true,
citecolor=black,
linkcolor=black,
anchorcolor=black,
filecolor=black,
menucolor=black,
urlcolor=black,
pdftitle={Pure metric geometry: introductory lectures},
pdfsubject={Geometry},
pdfauthor={Anton Petrunin}
]{hyperref}
\makeindex

\begin{document}

\title{Pure metric geometry:\\
introductory lectures}
\author{Anton Petrunin}
\date{}
\maketitle

\section*{Preface}

This text can serve as an introductory part for a variety of courses in metric geometry.
Here is a graph of essential dependencies of the lectures; some statements (mostly exercises) add more dependencies, but they can be ignored.
\begin{figure}[!ht]
\centering
\begin{tikzpicture}[->,>=stealth',shorten >=1pt,auto,scale=1.4,
  thick,main node/.style={circle,draw,font=\sffamily\bfseries,minimum size=8mm}]

  \node[main node] (1) at (1,0) {\ref{chap:defs}};
  \node[main node] (2) at (.5,-5/6){\ref{chap:urysohn}};
  \node[main node] (3) at (1.5,-5/6) {\ref{chap:injective}};
  \node[main node] (4) at (2,0) {\ref{chap:hausdorff}};
  \node[main node] (5) at (3,0) {\ref{chap:GH}};
  \node[main node] (6) at (4,0) {\ref{chap:ultralimits}};

  \path[every node/.style={font=\sffamily\small}]
   (1) edge node{}(2)
   (1) edge node{}(3)
   (1) edge node{}(4)
   (4) edge node{}(5)
   (5) edge node{}(6);
\end{tikzpicture}
\end{figure}
The necessary definitions are introduced in (\ref{chap:defs}).
In (\ref{chap:urysohn}) we discuss the Urysohn space.
In (\ref{chap:injective}) we discuss injective spaces.
In (\ref{chap:hausdorff}) we introduce the Hausdorff metric.
In (\ref{chap:GH}) and (\ref{chap:ultralimits}) we discuss two types of convergences of metric spaces --- the Gromov--Hausdorff limit and ultralimit.

Applications are given only as illustrations.
We stick to domestic affairs of metric spaces, keeping away from any extra structure. 
(Adding an extra structure brings an extra tool and often opens a huge field for development.
The examples include Alexandrov geometry,
geometric group theory,
metric-measure spaces and optimal transport.)

The only prerequisite is interest in the subject,
but any knowledge of classical geometry, differential geometry, topology, and real analysis will be useful. 

These notes are based on the minicourse given at SPbSU (Fall 2022) and the introductory part of a course at PSU (Spring 2020).
The latter included additional material from \cite{alexander-kapovitch-petrunin-2019,petrunin2020mnfld,nabutovsky}.
A part of the text is a compilation from \cite{alexander-kapovitch-petrunin-2019, alexander-kapovitch-petrunin-2025, petrunin-yashinski, petrunin-2022-PIGTIKAL, petrunin-zamorabarrera} and its drafts.

I want to thank
Sergei Ivanov,
Urs Lang,
Alexander Lytchak,
Rostislav Matveyev,
Julien Melleray,
and Sergio Zamora Barrera for help.
The present work is partially supported by NSF grant DMS-2005279,
the Simons Foundation grant \#584781,
and Minobrnauki of Russia, grant \#075-15-2022-289.

\thispagestyle{empty}
\tableofcontents
\thispagestyle{empty}

\include{metric}

\include{uryson}
\include{injective}
\include{converge}
\include{ultralimit}

\backmatter

\newgeometry{top=0.9in, bottom=0.9in,inner=0.5in, outer=0.5in}
\chapter{Semisolutions}

{

\footnotesize
\begin{multicols}{2}

\input{metric-sol}
\input{uryson-sol}
\input{injective-sol}
\input{converge-sol}
\input{ultralimit-sol}

\end{multicols}
}

\newgeometry{top=0.9in, bottom=0.9in,left=0.9in, right=0.9in, paperwidth=6in, paperheight=9in}

{\small\sloppy
\input{pure-metric.ind}

\def\emph{\textit}

\printbibliography[heading=bibintoc]
\fussy
}

\end{document}

%% file: metric.tex
\chapter{Definitions}\label{chap:defs}

In this lecture, we remind several definitions related to metric spaces and fix some conventions.

This lecture is self-contained, but it is written for students with some prior knowledge of metric spaces;
an introduction to general topology is sufficient but not necessary.
For a more detailed introduction, we recommend the first couple of chapters in the book by Dmitri Burago, Yuri Burago, and Sergei Ivanov \cite{burago-burago-ivanov}.

\section{Metric spaces}
\label{sec:metric spaces}

The distance between two points $x$ and $y$ in a metric space $\spc{X}$ will be denoted by $\dist{x}{y}{}$ or $\dist{x}{y}{\spc{X}}$.
The latter notation is used if we need to emphasize 
that the distance is taken in the space~${\spc{X}}$.

Let us recall the definition of metric. 

\begin{thm}{Definition}\label{def:metric}
A \index{metric}\emph{metric} on a set $\spc{X}$ is a real-valued function $(x,y)\z\mapsto\dist{x}{y}{\spc{X}}$ that satisfies the following conditions for any $x,y,z\z\in \spc{X}$:

\begin{subthm}{metric>=0}
$\dist{x}{y}{\spc{X}}\ge 0$,
\end{subthm}

\begin{subthm}{metric=0} $\dist{x}{y}{\spc{X}}= 0$ $\iff$ $x=y$,
\end{subthm}

\begin{subthm}{metric:sym} $\dist{x}{y}{\spc{X}}=\dist{y}{x}{\spc{X}}$,
\end{subthm}

\begin{subthm}{metric:triangle} $\dist{x}{y}{\spc{X}}+\dist{y}{z}{\spc{X}}\ge\dist{x}{z}{\spc{X}}$.
\end{subthm}

\end{thm}

Recall that a \index{metric space}\emph{metric space} is a set with a metric on it.
The elements of the set are called \index{point}\emph{points}. 
Most of the time we keep the same notation for the metric space and its underlying set;
the latter can be denoted by $\ushort{\spc{X}}$ if needed.

Given radius $R\in[0,\infty]$ and center $x\in \spc{X}$, the sets
\begin{align*}
\oBall(x,R)&=\set{y\in \spc{X}}{\dist{x}{y}{}<R},
\\
\cBall[x,R]&=\set{y\in \spc{X}}{\dist{x}{y}{}\le R}
\end{align*}
are called, respectively, the  \index{open ball}\emph{open} and  the \index{closed ball}\emph{closed  balls}.
The notations $\oBall(x,R)_{\spc{X}}$ and $\cBall[x,R]_{\spc{X}}$
might be used if we need to emphasize that these balls are taken in the metric space $\spc{X}$.

\begin{thm}{Exercise}\label{ex:quad-inq}
Show that the following inequality
\[\dist{p}{q}{\spc{X}}+\dist{x}{y}{\spc{X}}\le\dist{p}{x}{\spc{X}}+\dist{p}{y}{\spc{X}}+\dist{q}{x}{\spc{X}}+\dist{q}{y}{\spc{X}}\]
holds for any four points $p$, $q$, $x$, and $y$ in a metric space $\spc{X}$.
\end{thm}

\section{Topology}

The standard calculus definitions of \index{closed set}\emph{closed} and \index{open set}\emph{open sets}, \index{continuous function}\emph{continuous functions}, and \index{converging sequence}\emph{converging sequences} admit straightforward generalizations in the context of metric spaces.

\begin{thm}{Exercise}\label{ex:cont-dist}
Let $x$ be a point in a metric space $\spc{X}$.
Show that the \index{distance function}\emph{distance function} $\distfun_x\:\spc{X}\to \RR$ defined by
\[\distfun_x\:y\mapsto \dist{x}{y}{}\]
is continuous.
\end{thm}

\begin{thm}{Exercise}\label{ex:normal}
Let $A$ and $B$ be two disjoint closed sets in a metric space $\spc{X}$.
Construct a continuous function $f\:\spc{X}\to [0,1]$ such that $A=f^{-1}\{0\}$ and $B=f^{-1}\{1\}$.
\end{thm}

\begin{thm}{Advanced exercise}\label{ex:tietze}
Let $f\:A\to\RR$ be a continuous function defined on a closed set $A$ in a metric space $\spc{X}$.
Show that it admits a continuous extension to the whole space;
that is, there is a continuous function $F\:\spc{X}\to\RR$ such that $F(a)=f(a)$ for any $a\in A$.
\end{thm}

\section{Variations}

\parbf{Pseudometris.}
A metric for which the distance between two distinct points can be zero is called a \index{semimetric}\emph{semimetric} (also known as \index{pseudometric}\emph{pseudometric}).
In other words, to define semimetric, we need to remove condition \ref{SHORT.metric=0} from \ref{def:metric}.

Assume $\spc{X}$ is a semimetric space.
Consider an equivalence relation $\sim$ on $\spc{X}$ defined by
\[x\sim y\quad\iff\quad\dist{x}{y}{}=0.\] 
Note that if $x\sim x'$, then $\dist{y}{x}{}=\dist{y}{x'}{}$ for any $y\in\spc{X}$.
Thus, $\dist{*}{*}{}$ defines a metric on the
quotient $\spc{X}/{\sim}$.
The so-obtained metric space, say $\spc{X}'$, is called the 
\emph{corresponding metric space} for the semimetric space $\spc{X}$.

This construction shows that nearly any question about semimetric spaces can be reduced to a question about genuine metric spaces.
Often we do not distinguish between a semimetric space $\spc{X}$ and its corresponding metric space $\spc{X}'$.

\parbf{$\bm{\infty}$-metrics.}
One may also consider metrics with values in $[0,\infty]$;
that is, we allow infinite distance between points.
We might call them \index{metric!$\infty$-metric}\emph{$\infty$-metrics}, but most of the time we use the term \textit{metric}.

The following construction shows how to reduce questions about $\infty$-metrics to genuine metrics. 

Let 
\[x\approx y\quad\iff\quad\dist{x}{y}{}<\infty;\]
it defines another equivalence relation on $\spc{X}$.
The equivalence class of a point $x\in\spc{X}$ will be called the \index{metric component}\emph{metric component} 
 of $x$; it will be denoted by $\spc{X}_x$.
Note that 
\[\spc{X}_x=\oBall(x,\infty)_{\spc{X}};\]
that is, the metric component of $x$ is the open ball centered at $x$ and radius $\infty$.

If $\{\spc{X}_\alpha\}$ is a collection of metric spaces, then \index{disjoint union}\emph{disjoint union} $\bm{X}\z=\bigsqcup_\alpha\spc{X}_\alpha$ will be considered with a natural metric defined by
\[\dist{x}{y}{\bm{X}}
\df
\begin{cases}
\dist{x}{y}{\spc{X}_\alpha}&\text{if}\quad x,y\in \spc{X}_\alpha\quad\text{for some}\quad \alpha,
\\
\hfil \infty&\text{otherwise.}
\end{cases}
\]
It follows that any $\infty$-metric space is a disjoint union of genuine metric spaces --- the metric components of the original $\infty$-metric space.

\begin{thm}{Exercise}\label{ex:pseudo-infty-metric}
Given two sets $A$ and $B$ on the plane, set 
\[\dist{A}{B}{}=\mu(A\bigtriangleup B),\]
where $\mu$ denotes the Lebesgue measure and $\bigtriangleup$ denotes symmetric difference
\[A\bigtriangleup B
\df(A\cup B)\setminus(B\cap A)
=(A\setminus B)\cup(B\setminus A).\]

\begin{subthm}{ex:pseudo-infty-metric:pseudo}
Show that $\dist{*}{*}{}$ is a semimetric on the set of bounded closed subsets.
\end{subthm}

\begin{subthm}{ex:pseudo-infty-metric:infty}
Show that $\dist{*}{*}{}$ is an $\infty$-metric on the set of all open subsets.
\end{subthm}
\end{thm}

\section{Maximal metric and gluing}
\label{sec:max+glue}

\parbf{Maximal metric.}
Let $\{\,\dist{\ }{\ }{\alpha}\,\}$ be a family of $\infty$-semimetrics on a fixed set.
Observe that
\[\dist{x}{y}{}\df\sup_\alpha\{\,\dist{x}{y}{\alpha}\,\}\]
defines an $\infty$-semimetric; it is called the \index{maximal metric}\emph{maximal metric} of the family.

\parbf{Gluing.}
Suppose $\sim$ is an equivalence relation on an $\infty$-semimetric space $\spc{X}$.
Given $x\in\spc{X}$, 
denote by $[x]$ its equivalence class in the quotient $\spc{X}/\sim$.
Consider all $\infty$-semimetrics $|\ -\ |_\alpha$ on $\spc{X}/\sim$ such that the maps $\spc{X}\to \spc{X}/\sim$ defined by $x\mapsto [x]$ is \index{short map}\emph{short};
that is,
\[|[x]-[x']|_\alpha\le |x-x'|_{\spc{X}}\]
for any $x,x'\in \spc{X}$.
Let us equip $\spc{X}/\sim$ with the maximal metric of this family; in general, it is an $\infty$-semimetric.
The space $\spc{Z}$ that corresponds to the obtained $\infty$-semimetric space is called \index{gluing}\emph{gluing} of $\spc{X}$ along $\sim$.

This definition can be applied to a disjoint union of spaces; this way we can glue an arbitrary collection of metric spaces.

Note that any partially defined map $\phi$ from $\spc{X}$ to $\spc{Y}$ defines a minimal equivalence relation on $\spc{X}\sqcup\spc{Y}$ such that $x\sim \phi(x)$;
the corresponding gluing space is called \emph{gluing} along $\phi$.

The following exercise shows that metric gluing and the corresponding topological gluing
might have different topologies.
 
\begin{thm}{Exercise}\label{ex:gluing}
Construct a homeomorphism $\phi\:[0,1]\to [0,1]$ such that gluing of two unit intervals $[0,1]$ along $\phi$ is a one-point metric space.
\end{thm}

\section{Completeness}

A metric space $\spc{X}$ is called \index{complete space}\emph{complete} if every Cauchy sequence of points in $\spc{X}$ converges in $\spc{X}$.

\begin{thm}{Exercise}\label{ex:almost-min}
Suppose that $\rho$ is a positive continuous function on a complete metric space $\spc{X}$ and $\eps>0$.
Show that there is a point $x\in \spc{X}$ such that 
\[\rho(x)<(1+\eps)\cdot\rho(y)\]
for any point $y\in \oBall(x,\rho(x))$.
\end{thm}

Most of the time we will assume that a metric space is complete.
The following construction produces a complete metric space $\bar{\spc{X}}$ for any given metric space $\spc{X}$.

\parbf{Completion.}
Given a metric space $\spc{X}$, 
consider the set $\spc{C}$ of all Cauchy sequences in $\spc{X}$.
Note that for any two Cauchy sequences $(x_n)$ and $(y_n)$ the right-hand side in \ref{eq:cauchy-dist} is defined;
moreover, it defines a semimetric on~$\spc{C}$
\[\dist{(x_n)}{(y_n)}{\spc{C}}\df\lim_{n\to\infty}\dist{x_n}{y_n}{\spc{X}}.\eqlbl{eq:cauchy-dist}\]
The corresponding metric space is called the \index{completion}\emph{completion} of $\spc{X}$;
it will be denoted by $\bar{\spc{X}}$.
  
For each point $x\in\spc{X}$, one can consider a constant sequence $x_n=x$ which is Cauchy.
It defines a natural inclusion map $\spc{X}\hookrightarrow \bar{\spc{X}}$.
It is easy to check that this map is distance-preserving.
In particular, we can (and will) consider $\spc{X}$ as a subset of $\bar{\spc{X}}$.

Note that $\spc{X}$ is a dense subset in its completion $\bar{\spc{X}}$.

\begin{thm}{Exercise}\label{ex:complete-completion}
Show that the completion of a metric space is complete.
\end{thm} 

\section{G-delta sets}\label{sec:G-delta}

\begin{thm}{Baire's theorem}\label{thm:baire}
For any sequence $\Omega_1,\Omega_2,\dots$ of open dense subsets in a complete metric space, the intersection $\bigcap_{n\in \NN}\Omega_n$ is dense.
\end{thm}

A subset is called a \index{G-delta set}\emph{G-delta} if it can be presented as an intersection of a countable number of open subsets.
Note that by Baire's theorem, a countable intersection of dense G-delta sets is a dense G-delta set --- in particular it is nonempty.
Therefore we are allowed to say that a dense G-delta set contains \textit{most} of the points in a complete metric space.

\parit{Proof.}
We may assume that the space is nonempty; otherwise, there is nothing to prove.

Given a  closed ball $\cBall[p_0,R_0]$,
let us apply induction to construct a nested sequence of closed balls
\[\cBall[p_0,R_0]\supset\cBall[p_1,R_1]\supset \cBall[p_2,R_2]\supset\dots\]
such that $\cBall[p_n,R_n]\subset \Omega_n$ and $R_n>0$ for each $n\ge 1$. 
Assume $\cBall[p_{n-1},R_{n-1}]$ is already constructed.
Since $\Omega_n$ is dense we can choose a closed ball 
$\cBall[p_n,R_n]\subset\Omega_n\cap \cBall[p_{n-1},R_{n-1}]$.

Note that we can assume that $R_n<\tfrac1{n}$ for each $n\z\ge 1$.
In this case, the sequence $p_1$, $p_2,\dots$ is Cauchy; therefore, it is converging.
Observe that its limit $p_\infty$ belongs to each $\Omega_n$.
It follows that any closed ball $\cBall[p_0,R_0]$ contains a point in $\bigcap_{n\in \NN}\Omega_n$, hence the result.
\qeds

\section{Compact spaces}

Let us recall a few statements about compact metric spaces.

\begin{thm}{Definition}\label{def:compact}
A metric space $\spc{K}$ is compact if and only if one of the following equivalent conditions holds:

\begin{subthm}{}
 Every open cover of $\spc{K}$ has a finite subcover.
\end{subthm}

\begin{subthm}{}
 Every sequence of points in $\spc{K}$ has a subsequence that converges in $\spc{K}$.
\end{subthm}

\begin{subthm}{totally-bounded}
The space $\spc{K}$ is complete and \index{totally bounded space}\emph{totally bounded}; that is, for any $\eps>0$, the space $\spc{K}$ admits a finite cover by open $\eps$-balls.
\end{subthm}

\end{thm}

\begin{thm}{Lebesgue lemma}
Let $\spc{K}$ be a compact metric space.
Then for any open cover of $\spc{K}$, there is $\eps>0$ such that any $\eps$-ball in $\spc{K}$ lies in an element of the cover.

The value $\eps$ is called a \index{Lebesgue number}\emph{Lebesgue number} of the covering.
\end{thm}

A subset $N$ of a metric space $\spc{K}$ is called \index{net}\emph{$\eps$-net} if any point $x\in \spc{K}$ lies at the distance less than $\eps$ from a point in $N$.
More generally, a subset $N$ is called an \index{net}\emph{$\eps$-net} of a subset $S\subset\spc{K}$ if any point $x\in S$ lies at the distance less than $\eps$ from a point in $N$. 

Note that totally bounded spaces can be defined as spaces that admit a finite $\eps$-net for any $\eps>0$.

\begin{thm}{Exercise}\label{ex:compact-net}
Show that a space $\spc{K}$ is totally bounded if and only if it contains a compact $\eps$-net for any $\eps>0$. 
\end{thm}

Let $\pack_\eps\spc{X}$ be the exact upper bound on the number of points $x_1,\z\dots,x_n\in \spc{X}$ such that $\dist{x_i}{x_j}{}\ge\eps$ if $i\ne j$.

If $n=\pack_\eps\spc{X}<\infty$, then
the collection of points $x_1,\dots,x_n$ is called a \index{maximal packing}\emph{maximal $\eps$-packing}.
If $\spc{X}$ is a length space (see Section~\ref{sec:intrinsic}), then $n$ is the maximal number of disjoint open $\tfrac\eps2$-balls in $\spc{X}$.

\begin{thm}{Exercise}\label{ex:pack-net}
Show that any maximal $\eps$-packing is an $\eps$-net.
Conclude that a complete space $\spc{X}$ is compact if and only if $\pack_\eps\spc{X}\z<\infty$ for any $\eps>0$.
\end{thm}

\begin{thm}{Exercise}\label{ex:non-contracting-map}
Let $\spc{K}$  be a compact metric space and
$f\:\spc{K}\z\to \spc{K}$
be a distance-noncontracting map.
Prove that $f$ is an \index{isometry}\emph{isometry};
that is, $f$ is a distance-preserving bijection.
\end{thm}

A metric space $\spc{X}$ is called \index{locally compact space}\emph{locally compact} if any point in $\spc{X}$ admits a compact neighborhood;
equivalently, for any point $x\in\spc{X}$, a closed ball $\cBall[x,r]$ is compact for some $r>0$.

\section{Proper spaces}

A metric space $\spc{X}$ is called \index{proper space}\emph{proper} if all closed bounded sets in $\spc{X}$ are compact.
Note that $\spc{X}$ is proper if for some (and therefore any) point $p\in \spc{X}$ and any $R<\infty$, 
the closed ball $\cBall[p,R]_{\spc{X}}$ is compact.

Recall that a function $f\:\spc{X}\to\RR$ is \index{proper function}\emph{proper} if, for any compact set $K\subset \RR$, its inverse image $f^{-1}(K)$ is compact.
Observe that $\spc{X}$ is proper if and only if the function $\distfun_p\:\spc{X}\to\RR$ is \index{proper function}\emph{proper} for some (and therefore any) point $p\in \spc{X}$.

\begin{thm}{Exercise}\label{ex:loc-compact-not-proper}
Give an example of a metric space that is locally compact but not proper.
\end{thm}

\section{Geodesics}
\label{sec:geods}

Let $\spc{X}$ be a metric space 
and $\II$\index{$\II$} a real interval. 
A distance-preserving map $\gamma\:\II\to \spc{X}$ is called a \index{geodesic}\emph{geodesic}%
\footnote{Others call it differently: \textit{shortest path}, \textit{minimizing geodesic}.
Also, note that the meaning of the term \textit{geodesic} is different from what is used in Riemannian geometry, altho they are closely related.}; 
in other words, $\gamma\:\II\to \spc{X}$ is a geodesic if 
\[\dist{\gamma(s)}{\gamma(t)}{\spc{X}}=|s-t|\]
for any pair $s,t\in \II$.

If $\gamma\:[a,b]\to \spc{X}$ is a geodesic such that $p=\gamma(a)$, $q=\gamma(b)$, then we say that $\gamma$ is a geodesic from $p$ to $q$.
In this case, the image of $\gamma$ is denoted by $[p q]$\index{$[{*}{*}]$}, and, with abuse of notations, we also call it a \index{geodesic}\emph{geodesic}.
We may write $[p q]_{\spc{X}}$ 
to emphasize that the geodesic $[p q]$ is in the space  ${\spc{X}}$.

In general, a geodesic from $p$ to $q$ need not exist and if it exists, it need not  be unique.  
However, once we write $[p q]$ we assume that we have chosen such geodesic.

A \index{geodesic path}\emph{geodesic path} is a geodesic with constant-speed parameterization by the unit interval $[0,1]$.

A metric space is called \index{geodesic space}\emph{geodesic} if any pair of its points can be joined by a geodesic.

An $\infty$-metric space $\spc{X}$ is called {}\emph{geodesic} if each metric component of $\spc{X}$ is geodesic.

\begin{thm}{Exercise}\label{ex:pogorelov}
Let $f$ be a centrally symmetric positive continuous function on $\mathbb{S}^2$.
Given two points $x,y\in \mathbb{S}^2$,
set 
\[\|x-y\|=\int_{B(x,\frac \pi2)\setminus B(y,\frac\pi2)}f.\]

Show that $(\mathbb{S}^2,\|{*}-{*}\|)$ is a geodesic space,
and the geodesics in $(\mathbb{S}^2,\|{*}-{*}\|)$ run along great circles of $\mathbb{S}^2$.
\end{thm}

\section{Metric trees}

A geodesic space $\spc{T}$ is called a \index{metric tree}\emph{metric tree} if any two points in $\spc{T}$ are connected by a unique geodesic,
and the union of any two geodesics $[xy]_{\spc{T}}$, and $[yz]_{\spc{T}}$ contains the geodesic $[xz]_{\spc{T}}$.

{

\begin{wrapfigure}{r}{25 mm}
\vskip-6mm
\centering
\includegraphics{mppics/pic-105}
\end{wrapfigure}

The latter means that any triangle in $\spc{T}$ is a \index{tripod}\emph{tripod};
that is, for any three points $x$, $y$, and $z$ there is a point $p$ such that 
\[[xy]\cup[yz]\cup[zx]=[px]\cup[py]\cup[pz].\]

}

\begin{thm}{Exercise}\label{ex:4-point-trees}
Let $p$, $x$, $y$, and $z$ be points in a metric tree.

\begin{subthm}{ex:4-point-trees:diagonals}
Consider three numbers 
\begin{align*}
a&=|p-x|+|y-z|,
\\
b&=|p-y|+|z-x|,
\\
c&=|p-z|+|x-y|.
\end{align*}
Suppose that $a\le b\le c$.
Show that $b=c$.
\end{subthm}

\begin{subthm}{ex:4-point-trees:gromov-product}
Consider three numbers 
\begin{align*}
\alpha&=\tfrac12\cdot(|p-y|+|p-z|-|y-z|),
\\
\beta &=\tfrac12\cdot(|p-x|+|p-z|-|x-z|),
\\
\gamma&=\tfrac12\cdot(|p-x|+|p-y|-|x-y|).
\end{align*}
Suppose that $\alpha\le \beta\le \gamma$.
Show that $\alpha=\beta$.
\end{subthm}

\end{thm}

The set 
\[S(p,r)_{\spc{X}}=\set{x\in\spc{X}}{\dist{p}{x}{\spc{X}}=r}\]
will be called a \index{sphere}\emph{sphere} with center $p$ and radius $r$ in a metric space~$\spc{X}$.

\begin{thm}{Exercise}\label{ex:spheres-in-trees}
Show that spheres in metric trees are ultrametric spaces.
That is, 
\[\dist{x}{z}{}
\le
\max\{\,\dist{x}{y}{},\dist{y}{z}{}\,\}\]
for any $x,y,z\in S(p,r)_{\spc{T}}$.
\end{thm}

\section{Length}

A \index{curve}\emph{curve} is defined as a continuous map from a real interval $\II$ to a metric space.
If $\II=[0,1]$, then the curve is called a \index{path}\emph{path}.

\begin{thm}{Definition}
Let $\spc{X}$ be a metric space and
$\alpha\: \II\to \spc{X}$ be a curve.
We define the \index{length}\emph{length} of $\alpha$ as 
\[
\length \alpha \df \sup_{t_0\le t_1\le\ldots\le t_n}\sum_i \dist{\alpha(t_i)}{\alpha(t_{i-1})}{}.
\]

A curve $\alpha$ is called \index{rectifiable curve}\emph{rectifiable} if $\length \alpha<\infty$.
\end{thm}

\begin{thm}{Theorem}\label{thm:length-semicont}
Length is a lower semi-continuous with respect to the pointwise convergence of curves. 

More precisely, assume that a sequence
of curves $\gamma_n\:\II\to \spc{X}$ in a metric space $\spc{X}$ converges pointwise 
to a curve $\gamma_\infty\:\II\to \spc{X}$;
that is, for any fixed $t \in \II$ we have $\gamma_n(t)\z\to\gamma_\infty(t)$ as $n\to\infty$. 
Then 
$$\liminf_{n\to\infty} \length\gamma_n \ge \length\gamma_\infty.\eqlbl{eq:semicont-length}$$
\end{thm}

\begin{wrapfigure}{o}{20 mm}
\vskip-0mm
\centering
\includegraphics{mppics/pic-100}
\end{wrapfigure}

Note that the inequality \ref{eq:semicont-length} might be strict.
For example, the diagonal $\gamma_\infty$ of the unit square 
can be  approximated by stairs-like
polygonal curves $\gamma_n$
with sides parallel to the sides of the square ($\gamma_6$ is on the picture).
In this case
\[\length\gamma_\infty=\sqrt{2}\quad
\text{and}\quad \length\gamma_n=2\]
for any $n$.

\parit{Proof.}
Fix a sequence $t_0\le\dots\le t_k$ in $\II$.
Set 
\begin{align*}\Sigma_n
&\df
|\gamma_n(t_0)-\gamma_n(t_1)|+\dots+|\gamma_n(t_{k-1})-\gamma_n(t_k)|.
\\
\Sigma_\infty
&\df
|\gamma_\infty(t_0)-\gamma_\infty(t_1)|+\dots+|\gamma_\infty(t_{k-1})-\gamma_\infty(t_k)|.
\end{align*}

Note that for each $i$ we have 
\[|\gamma_n(t_{i-1})-\gamma_n(t_i)|\to|\gamma_\infty(t_{i-1})-\gamma_\infty(t_i)|\]
and therefore
\[\Sigma_n\to \Sigma_\infty\] 
as $n\to\infty$.
Note that 
\[\Sigma_n\le\length\gamma_n\]
for each $n$.
Hence,
$$\liminf_{n\to\infty} \length\gamma_n \ge \Sigma_\infty.$$

Since the partition was arbitrary, applying the definition of length, we get \ref{eq:semicont-length}.
\qeds

\begin{thm}{Exercise}\label{ex:1-Lip-G-delta}
Show that most of 1-Lipschitz paths in the plane have length 1.

More precisely, consider the space $\spc{P}$ of 1-Lipschitz paths in the plane;
that is, all paths $a\:[0,1]\to\RR^2$ such that $|a(t_0)-a(t_1)|\le |t_0-t_1|$ for any $t_0$ and $t_1$.
Equip $\spc{P}$ with the metric defined by
\[\dist{a}{b}{}\df\sup\set{\dist{a(t)}{b(t)}{}}{t\in[0,1]}.\]
Show that a dense G-delta set of paths in $\spc{P}$ has length 1.

\end{thm}

\section{Length spaces}\label{sec:intrinsic}

Let $\spc{X}$ be a metric space.
If for any $\eps>0$ and any pair of points $x,y\in\spc{X}$, there is a path $\alpha$ connecting $x$ to $y$ such that
\[\length\alpha< \dist{x}{y}{}+\eps,\]
then $\spc{X}$ is called a \index{length space}\emph{length space} and the metric on $\spc{X}$ is called a \index{length metric}\emph{length metric}.\label{page:length metric}

An $\infty$-metric space is a length space if each of its metric components is a length space.
In other words, if $\spc{X}$ is an $\infty$-metric space, then in the above definition we assume in addition that $\dist{x}{y}{\spc{X}}<\infty$.

Note that any geodesic space is a length space.
The following example shows that the converse does not hold.

\begin{thm}{Example}
Set $\II_n=[0,1+\tfrac1n]$ for every natural $n$.
Suppose a space $\spc{X}$ is obtained by gluing intervals $\{\II_n\}$, where the left ends are glued to $p$ and the right ends to~$q$.

Observe that the space $\spc{X}$ carries a natural complete length metric with respect to which $\dist{p}{q}{\spc{X}}=1$, but there is no geodesic connecting $p$ to~$q$.
\end{thm}

\begin{thm}{Exercise}\label{ex:no-geod}
Give an example of a complete length space $\spc{X}$ such that no pair of distinct points in $\spc{X}$ can be joined by a geodesic.
\end{thm}

Directly from the definition, it follows that if $\alpha\:[0,1]\to\spc{X}$ is a path from $x$ to $y$ 
(that is, $\alpha(0)=x$ and $\alpha(1)=y$), then 
\[\length\alpha\ge \dist{x}{y}{}.\]
Set 
\[\yetdist{x}{y}{}=\inf\{\,\length\alpha\,\}\]
where the greatest lower bound is taken for all paths from $x$ to~$y$.
It is straightforward to check that $(x,y)\mapsto \yetdist{x}{y}{}$ is an $\infty$-metric; moreover, $(\spc{X},\yetdist{*}{*}{})$ is a length space.
The metric $\yetdist{*}{*}{}$ is called the \index{induced length metric}\emph{induced length metric}.

\begin{thm}{Exercise}\label{ex:compact+connceted}
Let $\spc{X}$ be a complete length space.
Show that for any compact subset $K\subset\spc{X}$
there is a compact path-connected subset $K'\subset\spc{X}$ that contains $K$.  
\end{thm}

\begin{thm}{Exercise}\label{ex:compact=>complete}
Suppose $(\spc{X},\dist{*}{*}{})$ is a complete metric space.
Show that $(\spc{X},\yetdist{*}{*}{})$ is complete.
\end{thm}

Let $A$ be a subset of a metric space $\spc{X}$.
Given two points $x,y\in A$,
consider the value
\[\dist{x}{y}{A}=\inf_{\alpha}\{\,\length\alpha\,\},\]
where the greatest lower bound is taken for all paths $\alpha$ from $x$ to $y$ in~$A$.
In other words, $\dist{*}{*}{A}$ denotes the induced length metric on the subspace $A$.
(The notation $\dist{*}{*}{A}$ conflicts with the previously defined notation for distance $\dist{x}{y}{\spc{X}}$ in a metric space $\spc{X}$.
However, most of the time we will work with ambient length spaces where the meaning will be unambiguous.)

Let $x$ and $y$ be points in a metric space $\spc{X}$.

\begin{enumerate}[(i)]
\item A point $z\in \spc{X}$ is called a \index{midpoint}\emph{midpoint} between $x$ and $y$
if 
\[\dist{x}{z}{}=\dist{y}{z}{}=\tfrac12\cdot\dist[{{}}]{x}{y}{}.\]
\item Assume $\eps\ge 0$.
A point $z\in \spc{X}$ is called an \index{$\eps$-midpoint}\emph{$\eps$-midpoint} between $x$ and $y$
if 
\[\dist{x}{z}{}\le\tfrac12\cdot\dist[{{}}]{x}{y}{}+\eps
\quad\text{and}\quad
\dist{y}{z}{}\le\tfrac12\cdot\dist[{{}}]{x}{y}{}+\eps.\]
\end{enumerate}

Note that a $0$-midpoint is the same as a midpoint.

\begin{thm}{Menger's lemma}\label{lem:mid>geod}
Assume $\spc{X}$ is a complete metric space.
\begin{subthm}{lem:mid>length}
Suppose that for any two points in $\spc{X}$,  
and any positive $\eps$,
there is an $\eps$-midpoint.
Then $\spc{X}$ is a length space.
\end{subthm}

\begin{subthm}{lem:mid>geod:geod}
Suppose that for any two points in $\spc{X}$, 
there is a midpoint.
Then $\spc{X}$ is a geodesic space.
\end{subthm}
\end{thm}

The second part of this lemma was proved by Karl Menger \cite[Section 6]{menger}.

\parit{Proof; \ref{SHORT.lem:mid>length}.}
Choose $x,y\in \spc{X}$;
set $\eps_n=\frac\eps{4^n}$, $\alpha(0)=x$, and $\alpha(1)=y$.

Let $\alpha(\tfrac12)$ be an $\eps_1$-midpoint between $\alpha(0)$ and $\alpha(1)$.
Further, let $\alpha(\frac14)$ 
and $\alpha(\frac34)$ be $\eps_2$-midpoints between the pairs $(\alpha(0),\alpha(\tfrac12))$ 
and $(\alpha(\tfrac12),\alpha(1))$ respectively.
Continue the above procedure;
on the $n$-th step, we define $\alpha(\tfrac{k}{2^n})$,
for every odd integer $k$ such that $0<\tfrac k{2^n}<1$, 
as an $\eps_{n}$-midpoint of the already defined
$\alpha(\tfrac{k-1}{2^n})$ and $\alpha(\tfrac{k+1}{2^n})$.

This way we define $\alpha(t)$ for all dyadic rationals $t$ in $[0,1]$.
Moreover, $\alpha$ has Lipschitz constant $\dist{x}{y}{}+\eps$.
Since $\spc{X}$ is complete, the map $\alpha$ can be extended $(\dist{x}{y}{}+\eps)$-Lipschitz map $\alpha\:[0,1]\to \spc{X}$.
In particular
\[\begin{aligned}
\length\alpha&\le \dist{x}{y}{}+\eps.
\end{aligned}
\eqlbl{eq:eps-midpoint}
\]
Since $\eps>0$ is arbitrary, we get \ref{SHORT.lem:mid>length}.

\parit{\ref{SHORT.lem:mid>geod:geod}.} 
Apply the same argument 
with midpoints instead of $\eps_n$-midpoints.
In this case, \ref{eq:eps-midpoint} holds for $\eps_n=\eps\z=0$.
\qeds

In a compact space, a sequence of $\tfrac1n$-midpoints $z_n$ contains a convergent subsequence.
Therefore Menger's lemma (\ref{lem:mid>geod}) implies the following.

\begin{thm}{Proposition}\label{prop:length+proper=>geodesic}
Any proper length space is geodesic.
\end{thm}

\begin{thm}{Hopf--Rinow theorem}\label{thm:Hopf-Rinow}
Any complete, locally compact length space is proper.
\end{thm}

Before reading the proof, it is instructive to solve \ref{ex:loc-compact-not-proper}.
In the proof, we will use the following exercise.

\begin{thm}{Exercise}\label{ex:eps-nbhd(ball)}
Let $\spc{X}$ be a length space.
Show that $\oBall(x,R+\eps)_{\spc{X}}$ is the $\eps$-neighborhood of $\oBall(x,R)_{\spc{X}}$.
\end{thm}

\parit{Proof.}
Choose a point $x$ in a locally compact length space $\spc{X}$.
Let 
\[\rho(x)
\df
\sup\set{R}{\cBall[x,R]\ \text{is compact}}.\]
Since $\spc{X}$ is locally compact, 
$$\rho(x)>0
\quad\text{for any}\quad
x\in \spc{X}.\eqlbl{eq:rho>0}$$
It is sufficient to show that $\rho(x)=\infty$ for some (and therefore any) point $x\in \spc{X}$.

\begin{clm}{} If $\rho(x)<\infty$, then $B=\cBall[x,\rho(x)]$ is compact.
\end{clm}

Suppose  $\rho(x)>\eps>0$;
by \ref{ex:eps-nbhd(ball)}, 
the set $\cBall[x,\rho(x)-\eps]$ is a compact $2\cdot\eps$-net in~$B$.
Since $B$ is closed and hence complete, it must be compact; see \ref{totally-bounded} and \ref{ex:compact-net}.
\claimqeds

\begin{clm}{} $|\rho(x)-\rho(y)|\le \dist{x}{y}{\spc{X}}$ for any $x,y\in \spc{X}$;
in particular, $\rho\:\spc{X}\to\RR$ is a continuous function.
\end{clm}

Suppose $\rho(x)+|x-y|<\rho(y)$ for some $x,y\in \spc{X}$. 
Then 
$\cBall[x,\rho(x)+\eps]$ is a closed subset of $\cBall[y,\rho(y)]$ for some $\eps>0$.
Since $\cBall[y,\rho(y)]$ is compact, so is $\cBall[x,\rho(x)+\eps]$ --- a contradiction.\claimqeds

Let $\eps=\min\set{\rho(y)}{y\in B}$; the minimum is defined since $B$ is compact and $\rho$ is continuous.
By \ref{eq:rho>0}, we have $\eps>0$.

Choose a finite $\tfrac\eps{10}$-net $\{a_1,a_2,\dots,a_n\}$ in $B=\cBall[x,\rho(x)]$.
The union $W$ of the closed balls $\cBall[a_i,\eps]$ is compact.
By \ref{ex:eps-nbhd(ball)},
$\cBall[x,\rho(x)+\frac\eps{10}]\subset W$.
Therefore, $\cBall[x,\rho(x)+\frac\eps{10}]$ is compact,
a contradiction.
\qeds

\begin{thm}{Exercise}\label{exercise from BH}
Construct a geodesic space $\spc{X}$ that is locally compact,
but whose completion $\bar{\spc{X}}$ is neither geodesic nor locally compact.
\end{thm}

\begin{thm}{Advanced exercise}\label{ex:gross}
Show that for any compact connected space $\spc{X}$ there is a number $\ell$ such that for any finite collection of points there is a point $z$ that lies on average distance $\ell$ from the collection;
that is, for any $x_1,\dots,x_n\in \spc{X}$ there is $z\in \spc{X}$ such that
\[\tfrac1n\cdot\sum_i|x_i-z|_{\spc{X}}=\ell.\]
\end{thm}

%% file: uryson.tex
\chapter{Universal spaces}\label{chap:urysohn}

The Urysohn space is the main hero of this lecture.
It shares some fundamental properties with classical spaces (spheres, euclidean, and Lobachevsky spaces),
but also has many counterintuitive properties.

This space often serves as a counterexample to plausible conjectures;
so it is worth knowing it.
In addition, this space is beautiful.

\section{Embedding in a normed space}

Recall that a function $v\mapsto |v|$ on a vector space $\spc{V}$ is called \index{norm}\emph{norm} if it satisfies the following condition for any two vectors $v,w\in \spc{V}$ and a scalar $\alpha$:
\begin{itemize}
\item $|v|\ge 0$;
\item $|\alpha\cdot v|=|\alpha|\cdot |v|$;
\item $|v|+|w|\ge|v+w|$.
\end{itemize}

As an example, consider \index{$\ell^\infty$}$\ell^\infty$ --- the space of real sequences equipped with \index{sup-norm}\emph{sup-norm}; that is, the norm of $\bm{a}=(a_1,a_2,\dots)$ is defined by
\[|\bm{a}|_{\ell^\infty}
\df
\sup_n\{\,|a_n|\,\}.\]

It is straightforward to check that for any normed space the function $(v,w)\mapsto |v-w|$ defines a metric on it.
Therefore, any normed space is an example of metric space;
moreover, it is a geodesic space.
Often we do not distinguish normed space from the corresponding metric space.
(By the Mazur--Ulam theorem, the metric remembers the affine structure of the space; so, to recover the original normed space we only need to specify the origin.
A slick proof of this theorem was given by Jussi V\"{a}is\"{a}l\"{a} \cite{vaisala}.)

Recall that \index{diameter}\emph{diameter} of a metric space $\spc{X}$ (briefly $\diam \spc{X}$) is defined as the least upper bound on the distances between pairs of its points;
that is,
\[\diam \spc{X}
\df
\sup\set{\dist{x}{y}{\spc{X}}}{x,y\in \spc{X}}.\]
If $\diam\spc{X}<\infty$, then the space $\spc{X}$ is called \index{bounded space}\emph{bounded}.

\begin{thm}{Lemma}\label{lem:frechet}
Suppose $\spc{X}$ is a bounded \index{separable space}\emph{separable} metric space;
that is, $\spc{X}$ contains a countable, dense set, say $\{w_n\}$.
Given $x\in \spc{X}$, set $a_n(x)=\dist{w_n}{x}{\spc{X}}$.
Then 
\[\iota\:x\mapsto (a_1(x), a_2(x),\dots)\]
defines a distance-preserving embedding $\iota\:\spc{X}\hookrightarrow \ell^\infty$.
\end{thm}

\parit{Proof.} 
By the triangle inequality 
\[|a_n(x)-a_n(y)|\le \dist{x}{y}{\spc{X}}.\eqlbl{eq:a-a=<dist}\]
Therefore, $\iota$ is \index{short map}\emph{short} (in other words, $\iota$ is distance-nonexpanding).

Again by the triangle inequality we have 
\[|a_n(x)-a_n(y)|\ge \dist{x}{y}{\spc{X}}-2\cdot\dist{w_n}{x}{\spc{X}}.\]
Since the set $\{w_n\}$ is dense, we can choose $w_n$ arbitrarily close to $x$.
Whence 
\[\sup_n\{\,|a_n(x)-a_n(y)|\,\}\ge \dist{x}{y}{\spc{X}};\eqlbl{eq:a-a>=dist}\]
that is, $\iota$ is distance-noncontracting.

Finally, observe that \ref{eq:a-a=<dist} and \ref{eq:a-a>=dist} imply the lemma.
\qeds

\begin{thm}{Exercise}\label{ex:compact-length}
Show that any compact metric space $\spc{K}$ is isometric to a subspace of a compact geodesic space. 
\end{thm}

The following exercise generalizes the lemma to arbitrary separable spaces.

\begin{thm}{Exercise}\label{ex:frechet}
Suppose $\{w_n\}$ is a countable, dense set in a metric space $\spc{X}$.
Choose $x_0\in \spc{X}$;
given $x\in \spc{X}$, set 
\[a_n(x)=\dist{w_n}{x}{\spc{X}}-\dist{w_n}{x_0}{\spc{X}}.\]
Show that $\iota\:x\mapsto (a_1(x), a_2(x),\dots)$ defines a distance-preserving embedding $\iota\:\spc{X}\hookrightarrow \ell^\infty$.

Conclude that any separable metric space $\spc{X}$ admits a distance-preserving embedding $\iota\:\spc{X}\hookrightarrow \ell^\infty$.
\end{thm}

The following lemma implies that {}\textit{any metric space is isometric to a subset of a normed vector space};
its proof is nearly identical to the proof of \ref{ex:frechet}.
Given a set $\spc{X}$, denote by \index{$\ell^\infty(\spc{X})$}$\ell^\infty(\spc{X})$ the space of all bounded functions on $\spc{X}$ equipped with sup-norm; 
that is,
\[|f-g|_{\ell^\infty}=\sup\set{|f(x)-f(x)}{x\in \spc{X}}.\]

\begin{thm}{Lemma}\label{lem:kuratowski}
Let $x_0$ be a point in a metric space $\spc{X}$.
Then the map $\iota\:\spc{X}\to \ell^\infty(\spc{X})$ defined by 
\[\iota\:x\mapsto (\distfun_x-\distfun_{x_0})\]
is distance-preserving.

In particular, any metric space $\spc{X}$ admits a distance-preserving into~$\ell^\infty(\spc{X})$.
\end{thm}

\section{Extension property}
\label{sec:Extension property}

If a metric space $\spc{X}$ is a subspace of a semimetric space $\spc{X}'$, then we say that $\spc{X}'$ is an \index{extension}\emph{extension} of $\spc{X}$.
If in addition, $\diam\spc{X}'\le d$, then we say that $\spc{X}'$ is a {}\emph{$d$-extension}.

If the complement $\spc{X}'\setminus \spc{X}$ contains a single point, say $p$, then $\spc{X}'$ is called a \index{one-point extension}\emph{one-point extension} of $\spc{X}$.
In this case, to define a metric on $\spc{X}'$, it is sufficient to specify the distance function from $p$; that is, a function $f\:\spc{X}\to\RR$ defined by 
\[f(x)\df\dist{p}{x}{\spc{X}'}.\]
Any function $f$ of that type will be called an \index{extension function}\emph{extension function}\label{page:extension function} or {}\emph{$d$-extension function} respectively.

The extension function $f$ cannot be taken arbitrarily --- the triangle inequality implies that 
\[f(x)+f(y)\ge \dist{x}{y}{\spc{X}}\ge |f(x)-f(y)|\]
for any $x,y\in \spc{X}$.
In particular, $f$ is a non-negative 1-Lipschitz function on $\spc{X}$.
For a $d$-extension, we need to assume in addition that $\diam\spc{X}\z\le d$ and $f(x)\le d$ for any $x\in \spc{X}$.
A straightforward check shows that these conditions are necessary and sufficient.

\begin{thm}{Exercise}\label{ex:extension-of-extension}
Let $\spc{X}$ be a subspace of metric space $\spc{Y}$.
Assume $f$ is an extension function on $\spc{X}$.

\begin{subthm}{ex:extension-of-extension:a}
Show that 
\[\bar f(y)
\df
\inf_{x\in \spc{X}} \{\,f(x)+\dist{x}{y}{\spc{Y}}\,\}\]
defines an extension function on $\spc{Y}$.
\end{subthm}

\begin{subthm}{}
Assume that $\diam \spc{Y}\le d$ and $f(x)\le d$ for any $x\in  \spc{X}$.
Show that 
\[\bar f_d
\df
\min \{\, \bar f,d\,\}\]
is a $d$-extension function on $\spc{Y}$.
\end{subthm}

\end{thm}

The functions $\bar f$ and $\bar f_d$ in the above exercise are called \index{Katětov extensions}\emph{Katětov extensions} of $f$ and the minimal possible $\spc{X}$ is called its \index{support of extension function}\emph{support}, briefly \index{$\supp$}$\supp \bar f=\spc{X}$.

\begin{thm}{Definition}\label{def:finite+1}
A metric space $\spc{U}$ meets the \index{extension property}\emph{extension property}  if for any finite subspace $\spc{F}\subset\spc{U}$ and any extension function $f\:\spc{F}\to\RR$ there is a point $p\in \spc{U}$ such that $\dist{p}{x}{}=f(x)$ for any $x\in \spc{F}$.

If we assume in addition that $\diam \spc{U}\le d$ and instead of extension functions we consider only $d$-extension functions, then it defines the {}\emph{$d$-extension property}.

Furhter, if in addition, $\spc{U}$ is separable and complete, then it is called \index{Urysohn space}\emph{Urysohn space} or {}\emph{$d$-Urysohn space} respectively.
\end{thm}

\begin{thm}{Proposition}\label{prop:univeral-separable}
There is a separable metric space with the ($d$-) extension property (for any $d\ge 0$).
\end{thm}

\parit{Proof.}
Choose $d\ge 0$.
Let us construct a separable metric space with  the $d$-extension property.

Let $\spc{X}$ be a metric space such that $\diam\spc{X}\le d$.
Denote by $\spc{X}^d$ the space of all $d$-extension functions on $\spc{X}$ equipped with the metric defined by the sup-norm.
Note that the map $\spc{X} \to \spc{X}^d$ defined by $x\mapsto\distfun_x$ is a distance-preserving embedding,
so we can (and will) treat $\spc{X}$ as a subspace of $\spc{X}^d$; equivalently, $\spc{X}^d$ is an extension of $\spc{X}$.

Let us iterate this construction.
Start with a one-point space $\spc{X}_0$ and consider a sequence of spaces $(\spc{X}_n)$ defined by $\spc{X}_{n+1}\z\df\spc{X}_n^d$.
Note that the sequence is nested;
that is, $\spc{X}_0\subset \spc{X}_1\subset\dots$
and the union
\[\spc{X}_\infty=\bigcup_n\spc{X}_n;\]
comes with metric such that
$\dist{x}{y}{\spc{X}_\infty} = \dist{x}{y}{\spc{X}_n}$
if $x,y\in\spc{X}_n$.

Note that if $\spc{X}$ is compact, then so is $\spc{X}^d$.
It follows that each space $\spc{X}_n$ is compact.
In particular, $\spc{X}_\infty$ is a countable union of compact spaces;
therefore $\spc{X}_\infty$ is separable.

Any finite subspace $\spc{F}$ of $\spc{X}_\infty$ lies in some $\spc{X}_n$ for $n<\infty$.
By construction, given an extension function $f\:\spc{F}\to\RR$,
there is a point $p\in \spc{X}_{n+1}$ that meets the condition in \ref{def:finite+1}.
That is, $\spc{X}_\infty$ has the $d$-extension property.

The construction of a separable metric space with the extension property requires only two changes.
First, the sequence should be defined by $\spc{X}_{n+1}\z\df\spc{X}_n^{d_n}$, where $d_n$ is an increasing sequence such that $d_n\to\infty$.
Second, the point $p$ should be taken in $\spc{X}_{n+k}$ for sufficiently large $k$, so that $d_{n+k}>\max\{f(x)\}$
(here one has to apply \ref{ex:extension-of-extension:a}).%

(Alternatively, one can start with any separable space $\spc{X}_0$ and consider a nested sequence $\spc{X}_0\subset \spc{X}_1\subset{}\dots$ where $\spc{X}_{n+1}$ is the space of all extension functions on $\spc{X}_{n}$ with at most $n+1$ points in its support.
The last condition is needed to keep $\spc{X}_{n}$ separable.)
\qeds

Given a metric space $\spc{X}$, denote by $\spc{X}^\infty$ the space of all extension functions on $\spc{X}$ equipped with the metric defined by the sup-norm.

\begin{thm}{Exercise}\label{ex:inf-extension}
Construct a proper length space $\spc{X}$ such that $\spc{X}^\infty$ is not separable.
\end{thm}

\begin{thm}{Proposition}\label{prop:completion-univeral}
If a metric space $\spc{V}$ meets the ($d$-) extension property, then so does its completion.
\end{thm}

\parit{Proof.} 
Let us assume $\spc{V}$ meets the extension property.
We will show that its completion $\spc{U}=\bar{\spc{V}}$ meets the extension property as well.
The $d$-extension case can be proved along the same lines.

Note that $\spc{V}$ is a dense subset in a complete space $\spc{U}$.
Observe that $\spc{U}$ has the \index{approximate extension property}\emph{approximate extension property};
that is, if $\spc{F}\z\subset\spc{U}$ is a finite set, $\eps>0$, and $f\:\spc{F}\to \RR$ is an extension function, then
there exists $p\in \spc{U}$ such that
\[\dist{p}{x}{}\lg f(x)\pm\eps\eqlbl{eq:|p-x|><f(x)}\]
for any $x\in\spc{F}$.
Indeed, consider the Katětov extension $\bar f\:\spc{U}\to\RR$ of~$f$.
Since $\spc{V}$ is dense in $\spc{U}$, we can choose a finite set $\spc{F}'\in \spc{V}$ such that for any $x\in \spc{F}$ there is $x'\in \spc{F}'$ with $\dist{x}{x'}{}<\tfrac\eps2$.
Let $p$ be the point provided by the extension property for the restriction $\bar f|_{\spc{F}'}$.
It remains to observe $p$ meets \ref{eq:|p-x|><f(x)}.

It follows that there is a sequence of points $p_n\in \spc{U}$ such that for any $x\in \spc{F}$, 
\[\dist{p_n}{x}{}\lg f(x)\pm\tfrac1{2^n}.\]

Moreover, we can assume that 
\[\dist{p_n}{p_{n+1}}{} < \tfrac1{2^n}\eqlbl{eq:|pn-pn|}\]
for all large $n$.
Indeed, consider the sets $\spc{F}_n=\spc{F}\cup\{p_n\}$ and the functions $f_n\:\spc{F}_n\to\RR$ defined by $f_n(x)\df f(x)$ and
\[f_n(p_n)
\df
\max\set{\bigl|\dist{p_n}{x}{}- f(x)\bigr|}{x\in \spc{F}}\]
 if $x\ne p_n$.
Observe that $f_n$ is an extension function for large $n$ and
$f_n(p_n)\z<\tfrac1{2^n}$.
Therefore, applying the approximate extension property recursively we get~\ref{eq:|pn-pn|}.

Therefore, the sequence $p_n$ is Cauchy.
Note that its limit meets the condition in the definition of extension property (\ref{def:finite+1}).
\qeds

Note that \ref{prop:univeral-separable} and \ref{prop:completion-univeral} imply the following:

\begin{thm}{Theorem}\label{thm:urysohn-exists}
Urysohn space and $d$-Urysohn space exist for any $d>0$.
\end{thm}

Here is a slightly stronger statement:

\begin{thm}{Theorem}\label{thm:urysohn-exists+}
Any separable metric space $\spc{X}$ admits a distance-preserving embedding into an Urysohn space $\spc{U}$ such that any isometry of $\spc{X}$ can be extended to an isometry of $\spc{U}$.
\end{thm}

\parit{Sketch of proof.}
Start with $\spc{X}_0=\spc{X}$ and construct a nested sequence of spaces $\spc{X}_0\subset\spc{X}_1 \subset{}\dots$ as at the alternative end of the proof of~\ref{prop:univeral-separable}.
Note that 
any isometry $\spc{X}_n\to \spc{X}_n$ can be extended to a unique isometry $\spc{X}_{n+1}\to \spc{X}_{n+1}$.
It follows that any isometry of $\spc{X}$ can be extended to an isometry of $\spc{X}'=\bigcup_n\spc{X}_n$.

Now, consider a new nested sequence $\spc{X}\subset \spc{X}'\subset \spc{X}''\subset \dots$;
denote its union by $\spc{Y}$.
Arguing as in \ref{prop:univeral-separable} and \ref{prop:completion-univeral} we get that the completion of $\spc{Y}$ is an Urysohn space, say $\spc{U}$, that comes with a distance-preserving inclusion $\spc{X}\hookrightarrow \spc{U}$.

From above, any isometry of $\spc{X}$ can be extended to isometries of $\spc{X}'$, $\spc{X}''$, and so on.
They all define an isometry of $\spc{Y}$;
passing to its continuous extension, we get an isometry of $\spc{U}$.
\qeds

\section{Universality}

A metric space will be called \index{universal space}\emph{universal} if it has a subspace isometric to any given separable metric space.
In \ref{ex:frechet}, we proved that $\ell^\infty$ is a universal space. 
The following proposition shows that an Urysohn space is universal as well.
Unlike $\ell^\infty$, Urysohn spaces are separable;
so it might be considered as a \textit{better} universal space.
Theorem \ref{thm:compact-homogeneous} will give another reason why Urysohn spaces are better.

\begin{thm}{Proposition}\label{prop:sep-in-urys}
An Urysohn space is universal.
That is, if $\spc{U}$ is an Urysohn space, then any separable metric space $\spc{S}$ admits a distance-preserving embedding $\spc{S}\hookrightarrow\spc{U}$.

Moreover, for any finite subspace $\spc{F}\subset \spc{S}$,
any distance-preserving embedding $\spc{F}\hookrightarrow \spc{U}$ can be extended to a distance-preserving embedding $\spc{S}\hookrightarrow\spc{U}$.

A $d$-Urysohn space is $d$-universal;
that is, the above statements hold provided that $\diam\spc{S}\le d$.  
\end{thm}

\parit{Proof.}
We will prove the second statement;
the first statement is its partial case for $\spc{F}=\emptyset$.

The required isometry will be denoted by $x\mapsto x'$.

Choose a dense sequence of points $s_1,s_2,\dotsc\in\spc{S}$.
We may assume that $\spc{F}=\{s_1,\dots,s_n\}$, so $s_i'\in \spc{U}$ are defined for $i\le n$.

The sequence $s_i'$ for $i>n$ can be defined recursively using the extension property in $\spc{U}$.
Namely, suppose that $s_1',\dots,s_{i-1}'$ are already defined.
Since $\spc{U}$ meets the extension property, there is a point $s_i'\in \spc{U}$ such that
\[\dist{s_i'}{s_j'}{\spc{U}}=\dist{s_i}{s_j}{\spc{S}}\]
for any $j<i$.

The constructed map $s_i\mapsto s_i'$ is distance-preserving.
Therefore it can be continuously extended to the whole $\spc{S}$.
It remains to observe that the constructed map $\spc{S}\hookrightarrow\spc{U}$ is distance-preserving.
\qeds

\begin{thm}{Exercise}\label{ex:geodesics-urysohn}
Show that any two distinct points in an Urysohn space can be joined by an infinite number of distinct geodesics.
\end{thm}

\begin{thm}{Exercise}\label{ex:compact-extension}
Modify the proofs of \ref{prop:completion-univeral} and \ref{prop:sep-in-urys} to prove the following theorem.
\end{thm}

\begin{thm}{Theorem}\label{thm:compact-extension}
Let $K$ be a compact set in a separable space $\spc{S}$.
Then any distance-preserving map from $K$ to an Urysohn space can be extended to 
a distance-preserving map of the whole $\spc{S}$.
\end{thm}

\begin{thm}{Exercise}\label{ex:sc-urysohn}
Show that ($d$-) Urysohn space is simply-connected.
\end{thm}

\section{Uniqueness and homogeneity}

\begin{thm}{Theorem}\label{thm:urysohn-unique}
Suppose $\spc{F}\subset \spc{U}$ and $\spc{F}'\subset \spc{U}'$ be finite isometric subspaces in a pair of ($d$-)Urysohn spaces $\spc{U}$ and $\spc{U}'$.
Then any isometry $\iota\:\spc{F}\leftrightarrow \spc{F}'$ can be extended to an isometry $\spc{U}\leftrightarrow \spc{U}'$.

In particular, ($d$-)Urysohn space is unique up to isometry.
\end{thm}

Note that \ref{prop:sep-in-urys} implies that there are distance-preserving maps $\spc{U}\z\to \spc{U}'$ and $\spc{U}'\to \spc{U}$.
The next exercise shows that it does not solely imply the existence of an isometry $\spc{U}\leftrightarrow \spc{U}'$.

\begin{thm}{Exercise}\label{ex:no-isom}
Construct two metric spaces $\spc{X}$ and $\spc{Y}$ such that 
there are distance-preserving maps $\spc{X}\to \spc{Y}$ and $\spc{Y}\to \spc{X}$, but no isometry $\spc{X}\leftrightarrow \spc{Y}$.
\end{thm}

The following construction uses the idea of \ref{prop:sep-in-urys}, but it is applied \index{back-and-forth}\emph{back-and-forth} to ensure that the obtained distance-preserving map is onto.

\parit{Proof.}
Choose dense sequences $a_1,a_2,\dots{}\in \spc{U}$ and $b'_1,b'_2,\dots{}\in \spc{U}'$.
We can assume that $\spc{F}=\{a_1,\dots,a_n\}$, $\spc{F}'=\{b_1',\dots,b_n'\}$ and $\iota(a_i)=b_i'$ for $i\le n$.

The required isometry $\spc{U}\leftrightarrow \spc{U}'$ will be denoted by $u \leftrightarrow u'$.
Set $a_i=b_i$ and $a'_i=b'_i$ if $i\le n$.

Let us define recursively $a_{n+1}',b_{n+1}, a_{n+2}', b_{n+2},\dots$ --- on the odd step we define the images of $a_{n+1},a_{n+2},\dots$ and on the even steps we define inverse images of $b'_{n+1},b'_{n+2},\dots$
The same argument as in the proof of \ref{prop:sep-in-urys} shows that we can construct two sequences $a_1',a_2',\dots{}\in \spc{U}'$ and $b_1,b_2,\dots\in \spc{U}$ such that
\begin{align*}
\dist{a_i}{a_j}{\spc{U}}&=\dist{a_i'}{a_j'}{\spc{U}'},
\\
\dist{a_i}{b_j}{\spc{U}}&=\dist{a_i'}{b_j'}{\spc{U}'},
\\
\dist{b_i}{b_j}{\spc{U}}&=\dist{b_i'}{b_j'}{\spc{U}'}
\end{align*}
for all $i$ and $j$.

It remains to observe that the constructed distance-preserving bijection defined by $a_i\leftrightarrow a_i'$ and $b_i\leftrightarrow b_i'$ extends
continuously to an isometry $\spc{U}\leftrightarrow \spc{U}'$. 
\qeds

Observe that \ref{thm:urysohn-unique} implies that the Urysohn space (as well as the $d$-Urysohn space) is \index{homogeneous}\emph{finite-set-homogeneous}; that is,
\begin{itemize}
 \item any distance-preserving map from a finite subset to the whole space can be extended to an isometry.
\end{itemize}

Recall that $S(p,r)_{\spc{X}}$ denotes the sphere of radius $r$ centered at $p$ in a metric space $\spc{X}$;
that is, 
$$S(p,r)_{\spc{X}}=\set{x\in \spc{X}}{\dist{p}{x}{\spc{X}}=r}.$$

\begin{thm}{Exercise}\label{ex:sphere-in-urysohn}
Choose $d\in [0,\infty]$.
Denote by $\spc{U}_d$ the $d$-Urysohn space,
so $\spc{U}_\infty$ is the Urysohn space.

\begin{subthm}{ex:sphere-in-urysohn:sphere}
Assume that $L=S(p,r)_{\spc{U}_d}\ne \emptyset$.
Show that $L$ is isometric to $\spc{U}_{\ell}$; find $\ell$ in terms of $r$ and $d$.
\end{subthm}

\begin{subthm}{ex:sphere-in-urysohn:midpoint}
Let $\ell=\dist{p}{q}{\spc{U}_d}$.
Show that the subset $M\subset\spc{U}_d$ of midpoints between $p$ and $q$ is isometric to $\spc{U}_\ell$.
\end{subthm}

\begin{subthm}{ex:sphere-in-urysohn:homogeneous}
Show that $\spc{U}_d$ is \emph{not} countable-set-homogeneous;
that is, there is a distance-preserving map from a countable subset of $\spc{U}_d$ to $\spc{U}_d$ that cannot be extended to an isometry of $\spc{U}_d$.
\end{subthm}

\end{thm}

In fact, the Urysohn space is compact-set-homogeneous; more precisely the following theorem holds.

\begin{thm}{Theorem}\label{thm:compact-homogeneous}
Let $K$ be a compact set in a ($d$-)Urysohn space~$\spc{U}$.
Then any distance-preserving map $K\to \spc{U}$ can be extended to an isometry of $\spc{U}$.
\end{thm}

A proof can be obtained by modifying the proofs of \ref{prop:completion-univeral} and \ref{thm:urysohn-unique}
the same way as it is done in \ref{ex:compact-extension}.

\begin{thm}{Exercise}\label{ex:shere}
Let $S$ be a unit sphere in the Urysohn space $\spc{U}$.
Show that for any two distinct points $x,y\in \spc{U}$ there is a point $z\in S$ such that 
$\dist{x}{z}{}\ne \dist{y}{z}{}$.

Conclude that two isometries of $\spc{U}$ coincide if they coincide on $S$.
\end{thm}

\begin{thm}{Exercise}\label{ex:ext(shere)}
Let $B$ be an open unit ball in the Urysohn space $\spc{U}$.
Show that $\spc{U}\setminus B$ is isometric to $\spc{U}$.

Use it to construct an isometry of a unit sphere $S$ in $\spc{U}$ that cannot be extended to an isometry of $\spc{U}$.
\end{thm}

\begin{thm}{Exercise}\label{ex:katetov}

\begin{subthm}{ex:katetov:inclusion}
Show that there is a distance-preserving inclusion of the Urysohn space $\iota\:\spc{U}\hookrightarrow \spc{U}$ 
such that $\spc{U}'=\iota(\spc{U})$ is nowhere dense in $\spc{U}$ and any isometry of $\spc{U}'$ 
can be extended to an isometry of the whole~$\spc{U}$.
\end{subthm}

\begin{subthm}{ex:katetov:sol}
Consider a nested sequence $\spc{U}_0\subset \spc{U}_1\subset\dots$ of Urysohn spaces 
with each inclusion $\spc{U}_n\hookrightarrow \spc{U}_{n+1}$ as in \ref{SHORT.ex:katetov:inclusion}.
Show that the union $\bigcup_n\spc{U}_n$ is a noncomplete finite-set-homogeneous metric space that meets the extension property.
\end{subthm}

\end{thm}

{\sloppy

\begin{thm}{Exercise}\label{ex:homogeneous}
Which of the following metric spaces are 
one-point-homogeneous, finite-set-homogeneous, compact-set-homogeneous, countable-set-homogeneous?

\begin{subthm}{ex:homogeneous:euclidean}
Euclidean plane,
\end{subthm}

\begin{subthm}{ex:homogeneous:hilbert}
 Hilbert space $\ell^2$,
\end{subthm}

\begin{subthm}{ex:homogeneous:ell-infty}
 $\ell^\infty$,
\end{subthm}

\begin{subthm}{ex:homogeneous:ell-1}
 \index{$\ell^1$}$\ell^1$ --- the space of all real absolutely converging series $\bm{a}\z=(a_1,a_2,\dots)$ with the norm $|\bm{a}|_{\ell^1}=\sum_i|a_i|$.
 
\end{subthm}
\end{thm}

\begin{thm}{Exercise}\label{ex:homogeneous-tree}
Show that any separable one-point-homogeneous metric tree is isometric to the real line $\RR$ or the one-point space.
\end{thm}

}

\section{Remarks}

The statement in \ref{ex:frechet} was proved by Maurice René Fréchet in the paper where he first defined metric spaces \cite{frechet};
its extension \ref{lem:kuratowski} was given by Kazimierz Kuratowski~\cite{kuratowski}.
Both maps $x\mapsto (\distfun_x-\distfun_{x_0})$ and $x\mapsto \distfun_x$ can be called \index{Kuratowski embedding}\emph{Kuratowski embedding}.

\medskip

Let us describe a closely related construction introduced by Mikhael Gromov \cite{gromov-1981-hyperbolic}.
Suppose $\spc{X}$ be a proper metric space.
Denote by $C(\spc{X},\RR)$ the space of continuous functions $\spc{X}\to \RR$
equipped with the \emph{compact-open topology};
that is, for any compact set $K\subset \spc{X}$ and any open set $U\subset \RR$
the set of all continuous functions $f\: \spc{X}\to \RR$ such that $f(K)\subset U$
is declared to be open.

Choose a point $x_0\in \spc{X}$.
Consider the map $F_{\spc{X}}\:\spc{X}\to C(\spc{X},\RR)$
defined by
\[x\mapsto f_x\df -\dist{x}{x_0}{}+\distfun_x.\]

\begin{thm}{Exercise}\label{ex:horobry}
Show that if $\spc{X}$ is a proper length space, then $F_{\spc{X}}$ is an embedding.
Construct a proper metric space $\spc{Y}$ such that $F_{\spc{Y}}$ 
is not an embedding.

\end{thm}

Denote by $\hat{\spc{X}}$ 
the closure of $F_{\spc{X}}(\spc{X})$ in $C(\spc{X},\RR)$;
observe that $\hat {\spc{X}}$ is compact.
If $F_{\spc{X}}$ is an embedding, 
then $\hat {\spc{X}}$ is a compactification of $\spc{X}$,
and it is called the \index{horo-compactification}\emph{horo-compactification}.
In this case, the complement 
\[\partial_\infty \spc{X}\z=\hat {\spc{X}}\setminus F_{\spc{X}}(\spc{X})\] 
is called the {}\emph{horo-absolute} of~$\spc{X}$.
A variation of this construction for nonproper spaces was considered by Anders Karlsson \cite{karlsson-2023}.

\medskip

The following two exercises show that in this respect $\ell^\infty$ is very different from $\ell^1$.
For more on the subject, see \cite{deza-laurent}.

Let $S$ be a subset of $X$.
We say that $S$ \index{separating set}\emph{separates} $x$ and $y$ if 
$x\in S$ and $y\notin S$ or $x\notin S$ and $y\in S$.
The \index{cut metric}\emph{cut metric} $\delta_S$ on $X$ is a semimetric such that $\delta_S(x, y) = 1$ if $x$ and $y$ are separated
by $S$ and otherwise $\delta_S(x, y) = 0$.

\begin{thm}{Exercise}\label{ex:cut}
Show that a finite metric space $\spc{F}$ admits a distance-preserving embedding into $\ell^1$ if and only if the metric of $\spc{F}$ can be written as a nonnegative linear combination%
\footnote{that is, linear combination with nonnegative coefficients.} of cut metrics on $\spc{F}$.
\end{thm}

Recall that the vertex set of any graph comes with the \index{path metric}\emph{path metric} ---
the distance between two vertices is the minimal number of edges in a path connecting them.

{

\begin{wrapfigure}{r}{15mm}
\vskip-0mm
\centering
\includegraphics{mppics/pic-210}
\end{wrapfigure}

\begin{thm}{Exercise}\label{ex:K23}
Use \ref{ex:cut} to show that the metric for complete bipartite graph $K_{2,3}$ (see the diagram) does not admit a distance-preserving embedding into $\ell^1$.
\end{thm}

}

The question about the existence of a separable universal space was posed by Maurice René Fréchet and answered by
Pavel Urysohn~\cite{urysohn}.
Exercise \ref{ex:katetov} answers a question posed by Pavel Urysohn \cite[§$2(6)$]{urysohn}.
It was solved by Miroslav Katětov \cite{katetov},
but long after that, it was again mentioned as an open problem \cite[p. 83]{gromov-2007}.

The idea of Urysohn's construction was reused in graph theory; it produces the so-called \index{Rado graph}\emph{Rado graph},
also known as \textit{Erd\H{o}s--Rényi graph} or \textit{random graph}; see \cite{cameron}.
In fact, the Urysohn space is the random metric space in \textit{certain sense} \cite{vershik}.

\textit{The ($d$-) Urysohn space is homeomorphic to the Hilbert space};
the latter was proved by Vladimir Uspenskij \cite{uspenskij} using the so-called Toruńczyk criterion.

The finite-set-homogeneous spaces include euclidean spaces, hyperbolic spaces, and spheres all with standard length metrics and arbitrary finite dimensions.
In fact, these are the only examples of locally compact three-point-homogeneous length spaces.
The latter was proved by Herbert Busemann \cite{busemann-1942}; it also follows from the more general result of Jacques Tits about two-point-homogeneous spaces \cite{tits}.
The same conclusion holds for complete all-set-homogeneous geodesic spaces with local uniqueness of geodesics;
it was proved by Garrett Birkhoff \cite{birkhoff}.
The answer might be the same for complete separable all-set-homogeneous length spaces.
Without the separability condition, we also get the so-called \emph{universal metric trees} with finite valence \cite{dyubina-polterovich}; no other examples seem to be known \cite{lebedeva-petrunin2211.09671}.  

{\sloppy

\begin{thm}{Exercise}\label{ex:RP-not}
Show that the real projective plane $\RP^2$ with the standard metric is two-point-homogeneous, but not three-point-homogeneous.
\end{thm}

}

\begin{thm}{Exercise}\label{ex:hom-cube}
Let $Q$ be the set of vertices on the $n$-dimensional cube;
assume $n$ is large.
Show that $Q$ is three-point-homogeneous, but not four-point-homogeneous.
\end{thm}

\begin{thm}{Question}
Are there examples of metric spaces that are $n$-point-homogeneous, but not $(n+1)$-point-homogeneous for large $n$? See \cite{petrunin-431426}.
\end{thm}

%% file: injective.tex
\chapter{Injective spaces}\label{chap:injective}

Injective hull is a useful construction that provides a canonical choice of a specially nice (injective) space that includes a given metric space. 
This construction is similar to the convex hull in euclidean space.
The following exercise gives a bridge from the latter to the former.

\begin{thm}{Advanced exercise}\label{ex:conv-short}
Show that $A\subset \RR^n$ is a closed convex set if and only if for any  $B\subset \RR^n$ any short map $B\to A$ can be extended to a short map $\RR^n\to A$.
\end{thm}

\section{Definition}

\begin{thm}{Definition}\label{def:injective}
A metric space $\spc{Y}$ is called \index{injective space}\emph{injective} if, for any metric space $\spc{X}$ and any of its subspace $\spc{A}$,
any short map $f\:\spc{A}\to \spc{Y}$ can be extended to a short map $F\:\spc{X}\to \spc{Y}$;
that is, $f=F|_{\spc{A}}$.
\end{thm}

\begin{thm}{Exercise}\label{ex:inj=complete-geodesic-contractible}
Show that any injective space is 
\begin{multicols}{3}

\begin{subthm}{ex:inj=complete-geodesic-contractible:complete}
complete,
\end{subthm}

\begin{subthm}{ex:inj=complete-geodesic-contractible:geodesic}
geodesic, and
\end{subthm}

\begin{subthm}{ex:inj=complete-geodesic-contractible:contractible}
contractible.
\end{subthm}

\end{multicols}

\end{thm}

\begin{thm}{Exercise}\label{ex:bicombing}
Show that for any injective space $\spc{Y}$ there is a map $m\:\spc{Y}\times\spc{Y}\to\spc{Y}$ (the \index{midpoint map}\emph{midpoint map}) such that the inequality
\[2\cdot \dist{p}{m(x,y)}{\spc{Y}}\le\dist{p}{x}{\spc{Y}}+\dist{p}{y}{\spc{Y}}\]
holds for any $p,x,y\in \spc{Y}$.
\end{thm}

\begin{thm}{Exercise}\label{ex:injective-spaces}
Show that the following spaces are injective:
\begin{subthm}{ex:injective-spaces:R}
the real line;
\end{subthm}

\begin{subthm}{ex:injective-spaces:tree}
complete metric tree;
\end{subthm}

\begin{subthm}{ex:injective-spaces:ell-infty}
The space $\ell^\infty(\spc{S})$ for any set $\spc{S}$ (defined in \ref{lem:kuratowski}).
In particular, the coordinate plane with the metric induced by the $\ell^\infty$-norm.
\end{subthm}

\end{thm}

\begin{thm}{Exercise}\label{ex:extr-ball}
Let $\spc{Y}$ be an injective space.

\begin{subthm}{ex:extr-ball:one}
Show that any closed ball in $\spc{Y}$ is injective.
\end{subthm}

\begin{subthm}{ex:extr-ball:many}
Show that the intersection of an arbitrary collection of closed balls in $\spc{Y}$ is injective.
\end{subthm}

\end{thm}

\begin{thm}{Advanced exercise}\label{ex:extr-fixed}
Let $\spc{Y}$ be a bounded injective space.
Show that any short map $s\:\spc{Y}\to\spc{Y}$ has a fixed point. 
\end{thm}

\section{Admissible and extremal functions}

Let $\spc{X}$ be a metric space.
A function $r\:\spc{X}\to(-\infty,\infty]$ is called \label{page:admissible function}\index{admissible function}\emph{admissible} if the following inequality
\[r(x)+r(y)\ge \dist{x}{y}{\spc{X}}\eqlbl{eq:admissible}\]
holds for any $x,y\in \spc{X}$.

\begin{thm}{Observation}\label{obs:admissible}

\begin{subthm}{obs:admissible:nonnegative}
Any admissible function is nonnegative.
\end{subthm}

\begin{subthm}{obs:admissible:balls}
If $\spc{X}$ is a geodesic space, then a function $r\:\spc{X}\to\RR$ is admissible if and only if 
\[\cBall[x,r(x)]\cap\cBall[y,r(y)]\ne \emptyset\]
for any $x,y\in \spc{X}$.
\end{subthm}
 
\end{thm}

\parit{Proof; \ref{SHORT.obs:admissible:nonnegative}.} Apply \ref{eq:admissible} for $x=y$.

\parit{\ref{SHORT.obs:admissible:balls}.} Apply the triangle inequality and the existence of a geodesic $[xy]$.
\qeds

A minimal admissible function will be called \label{page:extremal function}\index{extremal function}\emph{extremal}.
More precisely, an admissible function $r\:\spc{X}\to\RR$ is extremal 
if for any admissible function $s\:\spc{X}\to\RR$ we have
\[s\le r\quad\Longrightarrow\quad s=r.\]

Applying Zorn's lemma, we get the following.

\begin{thm}{Observation}\label{obs:extremal:below}
For any admissible function $s$ there is an extremal function $r$ such that $r\le s$.
\end{thm}

\begin{thm}{Lemma}\label{lem:+-c}
Let $r$ be an extremal function and $s$ an admissible function on a metric space $\spc{X}$.
Suppose that $r\ge s-c$ for some constant~$c$.
Then $r\le s+c$; in particular, $c\ge 0$.
\end{thm}

\parit{Proof.}
Note that if $c<0$, then $r>s$.
The latter is impossible since $r$ is extremal and $s$ is admissible.

Observe that the function $\bar r=\min\{\,r,s+c\,\}$ is admissible.
Indeed, choose $x,y\in \spc{X}$.
If $\bar r(x)=r(x)$ and $\bar r(y)=r(y)$, then 
\[\bar r(x)+\bar r(y)=r(x)+ r(y)\ge \dist{x}{y}{}.\]
Further, if $\bar r(x)=s(x)+c$, then 
\begin{align*}
\bar r(x)+\bar r(y)&\ge [s(x)+c]+ [s(y)-c]= 
\\
&=s(x)+s(y) \ge 
\\
&\ge\dist{x}{y}{}.
\end{align*}

Since $r$ is extremal, we have $r=\bar r$;
that is, $r\le s+c$.
\qeds

\begin{thm}{Observations}\label{obs:extremal}
Let $\spc{X}$ be a metric space.

\begin{subthm}{obs:extremal:distfun}
For any point $p\in\spc{X}$ the distance function $r\z=\distfun_p$ is extremal.
\end{subthm}

\begin{subthm}{lem:extremal-lipschitz}
Any extremal function $r$ on $\spc{X}$ is \index{1-Lipschitz function}\emph{1-Lipschitz};
that is,
\[|r(p)-r(q)|\le \dist{p}{q}{}\]
for any $p,q\in\spc{X}$.
In other words, any extremal function is an extension function [see \ref{sec:Extension property}].
\end{subthm}

\begin{subthm}{lem:opposite}
An admissible function $r$ on $\spc{X}$ is extremal if and only if
for any point $p\in\spc{X}$ and any $\delta>0$, there is a point $q\in \spc{X}$
such that 
\[r(p)+r(q)<\dist{p}{q}{\spc{X}}+\delta.\]
\end{subthm}

\begin{subthm}{lem:opposite-compact}
Suppose $\spc{X}$ is compact.
Then an admissible function $r$ on $\spc{X}$ is extremal if and only if
for any point $p\in\spc{X}$ there is a point $q\in \spc{X}$
such that 
\[r(p)+r(q)=\dist{p}{q}{\spc{X}}.\]
\end{subthm}

\end{thm}

\parit{Proof; \ref{SHORT.obs:extremal:distfun}.}
By the triangle inequality, \ref{eq:admissible} holds;
that is, $r=\distfun_p$ is an admissible function.

Further, if $s\le r$ is another admissible function, then $s(p)=0$ and \ref{eq:admissible} implies that $s(x)\z\ge\dist{p}{x}{}$.
Whence $s=r$.

\parit{\ref{SHORT.lem:extremal-lipschitz}.}
By \ref{SHORT.obs:extremal:distfun}, $\distfun_p$ is admissible.
Since $r$ is admissible, we have that
\[r\ge \distfun_p-r(p).\]
Since $r$ is extremal, \ref{lem:+-c} implies that
\[r\le \distfun_p+r(p),\]
or, equivalently,
\[r(q)-r(p)\le \dist{p}{q}{}\]
for any $p,q\in\spc{X}$.
Whence the statement follows.

\parit{\ref{SHORT.lem:opposite}.}
Assume $r$ is extremal.
Arguing by contradiction, assume there is $\delta>0$ such that
\[r(q)\ge \distfun_p(q)-r(p)+\delta\]
for any $q$.
By \ref{SHORT.obs:extremal:distfun}, $\distfun_p$ is extremal; in particular, admissible.
Therefore \ref{lem:+-c} implies that
\[r(q)\le \distfun_p(q)+r(p)-\delta\]
for any $q$.
Taking $q=p$, we get $r(p)\le r(p)-\delta$, a contradiction.

Now suppose $r$ is not extremal; that is, there is an admissible function $s\le r$ such that $r(p)-s(p)=\delta>0$ for some $p$.
Then, for any $q$, we have
\[r(p)+r(q)\ge s(p)+s(q)+\delta\ge \dist{p}{q}{\spc{X}}+\delta\]
--- a contradiction.

\parit{\ref{SHORT.lem:opposite-compact}.}
The if part follows from \ref{SHORT.lem:opposite}.

Denote by $q_n$ the point provided by \ref{SHORT.lem:opposite} for $\delta=\tfrac1n$.
Let $q$ be a partial limit of $q_n$. 
Then 
\[r(p)+r(q)\le\dist{p}{q}{\spc{X}}.\]
Since $r$ is admissible, the opposite inequality holds;
whence the only-if part follows.
\qeds

\begin{thm}{Exercise}\label{ex:circle}
Consider the unit circle 
\[\mathbb{S}^1=\set{(x,y)}{x^2+y^2=1}\]
in the plane with induced length metric.
Show that $r\:\mathbb{S}^1\to\RR$ is extremal if and only if it is 1-Lipschitz and 
\[r(p)+r(-p)=\pi\] for any $p\in\mathbb{S}^1$.
\end{thm}

\begin{thm}{Exercise}\label{ex:retraction}
Given a real-valued function $s$ on a metric space $\spc{X}$,
consider the function
\[s^*(x)=\sup\set{\dist{x}{y}{\spc{X}}-s(y)}{y\in \spc{X}}\]
Show that the function $\tfrac12\cdot(s+s^*)$ is admissible for any $s$.
\end{thm}

\section{Equivalent conditions}

\begin{thm}{Theorem}\label{thm:injective=hyperconvex}
For any metric space $\spc{Y}$ the following conditions are equivalent:

\begin{subthm}{thm:injective=hyperconvex:injective}
$\spc{Y}$ is injective
\end{subthm}

\begin{subthm}{thm:injective=hyperconvex:extremal}
If $r\:\spc{Y}\to\RR$ is an extremal function, then there is a point $p\in \spc{Y}$ such that 
\[\dist{p}{x}{}= r(x)\]
for any $x\in \spc{Y}$.
\end{subthm}

\begin{subthm}{thm:injective=hyperconvex:balls}
$\spc{Y}$ is \index{hyperconvex space}\emph{hyperconvex};
that is, if $\set{\cBall[x_\alpha,r_\alpha]}{\alpha\in\IndexSet}$ is a family of closed balls in $\spc{Y}$ such that 
 \[r_\alpha+r_\beta\ge \dist{x_\alpha}{x_\beta}{}\]
 for any $\alpha,\beta\in \IndexSet$, then all the balls in the family $\{\cBall[x_\alpha,r_\alpha]\}_{\alpha\in\IndexSet}$ have a common point.
\end{subthm}

\end{thm}

\parit{Proof.} We will prove implications 
\ref{SHORT.thm:injective=hyperconvex:injective}$\Rightarrow$\ref{SHORT.thm:injective=hyperconvex:extremal}$\Rightarrow$\ref{SHORT.thm:injective=hyperconvex:balls}$\Rightarrow$\ref{SHORT.thm:injective=hyperconvex:injective}.

\parit{\ref{SHORT.thm:injective=hyperconvex:injective}$\Rightarrow$\ref{SHORT.thm:injective=hyperconvex:extremal}.}
By \ref{lem:extremal-lipschitz}, $r$ is an extension function.
Applying the definition of injective space to a one-point extension of $\spc{Y}$, we get a point $p\in \spc{Y}$ such that 
\[\dist{p}{x}{}=\distfun_p(x)\le r(x)\]
for any $x\in \spc{Y}$.
By \ref{obs:extremal:distfun}, the distance function $\distfun_p$ is extremal.
Since  $r$ is extremal, we get $\distfun_p= r$.

\parit{\ref{SHORT.thm:injective=hyperconvex:extremal}$\Rightarrow$\ref{SHORT.thm:injective=hyperconvex:balls}.}
By \ref{obs:admissible:balls}, part \ref{SHORT.thm:injective=hyperconvex:balls} is equivalent to the following statement:
\begin{itemize}
 \item If $r\:\spc{Y}\to\RR$ is an admissible function, then there is a point $p\in \spc{Y}$ such that 
\[\dist{p}{x}{}\le r(x)\eqlbl{eq:|p-x|=<r(x)}\]
for any $x\in \spc{Y}$.
\end{itemize}
Indeed, set $r(x)\df\inf\set{r_\alpha}{x_\alpha=x}$.
(If $x_\alpha\ne x$ for any $\alpha$, then $r(x)=\infty$.)
The condition in \ref{SHORT.thm:injective=hyperconvex:balls} implies that $r$ is admissible.
It remains to observe that $p\in \cBall[x_\alpha,r_\alpha]$ for every $\alpha$ if and only if \ref{eq:|p-x|=<r(x)} holds.

By \ref{obs:extremal:below}, for any admissible function $r$ there is an extremal function $\bar r\le r$;
hence \ref{SHORT.thm:injective=hyperconvex:extremal}$\Rightarrow$\ref{SHORT.thm:injective=hyperconvex:balls}.

\parit{\ref{SHORT.thm:injective=hyperconvex:balls}$\Rightarrow$\ref{SHORT.thm:injective=hyperconvex:injective}.}
Arguing by contradiction, suppose $\spc{Y}$ is not injective;
that is, there is a metric space $\spc{X}$ with a subset $\spc{A}$
such that a short map $f\:\spc{A}\to \spc{Y}$ cannot be extended to a short map $F\:\spc{X}\to \spc{Y}$.
By Zorn's lemma, we may assume that $\spc{A}$ is a maximal subset; that is, the domain of $f$ cannot be enlarged by a single point.%
\footnote{In this case, $\spc{A}$ must be closed, but we will not use it.}

Fix a point $p$ in the complement $\spc{X}\setminus \spc{A}$.
To extend $f$ to $p$, we need to choose $f(p)$ in the intersection of the balls 
$\cBall[f(x),r(x)]$, where $r(x)=\dist{p}{x}{}$.
Therefore, this intersection for all $x\in \spc{A}$ has to be empty.

Since $f$ is short, we have that 
\begin{align*}
r(x)+r(y)&\ge \dist{x}{y}{\spc{X}}\ge
\\
&\ge \dist{f(x)}{f(y)}{\spc{Y}}.
\end{align*}
By \ref{SHORT.thm:injective=hyperconvex:balls} the balls 
$\cBall[f(x),r(x)]$ have a common point --- a contradiction. 
\qeds

\begin{thm}{Exercise}\label{ex:one-point-gluing}
Suppose a length space $\spc{W}$ has two subspaces $\spc{X}$ and $\spc{Y}$ such that $\spc{X}\cup\spc{Y}=\spc{W}$ and $\spc{X}\cap\spc{Y}$ is a one-point set.
Assume $\spc{X}$ and $\spc{Y}$ are injective.
Show that  $\spc{W}$ is injective
\end{thm}

\begin{thm}{Exercise}\label{ex:Rm-ell-infty}
Show that an $m$-dimensional normed space is injective if and only if it is isometric to $\RR^m$ with $\ell^\infty$-norm; that is,
\[|(x_1,\dots,x_m)|=\max_i\{\,|x_i|\,\}.\]
\end{thm}

A metric space $\spc{Y}$ is called \index{finitely hyperconvex}\emph{finitely hyperconvex} or \index{countably hyperconvex}\emph{countably hyperconvex} if the condition in \ref{thm:injective=hyperconvex:balls} holds only for any finite or respectively countable family of balls.

\begin{thm}{Exercise}\label{ex:compact-hyperconvex}
Show that any proper finitely hyperconvex metric space is hyperconvex.
\end{thm}

\begin{thm}{Exercise}\label{ex:urysohn-hyperconvex}
Show that the $d$-Urysohn space is finitely hyperconvex, but not countably hyperconvex.
Conclude that the $d$-Urysohn space is not injective.

Try to do the same for the Urysohn space.
\end{thm}

\begin{thm}{Exercise}\label{ex:almost-hyperconvex}
Let $\spc{Y}$ be a complete metric space.
Suppose $\spc{Y}$ is \index{almost hyperconvex}\emph{almost hyperconvex},
that is, for any $\eps>0$ any family of closed balls $\set{\cBall[x_\alpha,r_\alpha+\eps]}{\alpha\in\IndexSet}$ has a common point if 
$r_\alpha+r_\beta\ge \dist{x_\alpha}{x_\beta}{}$ for all $\alpha,\beta\in \IndexSet$.
Show that $\spc{Y}$ is hyperconvex.
\end{thm}

\section{Space of extremal functions}
\label{sec:extremal-functions}

Let $\spc{X}$ be a metric space.
Consider the space $\Inj \spc{X}$ of extremal functions on $\spc{X}$ equipped with sup-norm; \label{page:InjX}
that is,
\[\dist{f}{g}{\Inj \spc{X}}\df\sup\set{|f(x)-g(x)|}{x\in \spc{X}}.\]

Recall that by \ref{obs:extremal:distfun}, any distance function is extremal.
It follows that the map $x\mapsto \distfun_x$ produces a distance-preserving embedding $\spc{X}\hookrightarrow\Inj \spc{X}$.
So we can (and will) treat $\spc{X}$ as a subspace of $\Inj \spc{X}$,
or, equivalently, $\Inj \spc{X}$ as an extension of $\spc{X}$.
In particular, from now on, a point $x\in\spc{X}$ can refer to the function $\distfun_x\:\spc{X}\to\RR$ and the other way around.

Since any extremal function is 1-Lipschitz, for any $f\in \Inj \spc{X}$ and $p\in \spc{X}$, we have that
$f(x)\le f(p)+\distfun_p(x)$.
By \ref{lem:+-c}, we also get $f(x)\ge -f(p)+\distfun_p(x)$.
Therefore
\[
\begin{aligned}
\dist{f}{p}{\Inj \spc{X}}&=\sup\set{|f(x)-\distfun_p(x)|}{x\in \spc{X}}=
\\
&=f(p).
\end{aligned}
\eqlbl{eq:f(p)=|f-p|}
\]
In particular, the statement in \ref{lem:opposite} can be written as 
\[\dist{f}{p}{\Inj\spc{X}}+\dist{f}{q}{\Inj\spc{X}}<\dist{p}{q}{\Inj\spc{X}}+\delta.\]

\begin{thm}{Exercise}\label{ex:Inj(compact)}
Show that $\Inj\spc{X}$ is compact if and only if so is $\spc{X}$.
\end{thm}

\begin{thm}{Exercise}\label{ex:tripod+square}
Describe the set of all extremal functions on a metric space $\spc{X}$ and the metric space $\Inj \spc{X}$ in each of the following cases:

\begin{subthm}{ex:tripod+square:2}
$\spc{X}$ is a metric space with exactly two points $v,w$ on distance 1 from each other.
\end{subthm}

\begin{subthm}{ex:tripod+square:tripod} 
$\spc{X}$ is a metric space with exactly three points $a,b,c$ such that 
\[\dist{a}{b}{\spc{X}}=\dist{b}{c}{\spc{X}}=\dist{c}{a}{\spc{X}}=1.\]
\end{subthm}

\begin{subthm}{ex:tripod+square:square}
$\spc{X}$ is  a metric space with exactly four points $p,q,x,y$ such that 
\[\dist{p}{x}{\spc{X}}=\dist{p}{y}{\spc{X}}=\dist{q}{x}{\spc{X}}=\dist{q}{y}{\spc{X}}=1\]
and
\[\dist{p}{q}{\spc{X}}=\dist{x}{y}{\spc{X}}=2.\]
\end{subthm}

\end{thm}

\begin{thm}{Exercise}\label{ex:kur-inj}
Assume $\spc{X}$ is a compact metric space.
Recall that the map $x\mapsto \distfun_x$ gives an isometric embedding $\spc{X}\hookrightarrow\ell^\infty(\spc{X})$; so we can think that $\spc{X}$ is a subset of $\ell^\infty(\spc{X})$.

Given two points $x,y\in \spc{X}$, denote by $G_{x,y}$ the union of all geodesics from $x$ to $y$ in $\ell^\infty(\spc{X})$.
Show that $\Inj\spc{X}$ is isometric to
\[G=\bigcap_{x\in \spc{X}}\left(\bigcup_{y\in \spc{X}}G_{x,y}\right).\]

\end{thm}

\begin{thm}{Proposition}\label{prop:InjX-is-injective}
$\Inj\spc{X}$ is injective for any metric space $\spc{X}$. 
\end{thm}

\begin{thm}{Lemma}\label{lem:r|X-extremal}
Let $\spc{X}$ be a metric space.
Then 
\[\sigma\in \Inj(\Inj \spc{X})
\quad\Longrightarrow\quad
\sigma|_\spc{X}\in \Inj \spc{X}.\]
\end{thm}

In other words, if $\sigma$ is an extremal function on $\Inj \spc{X}$,
then the restriction of $\sigma$ to $\spc{X}$ is an extremal function on $\spc{X}$.

\parit{Proof.}
Arguing by contradiction, suppose that there is an admissible function $s\:\spc{X}\to \RR$ such that $s(x)\le \sigma(x)$ for any $x\in\spc{X}$ and $s(p)\z< \sigma(p)$ for some point $p\in\spc{X}$.
Consider another function $\bar \sigma\:\Inj \spc{X}\to\RR$ such that $\bar \sigma(f)\df \sigma(f)$ if $f\ne p$ and $\bar \sigma(p)\df s(p)$.

Let us show that $\bar \sigma$ is admissible; that is, 
\[\dist{f}{g}{\Inj \spc{X}}\le\bar \sigma(f)+\bar \sigma(g)
\eqlbl{r-admissible}\]
for any $f,g\in \Inj \spc{X}$.

Since $\sigma$ is admissible and $\bar \sigma= \sigma$ on $(\Inj \spc{X})\setminus \{p\}$, it is sufficient to prove \ref{r-admissible} assuming $f\ne g=p$.
By \ref{eq:f(p)=|f-p|}, we have $\dist{f}{p}{\Inj \spc{X}}=f(p)$.
Therefore, \ref{r-admissible} boils down to the following inequality
\[\sigma(f)+s(p)\ge f(p).\eqlbl{eq:r(f)+s(p)>=f(p)}\]
for any $f\in\Inj \spc{X}$.

Fix small $\delta>0$. 
Let $q\in\spc{X}$ be the point provided by \ref{lem:opposite}.
Then
\begin{align*}
\sigma(f)+s(p)&\ge [\sigma(f)-\sigma(q)]+[\sigma(q)+s(p)]\ge
\intertext{since $\sigma$ is 1-Lipschitz, and $\sigma(q)\ge s(q)$, we can continue}
&\ge -\dist{q}{f}{\Inj \spc{X}}+[s(q)+s(p)]\ge
\intertext{by \ref{eq:f(p)=|f-p|} and since $s$ is admissible}
&\ge -f(q)+\dist{p}{q}{}>
\intertext{and by \ref{lem:opposite}}
&> f(p)-\delta.
\end{align*}
Since $\delta>0$ is arbitrary, \ref{eq:r(f)+s(p)>=f(p)} and \ref{r-admissible} follow.

Summarizing: the function $\bar \sigma$ is admissible, $\bar \sigma\le \sigma$ and $\bar \sigma(p)<\sigma(p)$;
that is, $\sigma$ is not extremal --- a contradiction.
\qeds

\parit{Proof of \ref{prop:InjX-is-injective}.}
Choose a function $\sigma\in\Inj(\Inj\spc{X})$.
By \ref{lem:r|X-extremal}, $s\z\df \sigma|_{\spc{X}}\in \Inj\spc{X}$;
that is, $s$ is extremal.
By \ref{thm:injective=hyperconvex:extremal},
it is sufficient to show that  
\[\sigma(f)\ge\dist{s}{f}{\Inj\spc{X}}
\eqlbl{eq:r(f)>=| r-f|}\]
for any $f\in\Inj\spc{X}$.

Since $\sigma$ is $1$-Lipschitz (\ref{lem:extremal-lipschitz}) we have that
\[
s(x)-f(x)=\sigma(x)-\dist{f}{x}{\Inj \spc{X}}\le \sigma(f).
\]
for any $x\in\spc{X}$.
By \ref{lem:+-c},
$
s(x)-f(x)\ge -\sigma(f)
$
for any $x\in\spc{X}$.
Whence \ref{eq:r(f)>=| r-f|} follows.
\qeds

\begin{thm}{Exercise}\label{ex:4-on-a-line}
Let $\spc{X}$ be a compact metric space.
Show that for any two points $f,g\in\Inj \spc{X}$ lie on a geodesic $[pq]$ with $p,q\in \spc{X}$.
\end{thm}

A metric space $\spc{X}$ is called \index{$\delta$-hyperbolic}\emph{$\delta$-hyperbolic} if 
\[\dist{p}{q}{}+\dist{x}{y}{}\le
\max\{\,\dist{p}{x}{}+\dist{q}{y}{},
\,
\dist{p}{y}{}+\dist{q}{x}{}\,\}+2\cdot\delta\]
for any $p,q,x,y\in \spc{X}$.

\begin{thm}{Advanced exercise}\label{ex:delta-hyp}
Show that $\Inj \spc{X}$ is $\delta$-hyperbolic if and only if so is $\spc{X}$.
\end{thm}

\section{Injective envelope}

An extension $\spc{E}$ of a metric space $\spc{X}$ will be called its \index{injective envelope}\emph{injective envelope} if $\spc{E}$ is an injective space, and there is no proper injective subspace of $\spc{E}$ that contains $\spc{X}$.

Two injective envelopes $e\:\spc{X}\hookrightarrow \spc{E}$ and $f\:\spc{X}\hookrightarrow \spc{F}$ are called  equivalent if there is an isometry $\iota\: \spc{E}\to\spc{F}$ such that $f=\iota\circ e$.

\begin{thm}{Theorem}\label{thm:inj-envelope}
For any metric space $\spc{X}$, its extension $\Inj\spc{X}$ is an injective envelope.

Moreover, any other injective envelope of $\spc{X}$ is equivalent to $\Inj\spc{X}$.
\end{thm}

\parit{Proof.} 
Suppose $S\subset \Inj\spc{X}$ is an injective subspace containing $\spc{X}$.
Since $S$ is injective, there is a short map $w\:\Inj\spc{X}\to S$ that fixes all points in $\spc{X}$.

Suppose that $w\:f\mapsto f'$; observe that $f(x)\ge f'(x)$ for any $x\in \spc{X}$.
Since $f$ is extremal, $f=f'$;
that is, $w$ is the identity map, and therefore $S=\Inj\spc{X}$.

Assume we have another injective envelope $e\:\spc{X}\hookrightarrow \spc{E}$.
Then there are short maps $v\:\spc{E}\to \Inj\spc{X}$ and $w\:\Inj\spc{X}\to \spc{E}$ such that $x=v\circ e(x)$ and $e(x)=w(x)$ for any $x\in\spc{X}$.
From above, the composition $v\circ w$ is the identity on $\Inj\spc{X}$.
In particular, $w$ is distance-preserving.

The composition $w\circ v\:\spc{E}\to \spc{E}$ is a short map that fixes points in $e(\spc{X})$.
Since $e\:\spc{X}\hookrightarrow \spc{E}$ is an injective envelope, the composition $w\circ v$ and therefore $w$ are onto.
Whence $w$ is an isometry.
\qeds

\begin{thm}{Exercise}\label{ex:inj-envelope}
Suppose $e\:\spc{X}\hookrightarrow \spc{E}$ and $f\:\spc{X}\hookrightarrow \spc{F}$ are two injective envelopes of $\spc{X}$.
Show that there is a unique isometry $\iota\:\spc{E}\to \spc{F}$ such that $\iota\circ e=f$.
\end{thm}

\begin{thm}{Exercise}\label{ex:d-p-inclusion}
Suppose $\spc{X}$ is a subspace of a metric space $\spc{U}$.
Show that the inclusion $\spc{X}\hookrightarrow\spc{U}$ can be extended to a distance-preserving inclusion $\Inj\spc{X}\hookrightarrow\Inj\spc{U}$.
\end{thm}

\begin{thm}{Exercise}\label{ex:hemisphere-inj}
Consider the hemisphere 
\begin{align*}
\mathbb{S}^2_+&=\set{(x,y,z)\in\RR^3}{x^2+y^2+z^2=1,\quad z\ge0}
\intertext{and its boundary}
\mathbb{S}^1&=\set{(x,y,z)\in\RR^3}{x^2+y^2+z^2=1,\quad z=0};
\end{align*}
 both with induced length metrics.
 
Show that there is unique isometric embedding $\iota\:\mathbb{S}^2_+\hookrightarrow\Inj\mathbb{S}^1$ such that $\iota(u)=u$ for any $u\in \mathbb{S}^1$.
\end{thm}

\section{Remarks}

Injective spaces were introduced by Nachman Aronszajn and Prom Panitchpakdi \cite{aronszajn-panitchpakdi}.
The injective envelope was introduced by John Isbell \cite{isbell}; it is also known as \index{tight span}\emph{tight span} and \index{hyperconvex hull}\emph{hyperconvex hull}.

It was observed by John Isbell \cite{isbell2} that \textit{if $\spc{X}$ is a Banach space, then its injective hull $\Inj\spc{X}$ has a natural structure of Banach space} (which is unique by the Mazur--Ulam theorem).
Moreover, $\spc{X}$ is a linear subspace of $\Inj\spc{X}$.
 
Let us mention that a metric space $\spc{X}$ is called \index{convex space}\emph{convex} if for any pair of points $x_1,x_2\in \spc{X}$ and any $r_1,r_2\ge 0$ we have 
\[r_1+r_2\ge \dist{x_1}{x_2}{\spc{X}}\qquad\Longrightarrow\qquad\cBall[x,r_1]_\spc{X}\cap \cBall[y,r_2]_\spc{X}\ne\emptyset;\]
in other words, a pair of balls intersect if the triangle inequality does not forbid it.
Clearly, hyperconvexity (\ref{thm:injective=hyperconvex:balls}) is stronger than convexity.
Note that \textit{any geodesic space is convex}.
The converse does not hold in general, but by Menger's lemma (\ref{lem:mid>geod:geod}) \textit{any complete convex space is geodesic}.

More generally, a metric space $\spc{X}$ is called \index{$n$-hyperconvex space}\emph{$n$-hyperconvex} if the condition in \ref{thm:injective=hyperconvex:balls} holds only for families with at most $n$ balls; so \textit{convex means $2$-hyperconvex}.

The following striking result was proved by Benjamin Miesch and Maël Pavón \cite{miesch-pavon2016}.

\begin{thm}{Theorem}
Any complete $4$-hyperconvex space is finitely hyperconvex.
\end{thm}

So, by \ref{ex:compact-hyperconvex}, it follows that \textit{any proper $4$-hyperconvex space is hyperconvex}.

\begin{thm}{Exercise}\label{ex:3-4-hypreconvex}
Show that $\ell^1$ is $3$- but not $4$-hyperconvex.
\end{thm}

Recall that if the following inequality
\[\dist{x}{z}{\spc{X}}
\le
\max\{\,\dist{x}{y}{\spc{X}},\dist{y}{z}{\spc{X}}\,\}\]
holds for any three points $x,y,z$ in a metric space $\spc{X}$,
then $\spc{X}$ is called an \index{ultrametric space}\emph{ultrametric space}.
In some sense, ultrametric spaces are dual to injective spaces.

\begin{thm}{Exercise}\label{ex:ultrametric}
Suppose that a metric space $\spc{X}$ satisfies the following property:
For any subspace $\spc{A}$ in $\spc{X}$ and any other metric space $\spc{Y}$, any short map $f\:\spc{A}\to \spc{Y}$ can be extended to a short map $F\:\spc{X}\to \spc{Y}$.

Show that $\spc{X}$ is an ultrametric space.
\end{thm}

A subspace $\spc{S}$ of a metric space $\spc{X}$ is called its \index{short retract}\emph{short retract} if there is a short map $\spc{X}\to \spc{S}$ that is the identity on $\spc{S}$.

\begin{thm}{Exercise}\label{ex:ultrametric-converse}
Show that any compact subspace $\spc{K}$ of an ultrametric space $\spc{X}$ is its short retract.

Construct an example of a complete ultrametric space $\spc{X}$ with a closed subspace $\spc{Q}$ that is not its short retract.
\end{thm}

The following exercise gives a sufficient condition for the existence of a short extension.

\begin{thm}{Exercise}\label{ex:petrunin-stadler}
Let $f\:A\z\to \spc{K}$ be a short map from a subset $A$ in a metric space $\spc{X}$ to compact metric space $\spc{K}$.
Assume that for any finite set $F\subset \spc{X}$ there is a short map $F\to \spc{K}$ that agrees with $f$ on $F\cap A$.
Show that there is a short map $\spc{X}\to \spc{K}$ that agrees with $f$ on $A$.
\end{thm}

%% file: converge.tex
\chapter{Space of subsets}\label{chap:hausdorff}

In this lecture we define and study Hausdorff metric on subsets of a given metric space.

\section{Hausdorff distance}

Let $\spc{X}$ be a metric space.
Given a subset $A\subset \spc{X}$,
consider the distance function to $A$
$$\distfun_A: \spc{X} \to [0,\infty)$$
defined as 
$$\distfun_A(x)
\df
\inf\set{\dist ax{\spc{X}}}{a\in A}.$$

Further, we define the so-called Hausdorff metric on all nonempty compact subsets of a given metric space $\spc{X}$.
The obtained metric space will be denoted as $\Haus\spc{X}$.

\begin{thm}{Definition}\label{def:hausdorff-convergence}
Let $A$ and $B$ be two nonempty compact subsets of a metric space $\spc{X}$.
Then the \index{Hausdorff distance}\emph{Hausdorff distance} between $A$ and $B$ is defined as 
$$|A-B|_{\Haus\spc{X}}
\df
\sup_{x\in \spc{X}}\{\,|\distfun_A(x)-\distfun_B(x)|\,\}.
$$

\end{thm}

The following observation gives a useful reformulation of the definition:

\begin{thm}{Observation}\label{obs:Haus-nbhds}
Suppose $A$ and $B$ be two compact subsets of a metric space $\spc{X}$.
Then $|A-B|_{\Haus\spc{X}}< R$ if and only if and only if 
$B$ lies in an $R$-neighborhood of $A$, 
and 
$A$ lies in an $R$-neighborhood of~$B$.
\end{thm}

\begin{thm}{Exercise}\label{ex:diam}
Let $\spc{X}$ be a metric space.
Given a subset $A\subset \spc{X}$, define its \index{diameter}\emph{diameter} as 
$$\diam A\df\sup_{a,b\in A} |a-b|.$$

Show that 
$$\diam\:\Haus\spc{X}\to \RR$$ 
is a $2$-Lipschitz function;
that is,
\[|\diam A-\diam B|\le 2\cdot\dist{A}{B}{\Haus\spc{X}}\]
for any two compact nonempty sets $A,B\subset\spc{X}$.
\end{thm}

\begin{thm}{Exercise}\label{ex:Hausdorff-bry}
Let $A$ and $B$ be two compact subsets in the euclidean plane $\RR^2$.
Assume $|A-B|_{\Haus\RR^2}<\eps$.

\begin{subthm}{ex:Hausdorff-bry:conv}
Show that $|\Conv A-\Conv B|_{\Haus\RR^2}<\eps$, where $\Conv A$ denoted the convex hull of $A$.
\end{subthm}
\begin{subthm}{ex:Hausdorff-bry:bry}
Is it true that
$|\partial A-\partial B|_{\Haus\RR^2}<\eps$,
where $\partial A$ denotes the boundary of $A$.

Does the converse hold? That is, assume $A$ and $B$ be two compact subsets in $\RR^2$
and $|\partial A-\partial B|_{\Haus\RR^2}<\eps$; 
is it true that $|A-B|_{\Haus\RR^2}\z<\eps$?
\end{subthm}

\end{thm}

Note that part \ref{SHORT.ex:Hausdorff-bry:conv} implies that $A\mapsto \Conv A$ defines a short map $\Haus\RR^2\to \Haus\RR^2$. 

\begin{thm}{Exercise}\label{ex:Haus-func}
Let $A$ and $B$ be compact subsets in metric space~$\spc{X}$.
Show that 
\[\dist{A}{B}{\Haus\spc{X}}=\sup_f\, \{\,\max_{a\in A}\{f(a)\}-\max_{b\in B}\{f(b)\,\},\]
where the least upper bound is taken for all $1$-Lipschitz functions $f$.

\end{thm}

Given a subset $A\subset \RR^n$,
the \index{support function}\emph{support function} $h_A\colon\mathbb{R}^n\to\mathbb{R}$ of a  nonempty closed set $A\subset \RR^n$ is defined as 
\[h_A(x)
\df
\sup\set{\langle x, a\rangle}{a\in A}.\]

\begin{thm}{Exercise}\label{ex:Haus-support}
Show that 
\[\dist{A}{B}{\Haus\RR^n}\ge \sup_{|u|=1}\{|h_A(u)-h_B(u)|\}\]
for any nonempty compact subsets $A,B\subset \RR^n$.

Moreover, equality holds if both $A$ and $B$ are convex.
\end{thm}

\begin{thm}{Advanced exercise}\label{ex:H-sections}
Suppose $C_t\subset \spc{X}$, $t\z\in [0,1]$ is a family of subsets.
A path $c\:[0,1]\to \spc{X}$ such that $c(t)\in C_t$ for all $t$ will be called a {}\emph{section} of $C_t$.

\begin{subthm}{ex:H-sections:S}
Construct a family of nonempty compact sets $C_t\subset\mathbb{S}^1$, $t\z\in [0,1]$ that is continuous in the Hausdorff topology, 
but does not admit a section.
\end{subthm}

\begin{subthm}{ex:H-sections:R}
Show that any family of nonempty compact sets $C_t\subset\RR$, $t\z\in [0,1]$ that is continuous in the Hausdorff topology, 
admits a section.
\end{subthm}

\end{thm}

\section{Hausdorff convergence}

\begin{thm}{Blaschke selection theorem}\label{thm:compact+Hausdorff}
A metric space $\spc{X}$ is compact if and only if
so is $\Haus\spc{X}$.
\end{thm}

The Hausdorff metric can be used to define convergence.
Namely, suppose $K_1,K_2,\dots$, and $K_\infty$ are compact sets in a metric space $\spc{X}$.
If $|K_\infty-K_n|_{\Haus\spc{X}}\to0$ as $n\to\infty$, then we say that 
the sequence $K_n$ {}\emph{converges} to $K_\infty$ \index{Hausdorff convergence}\emph{in the sense of Hausdorff};
equivalently, $K_\infty$ is the \index{Hausdorff limit}\emph{Hausdorff limit} of the sequence $K_n$.

Note that the theorem implies that from any sequence of nonempty compact sets in $\spc{X}$ one can select a convergent subsequence; 
for that reason, it is called a \textit{selection} theorem. 

\parit{Proof; if part.}
Consider the map $\iota$ that sends each point $x\in \spc{X}$ to the one-point subset $\{x\}$ of $\spc{X}$.
Note that $\iota\:\spc{X}\to \Haus\spc{X}$ is distance-preserving.

Suppose that $A\subset \spc{X}$.
Note that $\diam A=0$ if and only if $A$ is a one-point set.
By \ref{ex:diam}, $\iota(\spc{X})$ is a closed subset of the compact space $\Haus\spc{X}$.
It follows that $\iota(\spc{X})$, and therefore $\spc{X}$, are compact.
\qeds

Since the map $\iota$ above is distance-preserving, we can and will consider $\spc{X}$ as a subspace of $\Haus\spc{X}$.

\begin{thm}{Exercise}\label{ex:haus-contractible}
Let $\spc{X}$ be a compact length space.
Suppose that there is a short retraction $\Haus\spc{X}\to \spc{X}$.
Show that $\spc{X}$ is contractible.
\end{thm}

To prove the only-if part we will need the following two lemmas.

\begin{thm}{Monotone convergence}\label{lem:decreasing-converges}
Let $K_1\supset K_2\supset\dots$ be a nested sequence of nonempty compact sets in a metric space $\spc{X}$.
Then $K_\infty\z=\bigcap_n K_n$ is the Hausdorff limit of $K_n$;
that is, $|K_\infty-K_n|_{\Haus\spc{X}}\to0$ as $n\to\infty$.
\end{thm}

\parit{Proof.}
By finite intersection property, $K_\infty$ is a nonempty compact set.

Arguing by contradiction, assume that there is $\eps>0$ such that for each $n$ 
one can choose $x_n\in K_n$
such that $\distfun_{K_\infty}(x_n)\ge\eps$.
Note that $x_n\in K_1$ for each $n$.
Since $K_1$ is compact, 
there is 
a \index{partial limit}\emph{partial limit}
 $x_\infty$ of $x_n$;
that is, a limit of a subsequence.
Clearly, $\distfun_{K_\infty}(x_\infty)\ge \eps$.

On the other hand, since $K_n$ is closed and $x_m\in K_n$ for $m\ge n$,
we get $x_\infty\in K_n$ for each $n$.
It follows that $x_\infty\in K_\infty$ and therefore $\distfun_{K_\infty}(x_\infty)=0$ ---
a contradiction.\qeds

\begin{thm}{Lemma}\label{lem:complete+Hausdorff}
If $\spc{X}$ is a compact metric space, then $\Haus\spc{X}$
is complete.
\end{thm}

\parit{Proof.}
Let $Q_1,Q_2,\dots$ be a Cauchy sequence in $\Haus\spc{X}$.
Passing to a subsequence, we may assume that 
$$|Q_n-Q_{n+1}|_{\Haus\spc{X}}\le \tfrac1{10^n}\eqlbl{eq:eps=1/10}$$
for each $n$.

Denote by $K_n$ the closed $\tfrac2{10^n}$-neighborhood of $Q_n$;
that is,
\begin{align*}
K_n&= \set{x\in \spc{X}}{\distfun_{Q_n}(x)\le \tfrac2{10^n}}
\end{align*}
Since $\spc{X}$ is compact so is each $K_n$.

From \ref{eq:eps=1/10}, we get
$K_n\supset K_{n+1}$ 
for each $n$.
Set 
$$K_\infty=\bigcap_{n=1}^\infty K_n.$$
By the monotone convergence (\ref{lem:decreasing-converges}),
 $|K_n-K_\infty|_{\Haus\spc{X}}\to 0$ as $n\to\infty$.

By \ref{obs:Haus-nbhds}, $|Q_n-K_n|_{\Haus\spc{X}}\le \tfrac2{10^n}$.
Therefore, $|Q_n-K_\infty|_{\Haus\spc{X}}\to 0$ as $n\to\infty$ --- hence the lemma.
\qeds

\begin{thm}{Exercise}\label{ex:closure-union}
Let $\spc{X}$ be a complete metric space and $K_1,K_2,\dots$ be a sequence of compact sets 
that converges in the sense of Hausdorff.
Show that the union $K_1\cup K_2\cup\dots$ has compact closure.

Use this statement to show that in Lemma~\ref{lem:complete+Hausdorff} compactness of $\spc{X}$ can be exchanged to completeness.
\end{thm}

\parit{Proof of only-if part in \ref{thm:compact+Hausdorff}.}
According to Lemma~\ref{lem:complete+Hausdorff},
$\Haus\spc{X}$ is complete.
It remains to show that $\Haus\spc{X}$ is totally bounded (\ref{totally-bounded});
that is, given $\eps>0$ there is a finite $\eps$-net in $\Haus\spc{X}$.

Choose a finite $\eps$-net $A$ in $\spc{X}$.
Denote by $B$ the set of all nonempty subsets of $A$.
Note that  $B$ is a finite set in $\Haus\spc{X}$.
For each compact set $K\subset \spc{X}$, consider the subset $K'$ of all points $a\in A$
such that $\distfun_K(a)\le \eps$.
Observe that $K' \in B$ and $|K-K'|_{\Haus\spc{X}}\le\eps$.
In other words, $B$ is a finite $\eps$-net in $\Haus\spc{X}$.
\qeds

\begin{thm}{Exercise}\label{ex:Haus-length}
Let $\spc{X}$ be a complete metric space.
Show that $\spc{X}$ is a length space if and only if so is $\Haus\spc{X}$.
\end{thm}

\begin{thm}{Exercise}\label{ex:Haus-G-delta}

\begin{subthm}{ex:Haus-G-delta:closed}
Show that the set of all connected compact subsets of $\RR^2$ is closed in $\Haus\RR^2$. 
\end{subthm}

\begin{subthm}{ex:Haus-G-delta:curves}
Show that any connected compact subset of $\RR^2$ is a Hausdorff limit of a sequence of closed simple curves.
\end{subthm}

\end{thm}

\section{An application}

In this section, we will sketch a proof of the isoperimetric inequality in the plane that uses the Hausdorff convergence.

It is based on the following exercise.

\begin{thm}{Exercise}\label{ex:Huas-perimeter-area}
Let $\spc{C}$ be the set of all nonempty compact convex subsets in $\RR^2$.
Show that $\spc{C}$ is a closed subset of $\Haus\RR^2$ and 
perimeter and area are continuous on~$\spc{C}$.
(If the set degenerates to a line segment of length $\ell$, then its perimeter is defined as $2\cdot \ell$.)

More precisely, if a sequence of convex compact plane sets $X_n$ converges to $X_\infty$ in the sense of Hausdorff, then $X_\infty$ is convex,
\[\perim X_n\to \perim X_\infty,\quad\text{and}\quad\area X_n\to\area X_\infty\]
as $n\to\infty$.
\end{thm}

\begin{thm}{Isoperimetric inequality}\label{thm:isoperimetric}
Among the plane figures bounded by closed curves of length at most $\ell$, the round disk has the maximal area.
\end{thm}

\parit{Sketch.}
It is sufficient to consider only convex figures of the given perimeter; if a figure is not convex, pass to its convex hull and observe that it has a larger area and smaller perimeter.

Note that the selection theorem (\ref{thm:compact+Hausdorff}) together with the exercise implies the existence of figure $D$ with perimeter $\ell$ and maximal area.

It remains to show that $D$ is a round disk;
it will be done by means of elementary geometry.

Let us cut $D$ along a chord $[ab]$ into two lenses, $L_1$ and $L_2$.
Denote by $L_1'$ the reflection of $L_1$ across the perpendicular bisector of $[ab]$.
Note that $D$ and $D'=L_1'\cup L_2$ have the same perimeter and area.
That is, $D'$ has perimeter $\ell$ and maximal possible area;
in particular, $D'$ is convex.

The following exercise will finish the proof.
\qeds

{

\begin{wrapfigure}{r}{57 mm}
\vskip-0mm
\centering
\includegraphics{mppics/pic-405}
\end{wrapfigure}

\begin{thm}{Exercise}\label{ex:round-disc}
Suppose $D$ is a convex figure such that for any chord $[ab]$ of $D$ the above construction produces a convex figure $D'$.
Show that $D$ is a round disk.
\end{thm}

}

Another popular way to prove that $D$ is a round disk is given by the so-called \textit{Steiner's 4-joint method} \cite{blaschke}.

\section{Remarks}\label{sec:H-variation}

It seems that Hausdorff convergence was first introduced by Felix Hausdorff~\cite{hausdorff}.
A couple of years later an equivalent definition was given by Wilhelm Blaschke~\cite{blaschke}.

The following refinement was introduced by  Zdeněk Frolík \cite{frolik} and rediscovered by Robert Wijsman~\cite{wijsman}.  
This refinement is also called \index{Hausdorff convergence}\emph{Hausdorff convergence};
in fact, it takes an intermediate place between the original Hausdorff convergence and the so-called \textit{closed convergence}, also introduced by Hausdorff in \cite{hausdorff}.

\begin{thm}{Definition}\label{def:gen-Haus-conv}
Let $A_1,A_2,\dots$ be a sequence of closed sets in a metric space $\spc{X}$.
We say that the sequence $A_n$ converges to a closed set $A_\infty$ in the sense of Hausdorff if, for any $x\in\spc{X}$, we have
$\distfun_{A_n}(x)\z\to \distfun_{A_\infty}(x)$ as $n\to\infty$.
\end{thm}

For example, suppose $\spc{X}$ is the euclidean plane and $A_n$ is the circle with radius $n$ and center at the point $(0,n)$.
If we use the standard definition (\ref{def:hausdorff-convergence}), then the sequence $A_1,A_2,\dots$ diverges, but it converges to the $x$-axis in the sense of Definition~\ref{def:gen-Haus-conv}.

\begin{figure}[ht!]
\vskip-0mm
\centering
\includegraphics{mppics/pic-415}
\end{figure}

Further, consider the sequence of one-point sets $B_n=\{(n,0)\}$ in the euclidean plane.
It diverges in the sense of the standard definition, but, in the sense of \ref{def:gen-Haus-conv}, it converges to the empty set;
indeed, for any point $x$ we have $\distfun_{B_n}(x)\to\infty$ as $n\to \infty$ and $\distfun_{\emptyset}(x)= \infty$ for any~$x$.

The following exercise is analogous to the Blaschke selection theorem (\ref{thm:compact+Hausdorff}) for the modified Hausdorff convergence.

\begin{thm}{Exercise}\label{ex:generalized-selection}
Let $\spc{X}$ be a proper metric space
and $A_1,A_2,\dots$ be a sequence of closed sets in~$\spc{X}$.
Show that the sequence  $A_1,A_2,\dots$ has a convergent subsequence in the sense of Definition~\ref{def:gen-Haus-conv}.
\end{thm}

\chapter{Space of spaces}\label{chap:GH}

In this lecture we define and study the so-called Gromov--Hausdorff metric on the isometry classes of compact metric spaces.

\section{Gromov--Hausdorff metric}

The goal of this section is to cook up a metric space out of all compact metric spaces.
More precisely, we want to define the so-called  Gromov--Hausdorff metric on the set of \textit{isometry classes} of compact metric spaces.
(Being isometric is an equivalence relation, 
and an \index{isometry class}\emph{isometry class} is an equivalence class with respect to this relation.)

The obtained metric space will be denoted by $\GH$.
Given two metric spaces $\spc{X}$ and $\spc{Y}$,
denote by $[\spc{X}]$ and $[\spc{Y}]$ their isometry classes;
that is, $\spc{X}'\in [\spc{X}]$ if and only if $\spc{X}'\iso \spc{X}$.
Pedantically, the Gromov--Hausdorff distance from $[\spc{X}]$ 
to $[\spc{Y}]$ should be denoted as $|[\spc{X}]-[\spc{Y}]|_{\GH}$;
but we will write it as $|\spc{X}\z-\spc{Y}|_{\GH}$ and say (not quite correctly) that 
\textit{$|\spc{X}\z-\spc{Y}|_{\GH}$ is the Gromov--Hausdorff distance from  $\spc{X}$ 
to  $\spc{Y}$}.
In other words, from now on the term \textit{metric space} might also stand for its \textit{isometry class}.

The metric on $\GH$ is defined as the maximal metric such that \textit{the distance between subspaces in a metric space is not greater than the Hausdorff distance between them}.
Here is a formal definition:

\begin{thm}{Definition}\label{def:GH}
The \index{Gromov--Hausdorff distance}\emph{Gromov--Hausdorff distance} $|\spc{X}-\spc{Y}|_{\GH}$ between compact metric spaces $\spc{X}$ and $\spc{Y}$
is defined by the following
relation.
 
Given  $r > 0$, we have that $|\spc{X}-\spc{Y}|_{\GH} < r$ if and only if there exists a metric
space $\spc{W}$ and subspaces $\spc{X}'$ and $\spc{Y}'$ in $\spc{W}$ that are isometric to $\spc{X}$ and $\spc{Y}$
respectively such that $|\spc{X}'-\spc{Y}'|_{\Haus\spc{W}} < r$. 
(Here $|\spc{X}'-\spc{Y}'|_{\Haus\spc{W}}$ denotes the Hausdorff distance between sets $\spc{X}'$ and $\spc{Y}'$ in $\spc{W}$.)
\end{thm}

\begin{thm}{Theorem}\label{thm:GH-is-a-metric}
The set of isometry classes of compact metric spaces equipped with Gromov--Hausdorff metric forms a metric space (which is denoted by $\GH$).

In other words, for arbitrary compact metric spaces $\spc{X}$, $\spc{Y}$, and $\spc{Z}$ the following conditions hold

\begin{subthm}{GH-1} $|\spc{X}-\spc{Y}|_{\GH}\ge 0$;
\end{subthm}

\begin{subthm}{GH-2} $|\spc{X}-\spc{Y}|_{\GH}=0$ if and only if $\spc{X}$ is isometric to $\spc{Y}$;
\end{subthm}

\begin{subthm}{GH-3} $|\spc{X}-\spc{Y}|_{\GH}=|\spc{Y}-\spc{X}|_{\GH}$;
\end{subthm}

\begin{subthm}{GH-4} $|\spc{X}-\spc{Y}|_{\GH}+|\spc{Y}-\spc{Z}|_{\GH}\ge |\spc{X}-\spc{Z}|_{\GH}$.
\end{subthm}
\end{thm}

Note that \ref{SHORT.GH-1}, \ref{SHORT.GH-3},
and the if part of \ref{SHORT.GH-2} follow directly from \ref{def:GH}.
Part \ref{SHORT.GH-4} will be proved in Section~\ref{sec:GH-approx}.
The only-if part of \ref{SHORT.GH-2} will be proved in Section~\ref{sec:extfun=GH}.

Recall that $a\cdot\spc{X}$ denotes $\spc{X}$ \index{rescaled space}\emph{rescaled} by a factor $a>0$;
that is, $a\cdot\spc{X}$ is a metric space with the underlying set of $\spc{X}$ and the metric defined by
\[\dist{x}{y}{a\cdot\spc{X}}\df a\cdot\dist{x}{y}{\spc{X}}.\]

\begin{thm}{Exercise}\label{ex:d_GH-and-diam}
Let $\spc{X}$ be a compact metric space,
$\spc{O}$ be the one-point metric space.
Prove the following.

\begin{subthm}{ex:d_GH-and-diam:point}
$|\spc{X}-\spc{O}|_{\GH}=\tfrac12\cdot \diam \spc{X}.$
\end{subthm}

\begin{subthm}{ex:d_GH-and-diam:scale}
$|a\cdot\spc{X}-b\cdot \spc{X}|_{\GH}=\tfrac12\cdot|a-b|\cdot\diam\spc{X}.$
\end{subthm}

\begin{subthm}{ex:d_GH-and-diam:isometry}
$\iota[\spc{O}]=[\spc{O}]$ for any isometry $\iota\:\GH\to\GH$.
\end{subthm}

\end{thm}

\begin{thm}{Exercise}\label{ex:GH<H}
Find two subsets $A,B\subset\RR^2$ such that 
\[|A-B|_{\GH}>|A-\iota(B)|_{\Haus\RR^2}\]
for any isometry $\iota$ of $\RR^2$.
\end{thm}

\begin{thm}{Exercise}\label{ex:rectangle}
Let $\spc{A}_r$ be a rectangle $1$ by $r$ in the euclidean plane 
and $\spc{B}_r$ be a closed line interval of length $r$.
Show that 
\[|\spc{A}_r-\spc{B}_r|_{\GH}>\tfrac1{10}\]
for all large $r$.
\end{thm}

\begin{thm}{Advanced exercise}\label{ex:GH-inj}
Let $\spc{X}$ and $\spc{Y}$ be compact metric spaces;
denote by $\hat{\spc{X}}$ and $\hat{\spc{Y}}$ their injective envelopes (see \ref{sec:extremal-functions}).
Show that 
\[|\hat{\spc{X}}-\hat{\spc{Y}}|_{\GH}\le 2\cdot|\spc{X}- \spc{Y}|_{\GH}.\] 
In other words, $\spc{X}\mapsto \hat{\spc{X}}$ defines a $2$-Lipschitz map $\GH\to\GH$.

\end{thm}

\section{Approximations and almost isometries}\label{sec:GH-approx}

\begin{thm}{Definition}\label{ex:defGHR}
Let $\spc{X}$ and $\spc{Y}$ be two metric spaces.
A relation $\approx$ between points in $\spc{X}$ and $\spc{Y}$ is called \index{$\eps$-approximation}\emph{$\eps$-approximation} if the following conditions hold:
\begin{itemize}
\item For any $x\in  \spc{X}$ there is $y\in \spc{Y}$ such that $x\approx y$.
\item For any $y\in  \spc{Y}$ there is $x\in \spc{X}$ such that $x\approx y$.
\item If $x\approx y$ and $x'\approx y'$ for some $x, x'\in  \spc{X}$ and $y,y'\in \spc{Y}$, then 
\[\dist{x}{x'}{\spc{X}}\lg\dist{y}{y'}{\spc{Y}}\bigr|\pm2\cdot\eps.\]
\end{itemize}

\end{thm}

\begin{thm}{Exercise}\label{ex:H-R}
Let $\spc{X}$ and $\spc{Y}$ be two compact metric spaces.
Show that
\[\dist{\spc{X}}{\spc{Y}}{\GH}<\eps\]
if and only if there is an $\eps$-approximation between $\spc{X}$ and $\spc{Y}$.

In other words, $\dist{\spc{X}}{\spc{Y}}{\GH}$ is the greatest lower bound of values $\eps>0$ such that  there is an $\eps$-approximation between $\spc{X}$ and $\spc{Y}$.
\end{thm}

\parit{Proof of \ref{GH-4}.}
Suppose that 
\begin{itemize}
\item $\approx_1$ is a relation between points in $\spc{X}$ and $\spc{Y}$,
\item $\approx_2$ is a relation between points in $\spc{Y}$ and $\spc{Z}$.
\end{itemize}
Consider the relation $\approx_3$ between points in $\spc{X}$ and $\spc{Z}$ such that
$x\approx_3 z$ if and only if there is $y\in  \spc{Y}$ such that 
$x\approx_1 y$ and $y\approx_2 z$.

It is straightforward to check that if $\approx_1$ is an $\eps_1$-approximation and $\approx_2$ is an $\eps_2$-approximation, then $\approx_3$ is an $(\eps_1+\eps_2)$-approximation.

Applying \ref{ex:H-R}, we get that if 
\[|\spc{X}-\spc{Y}|_{\GH}<\eps_1
\quad\text{and}\quad
|\spc{Y}-\spc{Z}|_{\GH}<\eps_2,
\]
then 
\[|\spc{X}-\spc{Z}|_{\GH}<\eps_1+\eps_2.\]
Hence \ref{GH-4} follows.
\qeds

The following weakened version of isometry is closely related to $\eps$-approximations.

\begin{thm}{Definition} Let $\spc{X}$ and $\spc{Y}$ be metric spaces and $\eps>0$. 
A  map\footnote{possibly noncontinuous} $f\: \spc{X} \z\to \spc{Y}$ is called an \index{almost isometry}\emph{$\eps$-isometry} 
if $f(\spc{X})$ is an $\eps$-net in $\spc{Y}$ and
\[\dist{x}{x'}{\spc{X}}\lg\dist{f(x)}{f(x')}{\spc{Y}}\bigr|\pm\eps\]
for any $x,x'\in \spc{X}$.
\end{thm}

\begin{thm}{Exercise}\label{ex:eps-isom}
Let $\spc{X}$ and $\spc{Y}$ be compact metric spaces.

\begin{subthm}{ex:eps-isom:GH>isom}
If $\dist{\spc{X}}{\spc{Y}}{\GH}<\eps$, then there is a $2\cdot\eps$-isometry $f\:\spc{X}\to\spc{Y}$.
\end{subthm}

\begin{subthm}{ex:eps-isom:isom>GH}
If there is an $\eps$-isometry $f\:\spc{X}\to\spc{Y}$, then $\dist{\spc{X}}{\spc{Y}}{\GH}<\eps$.
\end{subthm}

\end{thm}

\section{Optimal realization}\label{sec:extfun=GH}

Note that
\[\dist{\spc{X}'}{\spc{Y}'}{\Haus\spc{W}}\ge \dist{\spc{X}}{\spc{Y}}{\GH},\]
where $\spc{X}$, $\spc{Y}$, $\spc{X}'$, $\spc{Y}'$, and $\spc{W}$ are as in \ref{def:GH}.
The following proposition states that equality holds for some choice of $\spc{X}'$, $\spc{Y}'$, and $\spc{W}$.

\begin{thm}{Proposition}\label{prop:GH=H}
For any two compact metric spaces $\spc{X}$ and $\spc{Y}$ there is a metric space $\spc{W}$
with subsets $\spc{X}'$ and $\spc{Y}'$ such that 
$\spc{X}'\iso\spc{X}$, $\spc{Y}'\iso\spc{Y}$, and 
\[\dist{\spc{X}'}{\spc{Y}'}{\Haus\spc{W}}=\dist{\spc{X}}{\spc{Y}}{\GH}.\]
\end{thm}

Let us introduce the so-called \textit{appropriate functions} and use them in a reinterpretation of the Gromov--Hausdorff distance.

Suppose $\spc{X}$, $\spc{Y}$, $\spc{X}'$, $\spc{Y}'$, and $\spc{W}$ are as in \ref{def:GH}.
By passing to the subspace $\spc{X}'\cup\spc{Y}'$ in $\spc{W}$, we can assume that $\spc{W}=\spc{X}'\cup\spc{Y}'$.
Note that in this case the metric on $\spc{W}$ is completely determined by the function $f\:\spc{X}\times \spc{Y}\to\RR$ defined by
\[f(x,y)
\df
\dist{x}{y}{\spc{W}};\]
a function $f$ that can appear this way will be called \index{appropriate function}\emph{appropriate}.

Note that a function $f\:\spc{X}\times\spc{Y}\to\RR$ is appropriate if and only if
$x\mapsto f(x,y)$ and $y\mapsto f(x,y)$ are extension functions [see \ref{sec:Extension property}];
that is, if
\[
\begin{aligned}
f(x,y)+f(x,y')
&\ge \dist{y}{y'}{\spc{Y}}\ge |f(x,y)-f(x,y')|,\quad\text{and}
\\
f(x,y)+f(x',y)
&\ge \dist{x}{x'}{\spc{X}}\ge |f(x,y)-f(x',y)|
\end{aligned}
\eqlbl{eq:appropriate}
\]
for any $x,x',\in\spc{X}$ and  $y,y'\in\spc{X}$.
In other words, the following defines a semimetric on $\spc{X}\sqcup\spc{Y}$
\[\dist{x}{y}{\spc{X}\sqcup\spc{Y}}\df
\begin{cases}
\dist{x}{y}{\spc{X}}&\text{if\ } x,y\in \spc{X},
\\
\dist{x}{y}{\spc{Y}}&\text{if\ } x,y\in \spc{Y},
\\
f(x,y)&\text{if\ } x\in \spc{X}\ \text{and}\ y\in \spc{Y},
\end{cases}
\]
and the corresponding metric space $\spc{W}$ contains isometric copies of $\spc{X}$ and $\spc{Y}$.

\begin{thm}{Observation}\label{obs:GH=min-appropriate}
Let $\spc{X}$, $\spc{Y}$ be metric spaces.
Given an appropriate function $f\:\spc{X}\times\spc{Y}\to\RR$, set 
\begin{align*}
a_f&=\max_{x\in \spc{X}}\{\min_{y\in\spc{Y}} \{f(x,y)\}\},
\\
b_f&=\max_{y\in \spc{Y}}\{\min_{x\in\spc{X}} \{f(x,y)\}\},
\\
c_f&=\max\{\,a_f,b_f\,\}.
\end{align*}
Then 
\[\dist{\spc{X}}{\spc{Y}}{\GH}=\inf\{\,c_f\,\},\]
where the greatest lower bound is taken for all appropriate functions $f\:\spc{X}\times\spc{Y}\to\RR$.
\end{thm}

\parit{Proof of \ref{prop:GH=H}.}
Equip the product $\spc{X}\times\spc{Y}$ with $\ell_1$-metric;
that is,
\[\dist{(x,y)}{(x',y')}{\spc{X}\times\spc{Y}}
\df
\dist{x}{x'}{\spc{X}}+\dist{y}{y'}{\spc{Y}}\]
Note that any appropriate functions $f\:\spc{X}\times\spc{Y}\to\RR$ is $1$-Lipschitz.

Let us equip the space of appropriate functions $\spc{X}\times\spc{Y}\to\RR$ with supnorm.
Observe that the functional $f\mapsto c_f$ is continuous.
By the Arzelà--Ascoli theorem, we can choose an appropriate function $f$ 
with minimal possible value $c_f$.
It remains to apply \ref{obs:GH=min-appropriate}.
\qeds

\begin{thm}{Exercise}\label{ex:XYZ}
Construct three compact metric spaces $\spc{X}$, $\spc{Y}$, and $\spc{Z}$
such that for any metric space $\spc{W}$
with subsets $\spc{X}'$, $\spc{Y}'$, and $\spc{Z}'$ such that 
$\spc{X}'\iso\spc{X}$, $\spc{Y}'\iso\spc{Y}$, and $\spc{Z}'\iso\spc{Z}$
at least one of the following three inequalities is strict
\begin{align*}
\dist{\spc{X}'}{\spc{Y}'}{\Haus\spc{W}}&\ge \dist{\spc{X}}{\spc{Y}}{\GH},
\\
\dist{\spc{Y}'}{\spc{Z}'}{\Haus\spc{W}}&\ge\dist{\spc{Y}}{\spc{Z}}{\GH},
\\
\dist{\spc{Z}'}{\spc{X}'}{\Haus\spc{W}}&\ge\dist{\spc{Z}}{\spc{X}}{\GH}.
\end{align*}
\end{thm}

\section{Convergence}

The Gromov--Hausdorff metric defines Gromov--Hausdorff \index{Gromov--Hausdorff convergence}\emph{convergence}.
Namely, a sequence of compact metric spaces $\spc{X}_n$ converges to compact metric spaces $\spc{X}_\infty$ in the sense of Gromov--Hausdorff if 
\[\dist{\spc{X}_n}{\spc{X}_\infty}{\GH}\to 0\quad\text{as}\quad n\to\infty.\]

This convergence is more important than the metric ---
in all applications, we use only the topology on $\GH$
and we do not care about the particular value of the Gromov--Hausdorff distance between spaces.
The following observation follows from \ref{ex:eps-isom}:

\begin{thm}{Observation}\label{obs:GH-e-isom}
A sequence of compact metric spaces $(\spc{X}_n)$ converges to  $\spc{X}_\infty$ in the sense of Gromov--Hausdorff if and only if there is a sequence $\eps_n\to0+$
and an $\eps_n$-isometry $f_n\:\spc{X}_n\to \spc{X}_\infty$ for each $n$.
\end{thm}

\pagebreak

\begin{thm}{Exercise}\label{ex:GH-SC}
\begin{subthm}{ex:GH-SC:circle}
Show that a circle is not a Gromov--Hausdorff limit of compact simply-connected length spaces.
\end{subthm}

\begin{subthm}{ex:GH-SC:nonsc-limit}
Construct a compact non-simply-connected metric space
that is a Gromov--Hausdorff limit of compact simply-connected length spaces.
\end{subthm}
\end{thm}

\begin{thm}{Exercise}\label{ex:sphere-to-ball}
\begin{subthm}{ex:sphere-to-ball:2}
Show that a sequence of length metrics on the 2-sphere cannot converge to the unit disk in the sense of Gromov--Hausdorff.
\end{subthm}

\begin{subthm}{ex:sphere-to-ball:3}
Construct a sequence of length metrics on the 3-sphere that converges to the unit 3-ball in the sense of Gromov--Hausdorff.
\end{subthm}

\end{thm}

\section{Uniformly totally bonded families}

\begin{thm}{Definition}\label{def:utb}
A family $\bm{Q}$ of (isometry classes) of compact metric spaces is called  \index{uniformly totally bonded family}\emph{uniformly totally bonded} if it meets the following two conditions:

\begin{subthm}{}
spaces in $\bm{Q}$ have uniformly bounded diameters; that is, there is $D\in\RR$ such that
\[\diam\spc{X}\le D\]
for any space $\spc{X}$ in $\bm{Q}$.
\end{subthm}

\begin{subthm}{}
For any $\eps>0$ there is $n\in\NN$ such that any space $\spc{X}$ in $\bm{Q}$ admits an $\eps$-net with at most $n$ points.
\end{subthm}
\end{thm}

\begin{thm}{Exercise}\label{ex:utb+pack}
Let $\bm{Q}$ be a family of compact spaces with uniformly bounded diameters.
Show that $\bm{Q}$ is uniformly totally bonded if for any $\eps>0$ there is $n\in\NN$ such that 
\[\pack_\eps\spc{X}\le n\]
for any space $\spc{X}$ in $\bm{Q}$.
\end{thm}

Fix a real constant $C$.
A Borel measure $\mu$ on a metric space $\spc{X}$ is called \index{doubling space}\emph{$C$-doubling} if
\[\mu[\oBall(p,2\cdot r)]< C\cdot\mu[\oBall(p,r)]\]
for any point $p\in \spc{X}$ and any $r>0$.
A Borel measure is called \index{doubling measure}\emph{doubling} if it is {}\emph{$C$-doubling} for some real constant $C$.

\begin{thm}{Exercise}\label{pr:doubling}
Let $\bm{Q}(C,D)$ be the set of all the compact metric spaces with diameter at most $D$ that admit a $C$-doubling measure.
Show that $\bm{Q}(C,D)$ is totally bounded.
\end{thm}

Given two metric spaces $\spc{X}$ and $\spc{Y}$, we will write $\spc{X}\le \spc{Y}$ if there is a distance-noncontracting map $f\:\spc{X}\to \spc{Y}$;
that is, if 
$$ |x-x'|_{\spc{X}}\le|f(x)-f(x')|_{\spc{Y}}$$
for any $x,x'\in \spc{X}$.

\begin{thm}{Exercise}\label{pr:under}

\begin{subthm}{pr:under:if}
Let $\spc{Y}$ be a compact metric space.
Show that the set of all spaces $\spc{X}$ such that $\spc{X}\le\spc{Y}$
is uniformly totally bounded.
\end{subthm}

\begin{subthm}{pr:under:only-if}
Show that for any uniformly totally bounded set $\bm{Q}\subset\GH$ there is a compact space $\spc{Y}$
such that $\spc{X}\le\spc{Y}$ for any $\spc{X}$ in $\bm{Q}$.
\end{subthm}

\end{thm}

\section{Gromov selection theorem}

The following theorem is analogous to Blaschke selection theorems (\ref{thm:compact+Hausdorff}).

\begin{thm}{Gromov selection theorem}\label{thm:gromov-compactness}
Let $\bm{Q}$ be a closed subset of $\GH$.
Then $\bm{Q}$ is compact if and only if the spaces in $\bm{Q}$ are uniformly totally bounded.
\end{thm}

\begin{thm}{Lemma}\label{lem:GH-complete}
The space $\GH$ is complete.
\end{thm}

Suppose 
$\spc{U}$ and $\spc{V}$ are metric spaces 
with isometric closed sets $A\subset\spc{U}$ and $A'\subset\spc{V}$;
let $\iota\:A\to A'$ be an isometry.
Consider the gluing $\spc{W}=\spc{U}\sqcup_\iota\spc{V}$ of $\spc{U}$ and $\spc{V}$ along $\iota$ [see \ref{sec:max+glue}].

Let us identify points of $\spc{U}$ and $\spc{V}$ with their images in $\spc{W}$.
It is straightforward to check that the metric on~$\spc{W}$ is defined by
\begin{align*}
\dist{u}{u'}{\spc{W}}&\df\dist{u}{u'}{\spc{U}},
\\
\dist{v}{v'}{\spc{W}}&\df\dist{v}{v'}{\spc{V}},
\\
\dist{u}{v}{\spc{W}}&\df\min\set{\dist{u}{a}{\spc{U}}+\dist{v}{\iota(a)}{\spc{V}}}{a\in A},
\end{align*}
where $u,u'\in \spc{U}$ and $v,v'\in \spc{V}$.

If one applies this construction to two copies of one space $\spc{U}$ with a set $A\subset \spc{U}$ and the identity map $\iota\:A\to A$, then the obtained space is called the \index{doubling}\emph{doubling} of $\spc{U}$ along~$A$; this space can be denoted by $\sqcup_A^2\spc{U}$.

Note that the inclusions $\spc{U}\hookrightarrow \spc{W}$ and $\spc{V}\hookrightarrow \spc{W}$ are distance-preserving.
Therefore we can and will consider $\spc{U}$ and $\spc{V}$ as the subspaces of $\spc{W}$;
this way the subsets $A$ and $A'$ will be identified and denoted further by~$A$.
Note that $A=\spc{U}\cap \spc{V}\subset \spc{W}$.

\parit{Proof.}
Let $\spc{X}_1,\spc{X}_2,\dots$ be a Cauchy sequence in $\GH$.
Passing to a subsequence if necessary, 
we can assume that $|\spc{X}_n-\spc{X}_{n+1}|_{\GH}<\tfrac1{2^n}$ for each~$n$.
In particular, for each $n$ there is a metric space $\spc{V}_n$ with distance-preserving inclusions $\spc{X}_n\hookrightarrow \spc{V}_n$ and $\spc{X}_{n+1}\hookrightarrow \spc{V}_n$ such that
\[|\spc{X}_n-\spc{X}_{n+1}|_{\Haus\spc{V}_n}<\tfrac1{2^n}\]
for each $n$.
Moreover, we may assume that $\spc{V}_n=\spc{X}_n\cup\spc{X}_{n+1}$.

Let us glue $\spc{V}_1$ to $\spc{V}_2$ along $\spc{X}_2$;
to the obtained space glue $\spc{V}_3$ along $\spc{X}_3$, and so on.
The obtained metric space $\spc{W}$
has an underlying set formed by the disjoint union of all $\spc{X}_n$ such that each inclusion $\spc{X}_n\z\hookrightarrow\spc{W}$ is distance-preserving and
\[|\spc{X}_n-\spc{X}_{n+1}|_{\Haus\spc{W}}<\tfrac1{2^n}\]
for each $n$.
In particular,
\[|\spc{X}_m-\spc{X}_n|_{\Haus\spc{W}}<\tfrac1{2^{n-1}}\eqlbl{eq:|x_m-X_n|}\] 
if $m>n$.

Denote by $\bar{\spc{W}}$ the completion of $\spc{W}$.
Observe that the union $\spc{X}_1\z\cup \spc{X}_2\cup\z\dots\cup \spc{X}_n$ is compact and \ref{eq:|x_m-X_n|} implies that it forms a $\tfrac1{2^{n-1}}$-net in $\bar{\spc{W}}$.
Whence $\bar{\spc{W}}$ is compact; see \ref{totally-bounded} and \ref{ex:compact-net}.

Applying the Blaschke selection theorem (\ref{thm:compact+Hausdorff}),
we can pass to a subsequence of $\spc{X}_n$ that converges in $\Haus\bar{\spc{W}}$; denote its limit by $\spc{X}_\infty$.
It remains to observe that $\spc{X}_\infty$ is the Gromov--Hausdorff limit of $\spc{X}_n$.
\qeds

\parit{Proof of \ref{thm:gromov-compactness}; only-if part.}
Suppose that there is no sequence $\eps_n\to0$ as described in \ref{def:utb}.
Observe that in this case
there is a sequence of spaces $\spc{X}_n\in\bm{Q}$ such that 
$$\pack_\delta \spc{X}_n\to\infty
\quad\text{as}\quad
n\to\infty$$
for some fixed $\delta>0$.

Since $\bm{Q}$ is compact, 
this sequence has a partial limit, say $\spc{X}_\infty\in\bm{Q}$.
Observe that $\pack_{\delta} \spc{X}_\infty=\infty$.
Therefore, $\spc{X}_\infty$ is not compact --- a contradiction.

\parit{If part.}
Given a positive integer $n$ consider the set of all nonempty metric spaces $\spc{W}_n$
with the number of points at most $n$ and diameter $\le D$.
Note that $\spc{W}_n$ is a compact set in $\GH$ for each $n$.

Let $D$ and $n=n(\eps)$ be as in the definition of uniformly totally bonded families (\ref{def:utb}).

Note that an $\eps$-net of any $\spc{X}\in\bm{Q}$ belongs to $\spc{W}_{n(\eps)}$.
Therefore, $\spc{W}_{n(\eps)}$ is a compact $\eps$-net of $\bm{Q}$ for any $\eps>0$.
Since $\bm{Q}$ is closed in a complete space $\GH$, it implies that $\bm{Q}$ is compact.
\qeds

\begin{thm}{Exercise}\label{ex:GH-G-delta}
Show that most of the compact metric spaces are homeomorphic to the Cantor set.

More precisely, suppose $\bm{Q}$ denotes all metric spaces homeomorphic to the Cantor set.
Show that $\bm{Q}$ is a dense G-delta set in $\GH$.
\end{thm}

\begin{thm}{Exercise}\label{ex-GH-length}
Show that the space $\GH$ is 
\begin{multicols}{3}

\begin{subthm}{ex-GH-length:separable}
separable,
\end{subthm}

\begin{subthm}{ex-GH-length:length}
length, and
\end{subthm}

\begin{subthm}{ex-GH-length:geodesic}
geodesic.
\end{subthm}

\end{multicols}

\end{thm}

\begin{thm}{Exercise}\label{ex:GH-po}
For two metric spaces $\spc{X}$ and $\spc{Y}$,
we write $\spc{X}\le \spc{Y}+\eps$ if
there is a map $f\:\spc{X}\to \spc{Y}$ such that 
\[\dist{x}{x'}{\spc{X}}\le \dist{f(x)}{f(x')}{\spc{Y}}+\eps\]
for any $x,x'\in \spc{X}$.

\begin{subthm}{ex:GH-po:a}
Show that 
$$\dist{\spc{X}}{\spc{Y}}{\GH'}
\df
\inf\set{\eps>0}{\spc{X}\le \spc{Y}+\eps
\quad\text{and}\quad
\spc{Y}\le \spc{X}+\eps}$$
defines a metric on the space of (isometry classes) of compact metric spaces.
\end{subthm}

\begin{subthm}{ex:GH-po:b}
Moreover, $\dist{*}{*}{\GH'}$ is equivalent to the Gromov--Hausdorff metric;
that is,
$$|\spc{X}_n-\spc{X}_\infty|_{\GH}\to 0 
\quad\iff\quad 
\dist{\spc{X}_n}{\spc{X}_\infty}{\GH'}\to 0$$ 
as $n\to\infty$.
\end{subthm}
\end{thm}

\section{Universal ambient space}

Recall that a metric space is called universal if it contains an isometric copy of any separable metric space (in particular, any compact metric space).
Examples of universal spaces include $\spc{U}_\infty$ --- the Urysohn space and $\ell^\infty$ --- the space of bounded infinite sequences with the metric defined by $\sup$-norm; see \ref{prop:sep-in-urys} and \ref{ex:frechet}.

The following proposition says that the space $\spc{W}$ in Definition~\ref{def:GH} can be exchanged to a fixed universal space.

\begin{thm}{Proposition}\label{prop:GH-with-fixed-Z}
Let $\spc{U}$ be a universal space.
Then for any compact metric spaces $\spc{X}$ and $\spc{Y}$ we have
$$|\spc{X}-\spc{Y}|_{\GH} = \inf \{|\spc{X}'-\spc{Y}'|_{\Haus\spc{U}}\},$$ 
where the greatest lower bound is taken over all pairs of sets $\spc{X}'$ and $\spc{Y}'$ in $\spc{U}$
which isometric to  $\spc{X}$ and $\spc{Y}$ respectively.  
\end{thm}

\parit{Proof of \ref{prop:GH-with-fixed-Z}.}
By the definition (\ref{def:GH}), we have that 
\[|\spc{X}-\spc{Y}|_{\GH} \le \inf \{|\spc{X}'-\spc{Y}'|_{\Haus\spc{U}}\};\]
it remains to prove the opposite inequality.

Suppse $|\spc{X}-\spc{Y}|_{\GH}<\eps$;
let $\spc{X}'$, $\spc{Y}'$ and $\spc{W}$ be as in \ref{def:GH}.
We can assume that $\spc{W}=\spc{X}'\cup\spc{Y}'$;
otherwise, pass to the subspace $\spc{X}'\cup\spc{Y}'$ of~$\spc{W}$.
In this case, $\spc{W}$ is compact;
in particular, it is separable.

Since $\spc{U}$ is universal, there is a distance-preserving embedding of $\spc{W}$ in $\spc{U}$;
let us keep the same notation for $\spc{X}'$, $\spc{Y}'$, and their images.
It follows that 
\[|\spc{X}'-\spc{Y}'|_{\Haus\spc{U}}<\eps,\]
--- hence the result.
\qeds

\begin{thm}{Exercise}\label{ex:GH-urysohn}
Let $\spc{U}_\infty$ be the Urysohn space.
Given two compact sets $A$ and $B$ in $\spc{U}_\infty$, define 
\[\|A-B\|\df\inf\{|A-\iota(B)|_{\Haus\spc{U}_\infty}\},\]
where the greatest lower bound is taken for all isometrics $\iota$ of $\spc{U}_\infty$.
Show that $\|{*}\z-{*}\|$ defines a semimetric 
on nonempty compact subsets of $\spc{U}_\infty$ and its corresponding metric space is isometric to $\GH$.
\end{thm}

The value $\|A-B\|$ is called Hausdorff distance \index{Hausdorff distance!up to isometry}\emph{up to isometry} from $A$ to $B$ in $\spc{U}_\infty$.

\section{Remarks}\label{sec:remarks-GH}

Suppose $\spc{X}_n\GHto \spc{X}_\infty$, then there is a metric on the disjoint union 
\[\bm{X}=\bigsqcup_{n\in \NN\cup\{\infty\}} \spc{X}_n\] 
that satisfies the following property:

\begin{thm}{Property}\label{propery:GH}
The restriction of metric on each $\spc{X}_n$ and $\spc{X}_\infty$ coincides with its original metric, 
and $\spc{X}_n\Hto \spc{X}_\infty$ as subsets in $\bm{X}$.
\end{thm}

Indeed, since $\spc{X}_n\GHto \spc{X}_\infty$, there is a metric on $\spc{V}_n=\spc{X}_n\sqcup \spc{X}_\infty$ such that the restriction of metric on each $\spc{X}_n$ and $\spc{X}_\infty$ coincides with its original metric, and $\dist{\spc{X}_n}{\spc{X}_\infty}{\Haus\spc{V}_n}<\eps_n$ for some sequence $\eps_n\to 0$.
Gluing all $\spc{V}_n$ along $\spc{X}_\infty$, we get the required space $\bm{X}$.

In other words, the metric on $\bm{X}$ \textit{defines} the convergence $\spc{X}_n\z\GHto \spc{X}_\infty$.
This metric makes it possible to talk about limits of sequences $x_n\in \spc{X}_n$ as $n\to\infty$, as well as weak limits of a sequence of Borel measures $\mu_n$ on $\spc{X}_n$ and so on.

For that reason, it is useful to define \index{Hausdorff convergence}\emph{convergence} by specifying the metric on $\bm{X}$ that satisfies the property
for the variation of Hausdorff convergence described in Section~\ref{sec:H-variation}.

This approach is more flexible;
in particular, it can be used to define the Gromov--Hausdorff convergence of arbitrary metric spaces (not necessarily compact).
A limit space for this generalized convergence is not uniquely defined.
For example, if each space $\spc{X}_n$ in the sequence is isometric to the half-line, then its limit might be isometric to the half-line or the whole line.
The first convergence is evident and the second could be guessed from the diagram.

\begin{figure}[ht!]
\vskip-0mm
\centering
\includegraphics{mppics/pic-500}
\end{figure}

Often the isometry class of the limit can be fixed by marking a point $p_n$ in each space $\spc{X}_n$, it is called \index{pointed convergence}\emph{pointed Gromov--Hausdorff convergence} --- we say that $(\spc{X}_n,p_n)$ converges to $(\spc{X}_\infty,p_\infty)$ if there is a metric on $\bm{X}$ as in \ref{propery:GH} such that $\spc{X}_n\Hto \spc{X}_\infty$ and $p_n\to p_\infty$.
For example, the sequence $(\spc{X}_n,p_n)=(\RR_+,0)$ converges to $(\RR_+,0)$, while $(\spc{X}_n,p_n)=(\RR_+,n)$ converges to $(\RR,0)$.

The pointed convergence works nicely for proper metric spaces;
the following theorem is an analog of Gromov's selection theorem for this convergence.

\begin{thm}{Theorem}\label{thm:pointed-gromov-compactness}%
Let $\bm{Q}$ be a set of isometry classes of pointed proper metric spaces.
Assume that for any $R>0$, the $R$-balls in the spaces centered at the marked points form a uniformly totally bounded family of spaces.
Then $\bm{Q}$ is precompact with respect to the pointed Gromov--Hausdorff convergence. 
\end{thm}

%% file: ultralimit.tex
\chapter{Ultralimits}\label{chap:ultralimits}

Ultralimits provide a very general way to pass to a limit.
This procedure works for \textit{any} sequence of metric spaces, its result reminds limit in the sense of Gromov--Hausdorff, but has some strange features; for example, the limit of a constant sequence of spaces $\spc{X}_n=\spc{X}$ is \textit{not} $\spc{X}$ in general (see \ref{ex:ultrapower:compact}).

In geometry, ultralimits are used mostly as a canonical way to pass to a convergent subsequence.
It is very useful in the proofs where one needs to repeat ``pass to convergent subsequence'' too many times.

This lecture is based on the introductory part of the paper by Bruce Kleiner and Bernhard Leeb \cite{kleiner-leeb}.

\section{Faces of ultrafilters}

\parbf{Measure-theoretic definition.}
Recall that $\NN=\{1,2,\dots\}$ is the set of natural numbers.

\begin{thm}{Definition}
A finitely additive measure $\omega$ 
on $\NN$ 
is called an \index{ultrafilter}\emph{ultrafilter} if it meets the following condition:
\begin{subthm}{}
$\omega(\NN)=1$ and 
$\omega(S)=0$ or $1$ for any subset $S\subset \NN$.
\end{subthm}
An ultrafilter $\omega$ is called 
\emph{nonprincipal}\index{ultrafilter!nonprincipal ultrafilter}\index{nonprincipal ultrafilter} if in addition 
\begin{subthm}{}
$\omega(F)=0$ for any finite subset $F\subset \NN$.
\end{subthm}
\end{thm}

If $\omega(S)=0$ for some subset $S\subset \NN$,
we say that $S$ is \index{$\omega$-small}\emph{$\omega$-small}. 
If $\omega(S)=1$, we say that $S$ contains \index{$\omega$-almost all}\emph{$\omega$-almost all} elements of $\NN$.

\begin{thm}{Advanced exercise}\label{ex:ultrakatetov}
Let $\omega$ be an ultrafilter on $\NN$ and $f\:\NN\z\to \NN$.
Suppose that $\omega(S)\le \omega(f^{-1}(S))$ for any set $S\subset \NN$.
Show that $f(n)=n$ for $\omega$-almost all $n\in\NN$.
\end{thm}

\parbf{Classical definition.}
More commonly, a nonprincipal ultrafilter is defined as a collection, say $\mathfrak{F}$, of subsets in $\NN$ such that
\begin{enumerate}
\item\label{filter:supset} if $P\in \mathfrak{F}$ and $Q\supset P$, then $Q\in \mathfrak{F}$,
\item\label{filter:cap} if $P, Q\in \mathfrak{F}$, then $P\cap Q\in \mathfrak{F}$,
\item\label{filter:ultra} for any subset $P\subset\NN$, either $P$ or its complement is an element of $\mathfrak{F}$.
\item\label{filter:non-prin} if $F\subset \NN $ is finite, then $F\notin \mathfrak{F}$.
\end{enumerate}
Setting $P\in\mathfrak{F}\Leftrightarrow\omega(P)=1$ makes these two definitions equivalent.

A nonempty collection of sets $\mathfrak{F}$ that does not include the empty set and satisfies only conditions \ref{filter:supset} and \ref{filter:cap} is called a \index{filter}\emph{filter}; 
if in addition, $\mathfrak{F}$ satisfies condition \ref{filter:ultra} it is called an \index{ultrafilter}\emph{ultrafilter}.
From Zorn's lemma, it follows that every filter contains an ultrafilter.
Thus there is an ultrafilter $\mathfrak{F}$ contained in the filter of all complements of finite sets.
Clearly, this ultrafilter $\mathfrak{F}$ is nonprincipal.

\parbf{Stone--\v{C}ech compactification.}
Given a set $S\subset \NN$, consider subset $\Omega_S$ of all ultrafilters $\omega$ such that $\omega(S)=1$.
It is straightforward to check that the sets $\Omega_S$ for all subsets $S\subset \NN$ form a topology on the set of ultrafilters on $\NN$. 
The obtained space was first considered by Andrey Tikhonov and called \index{Stone--\v{C}ech compactification}\emph{Stone--\v{C}ech compactification} of $\NN$;
it is usually denoted as $\beta\NN$.

Let $\omega_n$ be a principal ultrafilter such that $\omega_n(\{n\})=1$; that is, $\omega_n(S)=1$ if and only if $n\in S$.
Note that $n\mapsto\omega_n$ defines an embedding $\NN\hookrightarrow\beta\NN$;
so, we can (and will) consider $\NN$ as a subset of~$\beta\NN$.

The space $\beta\NN$ is the maximal compact Hausdorff space that contains $\NN$  as an everywhere dense subset.
More precisely, the inclusion $\NN\hookrightarrow\beta\NN$ has the following universal property:
\textit{for any compact Hausdorff space $\spc{X}$ 
and a map $f\:\NN\to \spc{X}$ there is a unique continuous map $\bar f\:\beta\NN\to X$ such that the restriction $\bar f|_\NN$ coincides with $f$.} 

\section{Ultralimits of points}
\label{ultralimits}\index{ultralimit}

Let us fix a nonprincipal  ultrafilter $\omega$ once and for all.

Assume $x_n$ is a sequence of points in a metric space $\spc{X}$. 
Let us define the \index{$\omega$-limit}\emph{$\omega$-limit} of a sequence $x_1,x_2,\dots$ as the point $x_\omega\in\spc{X}$ 
such that for any $\eps>0$, point $x_n$ lies in $\oBall(x_\omega,\eps)$ for $\omega$-almost all $n$; 
that is, if 
\[S_\eps=\set{n\in\NN}{\dist{x_\omega}{x_n}{}<\eps},\]
then $\omega(S_\eps)=1$ for any $\eps>0$.
In this case, we will write 
\[x_\omega=\lim_{n\to\omega} x_n
\ \ \text{or}\ \ 
x_n\to x_\omega\ \text{as}\ n\to\omega.\]

For example, if $\omega_n$ is the \textit{principal} ultrafilter defined by $\omega_n\{n\}=1$ for some $n\in\NN$, then
$x_{\omega_n}=x_n$.

The sequence $x_n$ can be regarded as a map $\NN\to\spc{X}$ defined by $n\mapsto x_n$.
If $\spc{X}$ is compact, then the map $\NN\to\spc{X}$ can be extended to a continuous map $\beta\NN\to\spc{X}$ from the Stone--\v{C}ech compactification $\beta\NN$ of $\NN$.
Then the $\omega$-limit $x_\omega$ is the image of $\omega$.

Note that the $\omega$-limits of a sequence and its subsequence may differ.
For example, sequence $y_n=-(-1)^n$ is a subsequence of $x_n=(-1)^n$, but for any ultrafilter $\omega$, we have
\[\lim_{n\to\omega}x_n
\ne
\lim_{n\to\omega}y_n.\] 

\begin{thm}{Proposition}\label{prop:ultra/partial}
Let $x_n$ be a sequence of points in a metric space $\spc{X}$.
Assume $x_n\to x_\omega$ as $n\to\omega$.
Then $x_\omega$ is a \index{partial limit}\emph{partial limit} of $x_n$;
that is, there is a subsequence of $x_n$ that converges to $x_\omega$ in the usual sense.
\end{thm}

\parit{Proof.}
Given $\eps>0$, 
let $S_\eps=\set{n\in\NN}{\dist{x_n}{x_\omega}{}<\eps}$.
Recall that $\omega(S_\eps)=1$ for any $\eps>0$.

Since $\omega$ is nonprincipal, the set $S_\eps$ is infinite for any $\eps>0$.
Therefore, we can choose an increasing sequence $n_k$
such that $n_k\in S_{\frac1k}$ for each $k\in \NN$.
Clearly, $x_{n_k}\to x_\omega$ as $k\to\infty$.
\qeds

\begin{thm}{Proposition}\label{prop:ultra/compact}
Any sequence $x_n$ of points in a compact metric space $\spc{X}$ has a unique $\omega$-limit $x_\omega$.

In particular, a bounded sequence of real numbers has a unique $\omega$-limit.
\end{thm}

The proposition is  analogous to the Bolzano--Weierstrass theorem,
and it can be proved the same way.
The following lemma is an ultralimit analog of the Cauchy convergence test.

\begin{thm}{Lemma}\label{lem:X-X^w}
A sequence of points in a metric space converges if and only if all its subsequences 
have the same $\omega$-limit.
\end{thm}

\parit{Proof.} 
The only-if part is evident; it remains to prove the if part.
Suppose $z$ is a $\omega$-limit of all subsequences of $x_1$, $x_2,\dots$
By \ref{prop:ultra/partial}, $z$ is a partial limit of $x_n$.
If $x_1$, $x_2,\dots$ is Cauchy, then $x_n\to z$, and the lemma is proved.

Assume $x_1$, $x_2,\dots$ is not Cauchy.
Then for some $\eps>0$, there is a subsequence $y_n$ of $x_n$ such that $\dist{x_n}{y_n}{}\ge\eps$ for all $n$.
Therefore $\dist{x_\omega}{y_\omega}{}\ge \eps$ --- a contradiction. 
\qeds

Recall that $\ell^\infty$ denotes the space of bounded sequences of real numbers equipped with the sup-norm.

\begin{thm}{Exercise}\label{ex:linear}
Construct a linear functional $L\:\ell^\infty\to\RR$ such that
for any sequence $\bm{s}=(s_1,s_2,\dots)\in \ell^\infty$ the image $L(\bm{s})$ is a partial limit of $s_1,s_2,\dots$
\end{thm}

\begin{thm}{Exercise}\label{ex:ultrakatetov+}
Suppose that $f\:\NN\to\NN$ is a map such that 
\[\lim_{n\to\omega}x_n=\lim_{n\to\omega}x_{f(n)}\]
for any bounded sequence $x_n$ of real numbers.
Show that $f(n)=n$ for $\omega$-almost all $n\in\NN$.
\end{thm}

\section{An illustration}

In this section, we illustrate the power of ulralimits by proving the following simple claim.

\begin{thm}{Claim}
Let $\spc{X}$ and $\spc{Y}$ be compact spaces.
Suppose that for every $n\in\NN$ there is a $\tfrac1n$-isometry $f_n\:\spc{X}\to \spc{Y}$.
Then there is an isometry $\spc{X}\to \spc{Y}$.
\end{thm}

\parit{Proof.}
Consider the $\omega$-limit $f_\omega$ of~$f_n$;
according to \ref{prop:ultra/compact}, $f_\omega$ is defined.
Since 
\[|f_n(x)-f_n(x')|\lege |x-x'|\pm\tfrac1n\]
we get that 
\[|f_\omega(x)-f_\omega(x')|= |x-x'|\]
for any $x,x'\in \spc{X}$;
that is, $f_\omega$ is distance-preserving.

Further, since $f_n$ is a $\tfrac1n$-isometry,
for any $y\in \spc{Y}$ there is a sequence $x_1,x_2,\dots\in \spc{X}$ such that $|f_n(x_n)-y|\le \tfrac1n$ for any $n$.
Therefore,
\[f_\omega(x_\omega)=y,\]
where $x_\omega$ is the $\omega$-limit of $x_n$;
that is, $f_\omega$ is onto.

It follows that $f_\omega\:\spc{X}\to\spc{Y}$ is an isometry.
\qeds

\section{Ultralimits of spaces}\label{sec:Ultralimit of spaces}

Recall that $\omega$ is a fixed nonprincipal ultrafilter on $\NN$.

Let $\spc{X}_n$ be a sequence of metric spaces.
Consider all sequences of points $x_n\in \spc{X}_n$.
On the set of all such sequences, define an $\infty$-semimetric by
\[\dist{(x_n)}{(y_n)}{}
\df
\lim_{n\to\omega} \dist{x_n}{y_n}{\spc{X}_n}.
\eqlbl{eq:olim-dist}\]
Note that the $\omega$-limit on the right-hand side is always defined 
and takes a value in $[0,\infty]$. 
(The $\omega$-convergence to $\infty$ is defined analogously to the usual convergence to $\infty$; that is, $\lim x_n=\infty$ $\Longleftrightarrow$ $\lim\tfrac1{x_n}=0$).

Let $\spc{X}_\omega$ be the corresponding metric space; 
that is, the underlying set of $\spc{X}_\omega$ is formed by classes of equivalence of sequences of points $x_n\in\spc{X}_n$ 
defined by 
\[(x_n)\sim(y_n)
\ \Leftrightarrow\ 
\lim_{n\to\omega} \dist{x_n}{y_n}{}=0\]
and the distance is defined by \ref{eq:olim-dist}.

The space $\spc{X}_\omega$ is called the \index{$\omega$-limit space}\emph{$\omega$-limit} of $\spc{X}_n$.
(It is also called \index{$\omega$-product}\emph{$\omega$-product}; this term is motivated by the fact that   
$\spc{X}_\omega$ is a quotient of the product $\prod\spc{X}_n$)
Typically  $\spc{X}_\omega$ will denote the  
$\omega$-limit of sequence $\spc{X}_n$;
we may also write  
\[\spc{X}_n\to\spc{X}_\omega\ \ \text{as}\ \  n\to\omega\ \ \text{or}\ \ \spc{X}_\omega=\lim_{n\to\omega}\spc{X}_n.\]

Given a sequence $x_n\in \spc{X}_n$,
we will denote by $x_\omega$ its equivalence class which is a point in $\spc{X}_\omega$;
it can be written as
\[x_n\to x_\omega \ \ \text{as}\ \  n\to\omega,\ \ \text{or}\ \ x_\omega=\lim_{n\to\omega} x_n.\]

\begin{thm}{Observation}\label{obs:ultralimit-is-complete}
The $\omega$-limit of any sequence of metric spaces is complete. 
\end{thm}

We will repeat the proof of \ref{ex:complete-completion} using a slightly different language.

\parit{Proof.}
Let $\spc{X}_n$ be a sequence of metric spaces and $\spc{X}_n\to\spc{X}_\omega$ as $n\to\omega$.
Choose a Cauchy sequence $x_1,x_2,\dots{}\in\spc{X}_\omega$.
Passing to a subsequence, we can assume that $\dist{x_k}{x_{m}}{\spc{X}_\omega}<\tfrac1{k}$ if $k<m$.

Choose a double sequence $x_{n,m}\in \spc{X}_n$ such that for any fixed $m$ we have $x_{n,m}\to x_m$ as $n\to\omega$.
Note that for any $k<m$ the inequality $\dist{x_{n,k}}{x_{n,m}}{}<\tfrac1{k}$ holds for $\omega$-almost all $n$.

Given $m\in\NN$, consider the subset $S_m\subset\NN$ defined by
\[S_m=\set{n\ge m}{\dist{x_{n,k}}{x_{n,l}}{}<\tfrac1{k} \quad\text{for all}\quad k<l\le m}.\]
Note that 
\begin{itemize}
\item $\NN= S_1\supset S_2\supset\dots$
\item $\omega(S_m)=1$ for each $m$, and
\item $\min S_m\ge m$.
\end{itemize}

Consider the sequence $y_n=x_{n,m(n)}$, where $m(n)$ is the largest value such that $n\in S_{m(n)}$;
from above, $m(n)\le n$.
Denote by $y_\omega\in \spc{X}_\omega$ the $\omega$-limit of $y_n$.

Observe that $|y_m-x_{n,m}|<\tfrac1{m}$ for $\omega$-almost all $n$.
It follows that $|x_m-y_\omega|\le \tfrac1{m}$ for any $m$.
Therefore, $x_n\to y_\omega$ as $n\to \infty$.
That is, any Cauchy sequence in $\spc{X}_\omega$ converges.
\qeds

\begin{thm}{Observation}\label{obs:ultralimit-is-geodesic}
The $\omega$-limit of any sequence of length spaces is geodesic. 
\end{thm}

\parit{Proof.}
If $\spc{X}_n$ is a sequence of length spaces, then for any sequence of pairs $x_n, y_n\in X_n$ there is a sequence of $\tfrac1n$-midpoints $z_n$.

Let $x_n\to x_\omega$, $y_n\to y_\omega$ and $z_n\to z_\omega$ as $n\to \omega$.
Note that $z_\omega$ is a midpoint of $x_\omega$ and $y_\omega$ in $\spc{X}_\omega$.

By \ref{obs:ultralimit-is-complete}, $\spc{X}_\omega$ is complete.
Applying Menger's lemma (\ref{lem:mid>geod}), we get the statement.
\qeds

\begin{thm}{Exercise}\label{ex:lim(tree)}
Show that an ultralimit of metric trees is a metric tree. 
\end{thm}

\begin{thm}{Exercise}\label{ex:ultracompact}
Suppose that $\spc{X}_\infty$ and $\spc{X}_1,\spc{X}_2,\dots$ are compact metric spaces.
Assume $\spc{X}_n\GHto\spc{X}_\infty$.
Show that $\spc{X}_\omega\iso\spc{X}_\infty$.
\end{thm}

\parbf{Pointed limit.}
If $\diam \spc{X}_n\to \infty$ as $n\to\omega$,
then the metric on $\spc{X}_\omega$ takes value $\infty$;
so $\spc{X}_\omega$ has at least two metric components.

To specify a metric component in $\spc{X}_\omega$,
we may choose a sequence of marked points $p_n$ in $\spc{X}_n$, and pass to the metric component of its $\omega$-limit $p_\omega$ in $\spc{X}_\omega$.
The obtained metric component $\spc{Z}=(\spc{X}_\omega)_{p_\omega}$ with marked $p_\omega$ is called the \index{pointed ultralimi}\emph{pointed $\omega$-limit} of $(\spc{X}_n,p_n)$.
Note that $\spc{Z}$ is a genuine metric space.

If, in the definition of ultralimit, we consider only sequences $x_n\z\in \spc{X}_n$ with such that $\dist{p_n}{x_n}{}$ is bounded, 
then arrive at $\spc{Z}$.  

For proper metric spaces, there is a relation between the pointed ultralimit and the pointed Gromov--Hausdorff limit introduced in~\ref{sec:remarks-GH}.
Namely, \textit{if $(\spc{X}_\infty,p_\infty)$ and $(\spc{X}_1,p_1),(\spc{X}_2,p_2),\dots$ are proper pointed metric spaces such that $(\spc{X}_n,p_n)\GHto(\spc{X}_\infty,p_\infty)$ as defined in \ref{sec:remarks-GH}, then $(\spc{X}_\infty,p_\infty)$ is isometric to the pointed $\omega$-limit of $(\spc{X}_n,p_n)$.}
The proof is the same as for \ref{ex:ultracompact}.

\section{Ultrapower}

If all the metric spaces in the sequence are identical $\spc{X}_n=\spc{X}$, 
its $\omega$-limit 
$\lim_{n\to\omega}\spc{X}_n$
is denoted by $\spc{X}^\omega$
and called \index{ultrapower}\index{$\omega$-power}\emph{$\omega$-power} of $\spc{X}$ (also known as \index{ultracompletion}\index{$\omega$-completion}\emph{$\omega$-completion}).

\begin{thm}{Exercise}\label{ex:ultrapower}
For any point $x\in \spc{X}$, consider the constant sequence $x_n=x$
and set $\iota(x)=\lim_{n\to\omega}x_n\in\spc{X}^\omega$.

\begin{subthm}{ex:ultrapower:a}
Show that $\iota\:\spc{X}\to\spc{X}^\omega$ is distance-preserving embedding. (So we can and will consider $\spc{X}$ as a subset of $\spc{X}^\omega$.)
\end{subthm}

\begin{subthm}{ex:ultrapower:compact} 
Show that $\iota$ is onto if and only if $\spc{X}$ is compact.
\end{subthm}

\begin{subthm}{ex:ultrapower:proper} 
Show that if $\spc{X}$ is proper, then $\iota(\spc{X})$ forms a metric component of $\spc{X}^\omega$; that is, a subset of $\spc{X}^\omega$ that lies at a finite distance from a given point.
\end{subthm}

\end{thm}

If $\spc{X}$ is a genuine metric space, then the metric component of $\spc{X}$ in $\spc{X}^\omega$ is also called the ultrapower of $\spc{X}$;
if needed, we may call it the \index{small ultrapower}\emph{small ultrapower}, and the whole space $\spc{X}^\omega$ could be called the \index{big ultrapower}\emph{big ultrapower} of $\spc{X}$.
Note that the small ultrapower of genuine metric space is a genuine metric space.
Further, according to \ref{SHORT.ex:ultrapower:proper}, \textit{proper metric space is isometric to its small ultrapower}.

Note that \ref{SHORT.ex:ultrapower:compact} implies that the inclusion $\spc{X}\hookrightarrow\spc{X}^\omega$ is not onto if the space $\spc{X}$ is not compact.
However, the spaces $\spc{X}$ and $\spc{X}^\omega$ might be isometric; here is an example:

\begin{thm}{Exercise}\label{ex:isom-ultrapower}
Let $\spc{X}$ be an infinite countable set with \index{discrete metric}\emph{discrete metric};
that is $\dist{x}{y}{\spc{X}}=1$ if $x\ne y$.
Show that 

\begin{subthm}{ex:isom-ultrapower:no}
$\spc{X}^\omega$ is not isometric to $\spc{X}$, but
\end{subthm}

\begin{subthm}{ex:isom-ultrapower:yes}
$\spc{X}^\omega$ is  isometric to $(\spc{X}^\omega)^\omega$.
\end{subthm}

\end{thm}

\begin{thm}{Exercise}\label{ex:ultrapower(ultrapower)}
Given a nonprincipal ultrafilter $\omega$, construct an ultrafilter $\omega_1$ such that 
\[\spc{X}^{\omega_1}\iso(\spc{X}^\omega)^\omega\]
for any metric space~$\spc{X}$.

\end{thm}

\begin{thm}{Observation}\label{obs:ultrapower-is-geodesic}
Let $\spc{X}$ be a complete metric space. 
Then $\spc{X}^\omega$ is geodesic space if and only if $\spc{X}$ is a length space.
\end{thm}

\parit{Proof.}
The if part follows from \ref{obs:ultralimit-is-geodesic}; it remains to prove the only-if part

Assume $\spc{X}^\omega$ is geodesic.
Then any pair of points $x,y\in \spc{X}$ has a midpoint $z_\omega\in\spc{X}^\omega$.
Fix a sequence of points $z_n\in  \spc{X}$ such that $z_n\to z_\omega$ as $n\to \omega$.

Note that 
$\dist{x}{z_n}{\spc{X}}\to \tfrac12\cdot \dist{x}{y}{\spc{X}}$
and 
$\dist{y}{z_n}{\spc{X}}\to \tfrac12\cdot \dist{x}{y}{\spc{X}}$
as 
$n\to\omega$.
In particular, for any $\eps>0$, the point $z_n$ is an $\eps$-midpoint of $x$ and $y$ for $\omega$-almost all $n$.
It remains to apply Menger's lemma (\ref{lem:mid>geod}).
\qeds

\begin{thm}{Exercise}\label{ex:two-geodesics-in-ultrapower}
Assume $\spc{X}$ is a complete length space 
and $p,q\in\spc{X}$ cannot be joined by a geodesic in $\spc{X}$.  
Show that there is at least a continuum of distinct geodesics between $p$ and $q$ 
in the ultrapower $\spc{X}^\omega$.
\end{thm}

\begin{thm}{Exercise}\label{ex:notproper-limit}
Construct a proper metric space $\spc{X}$ such that its big ultrapower $\spc{X}^\omega$ is not locally compact.
\end{thm}

\section{Tangent and asymptotic spaces}
\label{sec:tan+asymptotic}

Choose a space $\spc{X}$ and a sequence $\lambda_n$ of positive numbers.
Consider the sequence of \index{rescaled space}\emph{rescalings} $\spc{X}_n=\lambda_n\cdot\spc{X}=(\spc{X},\lambda_n\cdot\dist{*}{*}{\spc{X}})$.

Choose a point $p\in \spc{X}$ and denote by $p_n$ the corresponding point in $\spc{X}_n$.
Consider the $\omega$-limit $\spc{X}_\omega$ of $\spc{X}_n$ (one may denote it by $\lambda_\omega\cdot \spc{X}$);
set $p_\omega$ to be the $\omega$-limit of $p_n$.

If $\lambda_n\to \infty$ as $n\to\omega$, then the metric component of $p_\omega$ in $\spc{X}_\omega$ is called \index{$\lambda_\omega$-tangent space}\emph{$\lambda_\omega$-tangent space} at $p$ and denoted by $\T_p^{\lambda_\omega}\spc{X}$ (or $\T_p^{\omega}\spc{X}$ if $\lambda_n=n$).\label{page:ultratangent space}

If $\lambda_n\to 0$ as $n\to\omega$, then the metric component of $p_\omega$ is called \index{$\lambda_\omega$-asymptotic space}\emph{$\lambda_\omega$-asymptotic space}%
\footnote{Often it is called an \textit{asymptotic cone} despite that it is not a cone in general; this name is used since in good cases it has a cone structure.} and denoted by $\Asym\spc{X}$ or $\Asym^{\lambda_\omega}\spc{X}$.
Note that the space $\Asym\spc{X}$ and its point $p_\omega$ do not depend on the choice of $p\in \spc{X}$.

The following exercise states that the constructions above depend on the sequence $\lambda_n$ and a nonprincipal ultrafilter $\omega$.

\begin{thm}{Exercise}\label{ex:ultraT}
Construct a metric space $\spc{X}$ with a point $p$ such that the tangent space
$\T_p^{\lambda_\omega}\spc{X}$ (or its asymptotic cone $\Asym^{\lambda_\omega}\spc{X}$) depends on the sequence $\lambda_n$ and/or ultrafilter~$\omega$.
\end{thm}

For nice spaces, different choices of the sequence of coefficients and ultrafilter may give the same space; 
some examples are given in the following exercise.

\begin{thm}{Exercise}\label{ex:Asym(Lob)}
Let $\spc{T}=\Asym\spc{L}$, where $\spc{L}$ is one of the following spaces
\begin{enumerate}[(i)]
 \item Lobachevsky plane,
 \item Lobachevsky space, or
 \item 3-regular metric tree; that is, the degree of any vertex. Assume that each edge has unit length.
\end{enumerate}

\begin{subthm}{ex:Asym(Lob):metric-tree}
Show that $\spc{T}$ is a complete metric tree.
\end{subthm}

\begin{subthm}{ex:Asym(Lob):homogeneous}
Show that $\spc{T}$ is one-point-homogeneous; that is, given two points $s,t\in \spc{T}$ there is an isometry of $\spc{T}$ that maps $s$ to $t$.
\end{subthm}

\begin{subthm}{ex:Asym(Lob):continuum}
Show that $\spc{T}$ has \index{degree}\emph{continuum degree} at any point;
that is, for any point $t\in \spc{T}$ the set of connected components of the complement $\spc{T}\setminus\{t\}$ has cardinality continuum.
\end{subthm}

\end{thm}

\begin{thm}{Exercise}\label{ex:T(Sx[0,1]/Sx0)}
Consider the cylinder $\mathbb{S}^1\times[0,1]$ with the standard metric.
Let $\spc{X}$ be the quotient space $\mathbb{S}^1\times[0,1]/\mathbb{S}^1\times\{0\}$;
that is, 
\[\dist{(u_1,t_1)}{(u_2,t_2)}{\spc{X}}
\df
\min\{\,\dist{(u_1,t_1)}{(u_2,t_2)}{\mathbb{S}^1\times[0,1]},t_1+t_2\,\}.\]
Describe the ultratangent space $\T_o^{\omega}\spc{X}$, where $o\in\spc{X}$ is the point that corresponds to $\mathbb{S}^1\times\{0\}$.
\end{thm}

\section{Remarks}

A nonprincipal ultrafilter $\omega$ is called 
\emph{selective}\index{ultrafilter!selective ultrafilter}\index{selective ultrafilter} if for any partition of $\NN$ into sets $\{C_\alpha\}_{\alpha\in\IndexSet}$ such that $\omega(C_\alpha)\z=0$ for each $\alpha$, 
there is a set $S\subset \NN$ such that $\omega(S)=1$ and $S\cap C_\alpha$ is a one-point set for each $\alpha\in\IndexSet$.

The existence of a selective ultrafilter follows from the continuum hypothesis \cite{rudin}.

If needed, we may assume that the chosen ultrafilter $\omega$ is selective.
In this case, \textit{the subsequence $(x_n)_{n\in S}$ in \ref{prop:ultra/partial} can be chosen so that $\omega(S)=1$}.

%% file: metric-sol.tex
\refstepcounter{chapter}
\setcounter{eqtn}{0}

\parbf{\ref{ex:quad-inq}.}
Add four triangle inequalities (\ref{metric:triangle}).

\parbf{\ref{ex:normal}.}
Consider the function 
\[f(x)=\frac{\distfun_Ax}{\distfun_Ax+\distfun_Bx},\]
where $\distfun_Ax\z\df\inf_{a\in A}\dist{a}{x}{}$.
Show that $f$ is continuous and satisfies the needed property.

\parbf{\ref{ex:tietze}.}
Use \ref{ex:normal} to construct an approximation of the needed function and pass to a limit or find a proof of the \index{Tietze extension theorem}\emph{Tietze extension theorem}.

\parbf{\ref{ex:pseudo-infty-metric}};\ref{SHORT.ex:pseudo-infty-metric:pseudo}.
Note that if $\mu(A)=\mu(B)=0$, then $|A-B|=0$.
Therefore, \ref{metric=0} does not hold for bounded closed subsets.
It is straightforward to check the remaining conditions in~\ref{def:metric} hold true.

\parit{\ref{SHORT.ex:pseudo-infty-metric:infty}.}
Note that the distance from the empty set to the whole plane is infinite; so the value $|A-B|$ might be infinite.
It is straightforward to check the remaining conditions in~\ref{def:metric}.

\parit{Remark.}
Metrics of the form $\dist{A}{B}{}=\mu(A\bigtriangleup B)$ are very special.
In particular, they satisfy the so-called \index{hypermetric inequality}\emph{hypermetric inequalities}; that is, for any sequence of sets $A_1,\dots, A_n$ and any sequence of integers $b_1,\dots,b_n$ such that $\sum_ib_i=1$ we have
\[\sum_{i,j}b_i\cdot b_j\cdot \dist{A_i}{A_j}{}\le 0.\]
Note that for $n=3$ and $b_1=b_2=-b_3=1$ we get the usual triangle inequality.
For more on the subject, see \cite{deza-laurent}.

\parbf{\ref{ex:gluing}.}
Choose $\delta>0$ and an increasing linear bijection $\ell\:[a,b]\to [c,d]$.
Show that $\ell$ has arbitrarily close increasing piecewise-linear bijection $s\:[a,b]\to [c,d]$ such that derivative at any point is either $<\delta$ or $>\tfrac1\delta$.

Start with the identity map $[0,1]\to [0,1]$;
iterate the above construction for smaller and smaller $\delta$ and pass to the limit.
This way we obtain an increasing  bijection $x\leftrightarrow x'$ from $[0,1]$ to itself
such that for any $\eps>0$ there is a partition $0=t_0<t_1<\dots <t_{2\cdot n}=1$ of $[0,1]$ with 
\begin{align*}
\eps&>|t_0-t_1|+|t_1'-t_2'|+|t_2-t_3|+\dots
\\
&\dots+|t_{2\cdot n-2}-t_{2\cdot n-1}|+|t_{2\cdot n-1}'-t_{2\cdot n}'|.
\end{align*}
Make a conclusion.

\parbf{\ref{ex:almost-min}.}
Assume the statement is wrong. 
Then for any point $x\in \spc{X}$, there is a point $x'\in \spc{X}$ such that 
\begin{align*}
\dist{x}{x'}{}&<\rho(x)
\intertext{and}
\rho(x')&\le\frac{\rho(x)}{1+\eps}.
\end{align*}
Consider a sequence $x_1,x_2,\dots\in \spc{X}$ such that $x_{n+1}\z=x_n'$.
Show that this is a Cauchy sequence.
Since $\spc{X}$ is complete, $x_n$ converges;
denote its limit by $x_\infty$.
Since $\rho$ is a continuous function we get
\begin{align*}
\rho(x_\infty)&=\lim_{n\to\infty}\rho(x_n)=0.
\end{align*}

The latter contradicts that $\rho>0$.

\parbf{\ref{ex:complete-completion}.}
Let $\bar {\spc{X}}$ be the completion of $\spc{X}$.
By the definition, for any $y\in \bar {\spc{X}}$ there is a Cauchy sequence $x_n$ in  $\spc{X}$ that converges to $y$.

Choose a Cauchy sequence $y_m$ in $\bar {\spc{X}}$.
From above, we can choose points $x_{n,m}\in \spc{X}$ such that $x_{n,m}\to y_m$ for any $m$.
Choose $z_m=x_{n_m,m}$ such that $|y_m-z_m|<\tfrac1m$.
Observe that $z_m$ is Cauchy.
Therefore, its limit $z_\infty$ lie in $\bar{\spc{X}}$.
Finally, show that $x_m\to z_\infty$.

\parbf{\ref{ex:compact-net}.}
A compact $\eps$-net $N$ in $\spc{K}$ contains a finite $\eps$ net $F$.
Show and use that $F$ is a $2\cdot\eps$-net of $\spc{K}$.

\parbf{\ref{ex:non-contracting-map}.}
Given a pair of points $x_0,y_0\in \spc{K}$, 
consider two sequences $x_0,x_1,\dots$ and $y_0,y_1,\dots$
such that $x_{n+1}=f(x_n)$ and $y_{n+1}\z=f(y_n)$ for each $n$.

Since $\spc{K}$ is compact, 
we can choose an increasing sequence of integers $n_k$
such that both sequences $(x_{n_i})_{i=1}^\infty$ and $(y_{n_i})_{i=1}^\infty$
converge.
In particular, both are Cauchy;
that is,
\[
|x_{n_i}-x_{n_j}|_{\spc{K}}\to 0 
\quad\text{and}\quad
|y_{n_i}-y_{n_j}|_{\spc{K}}\to 0
\]
as $\min\{i,j\}\to\infty$.

Since $f$ is distance-noncontracting, 
\[
|x_0-x_{|n_i-n_j|}|
\le 
|x_{n_i}-x_{n_j}|
\]
for any $i$ and $j$.
Therefore, there is a sequence $m_i\to\infty$ such that
\[
x_{m_i}\to x_0\quad\text{and}\quad y_{m_i}\to y_0
\leqno({*})\]
as $i\to\infty$.

Since $f$ is distance-noncontracting, the sequence $\ell_n=|x_n-y_n|_{\spc{K}}$ is nondecreasing.
By $({*})$,  $\ell_{m_i}\to\ell_0$ as $m_i\to\infty$.
It follows that 
\[\ell_0=\ell_1=\dots\]
In particular, 
\[|x_0-y_0|_{\spc{K}}=\ell_0=\ell_1=|f(x_0)-f(y_0)|_{\spc{K}}\]
for any pair of points $(x_0,y_0)$ in $\spc{K}$.
That is, the map $f$ is distance-preserving; hence $f$ is injective.
From $({*})$, we also get that $f(\spc{K})$ is everywhere dense.
Since $\spc{K}$ is compact $f\:\spc{K}\to \spc{K}$ is surjective --- hence the result.

\parit{Remarks.}
This is a basic lemma in the introduction to Gromov--Hausdorff distance \cite[see 7.3.30 in][]{burago-burago-ivanov}.
The presented proof is not quite standard,
I learned it from Travis Morrison, 
a student in my MASS class at Penn State, Fall 2011.

Note that this exercise implies that \textit{any surjective non-expanding map from a compact metric space to itself is an isometry}. 

\parbf{\ref{ex:loc-compact-not-proper}.}
Check an infinite set with a discrete metric.

\parbf{\ref{ex:pogorelov}.}
Set $B_p=B(x,\tfrac \pi2)_{\mathbb{S}^2}$.
The triangle inequality follows since
\[
(B_x\setminus B_y)
\cup 
(B_y\setminus B_z)
\supseteq
B_x\setminus B_z.
\leqno(*)\]
The remaining conditions in Definition \ref{def:metric} are evident.

Observe that
$B_x\setminus B_y$
does not overlap with
$B_y\z\setminus B_z$ and  we get equality in $(*)$ if and only if $y$ lies on the great circle arc from $x$ to $z$.
Therefore, the second statement follows.

\begin{wrapfigure}{r}{24 mm}
\vskip-0mm
\centering
\includegraphics{mppics/pic-29}
\end{wrapfigure}

\parit{Remarks.}
This construction is due to 
Aleksei Pogorelov \cite{pogorelov}.
It is closely related to the construction given 
by David Hilbert \cite{hilbert}
which was the motivating example for his fourth problem. 
See also the remark after the solution of~\ref{ex:pseudo-infty-metric}.

\parbf{\ref{ex:4-point-trees}.}
We may assume that non of the points $p,x,y,z$ lie on a geodesic between the other two.

Let $K$ be the set in the tree covered by all six geodesics with the endpoints $p,x,y,z$.
Observe that $K$ looks like an H or like an X; make a conclusion.

\parit{Remarks.}The value $\tfrac12\cdot(|p-x|+|p-y|-|x-y|)$ is called \index{Gromov's product}\emph{Gromov's product} of $x$ and $y$ with the origin at $p$;
usually it is denoted by $(x|y)_p$.

Note that a four-point metric space admits an isometric embedding into a metric tree if and only if one of these two equivalent conditions holds.
Moreover, a metric space admits an isometric embedding into a metric tree if every its four-point subspace admits such embedding.

\parbf{\ref{ex:spheres-in-trees}.}
Apply \ref{ex:4-point-trees}.

\parbf{\ref{ex:1-Lip-G-delta}.}
Note that $\spc{P}$ is complete.
Choose $\eps>0$.
Use \ref{thm:length-semicont} to show that the set of paths of length $>1-\eps$ is open in~$\spc{P}$;
show that this set is also dense in~$\spc{P}$.
Apply Baire's theorem (\ref{thm:baire}).

\parit{Remark.}
You might find it surprising that \textit{most of the short maps from the sphere to the plane are \index{length-preserving map}\emph{length-preserving}}; that is, they preserve lengths of all curves.
The latter follows from the result of 
Bernd Kirchheim, 
Emanuele Spadaro,
and 
L{\'a}szl{\'o} Sz{\'e}kelyhidi \cite{KSS}.
(While most of the maps have this property, it is not at all easy to construct a single such example.)

\parbf{\ref{ex:no-geod}.}
\textit{Formally speaking, a one-point space is a solution,
but we will construct a nontrivial example.}

Recall that $c_0\subset\ell^\infty$ denotes the space of all real sequences converging to zero.
Consider the unit ball $B$ in $c_0$;
denote by $\rho_0$ the metric on $B$.

Let \[\phi(\bm{x})=2+\tfrac{1}2\cdot x_1+\tfrac{1}4\cdot x_2+\tfrac{1}8\cdot x_3+\dots,\]
where $\bm{x}=(x_1,x_2\,\dots)\in B$.
Consider another length metric $\rho_1$ on $B$ that is different from $\rho_0$ by the conformal factor $\phi$;
that is, if $t\mapsto\bm{x}(t)$ for $t\in[0,\ell]$ is a curve parametrized by $\rho_0$-length,
then its $\rho_1$-length is defined by
\[\length_{\rho_1}\bm{x}\df\int\limits_0^\ell\phi\circ\bm{x}(t)\cdot dt.\]
Note that the metric $\rho_1$ is bilipschitz to~$\rho_0$.

Assume $t\mapsto \bm{x}(t)$ and $t\mapsto \bm{x}'(t)$ are two curves parametrized by $\rho_0$-length that differ only in the $m$-th coordinate; denote them by $x_m(t)$ and $x_m'(t)$ respectively.
Show that if $x'_m(t)\le x_m(t)$ for any $t$ and 
the function $x'_m(t)$ is locally $1$-Lipschitz at all $t$ such that $x'_m(t)< x_m(t)$, then 
\[\length_{\rho_1}\bm{x}'\le \length_{\rho_1}\bm{x}.\]
Moreover, this inequality is strict if $x'_m(t)\z< x_m(t)$ for some~$t$.

Fix a curve $\bm{x}(t)$, $t\in[0,\ell]$, parametrized by  $\rho_0$-length.
We can choose large $m$ so that $x_m(t)$ is sufficiently close to $0$ for any~$t$.
In this case, it is easy to construct a function $t\mapsto x'_m$ that meets the above conditions.
It follows that for any curve $\bm{x}(t)$ in $(B,\rho_1)$, we can find a shorter curve $\bm{x}'(t)$ with the same endpoints.
In particular, $(B,\rho_1)$ has no geodesics.

\parit{Remark.}
This solution was suggested by Fedor Nazarov~\cite{nazarov}.

\parbf{\ref{ex:compact+connceted}.}
Choose a sequence of positive numbers $\varepsilon_n\to 0$ and a finite $\varepsilon_n$-net $N_n$ of $K$ for each $n$.
We can assume that $\eps_0>\diam K$, and $N_0$ is a one-point set.
If $\dist{x}{y}{}<\eps_k$ for some $x\in N_{k+1}$ and $y\in N_{k}$, then connect them by a curve of length at most $\eps_k$.

Let $K'$ be the union of all these curves and $K$.
Show that $K'$ is compact and path-connected.

\parit{Source:} This problem is due to Eugene Bilokopytov \cite{bilokopytov}.

\parbf{\ref{ex:compact=>complete}.}
Choose a Cauchy sequence $x_n$ in $(\spc{X},\|*\z-*\|)$; it is sufficient to show that a subsequence of $x_n$ converges.

Observe that the sequence $x_n$ is Cauchy in $(\spc{X},|*-*|)$;
denote its limit by $x_\infty$.

Passing to a subsequence, we can assume that $\|x_n-x_{n+1}\|\z<\tfrac1{2^n}$.
It follows that there is a 1-Lipschitz path $\gamma$ in $(\spc{X},\|*-*\|)$ such that $x_n=\gamma(\tfrac1{2^n})$ for each $n$ and $x_\infty=\gamma(0)$.
Therefore,
\begin{align*}
\|x_\infty-x_n\|&\le \length\gamma|_{[0,\frac1{2^n}]}\le \tfrac1{2^n}.
\end{align*}
In particular, $x_n$ converges to $x_\infty$ in $(\spc{X},\|*\z-*\|)$.

\parit{Source:} \cite[Corollary]{hu-kirk}; see also \cite[Lemma 2.3]{petrunin-stadler}.

\parbf{\ref{ex:eps-nbhd(ball)}.}
Let $U$ be the $\eps$-neighborhood of $\oBall(x,R)_{\spc{X}}$.
By the triangle inequality, $U\z\subset \oBall(x,R+\eps)_{\spc{X}}$;
this inclusion holds in any metric space.

Choose $y\in \oBall(x,R+\eps)_{\spc{X}}$, so $\dist{x}{y}{{\spc{X}}}\z<R+\eps$.
Since ${\spc{X}}$ is a length space, there is a curve $\gamma$ from $x$ to $y$ with length less than $R+\eps$.
Show and use that $\gamma$ contains a point $m$ such that $\dist{x}{m}{{\spc{X}}}<R$ and $\dist{y}{m}{{\spc{X}}}<\eps$.

\parbf{\ref{exercise from BH}.}
Consider the following subset of $\RR^2$ equipped with the induced length metric
\[
\spc{X}
=
\bigl((0,1]\times\{0,1\}\bigr)
\cup
\bigl(\{1,\tfrac12,\tfrac13,\dots\}\times[0,1]\bigr)
\]
Note that $\spc{X}$ is locally compact and geodesic.

Its completion $\bar{\spc{X}}$ is isometric to the closure of $\spc{X}$ equipped with the induced length metric.
Note that $\bar{\spc{X}}$ is obtained from $\spc{X}$ by adding two points $p=(0,0)$ and $q\z=(0,1)$.

{

\begin{wrapfigure}{r}{20 mm}
\vskip-4mm
\centering
\includegraphics{mppics/pic-1}
\end{wrapfigure}

Observe that $p$ admits no compact neighborhood in $\bar{\spc{X}}$ and there is no geodesic connecting $p$ to $q$ in~$\bar{\spc{X}}$.

\parit{Source:} \cite[I.3.6(4)]{bridson-haefliger}.

}

\parbf{\ref{ex:gross}.}
Suppose this number does not exist.
Show that there are two point-arrays $(x_1,\z\dots,x_n)$ and $(y_1,\dots,y_m)$
such that
\[
\min_{z\in \spc{X}}\{\,f(z)\,\}>\max_{z\in \spc{X}}\{\,h(z)\,\},
\leqno({*})
\]
where
\begin{align*}
f(z)&=\tfrac1n\cdot\sum_i|x_i-z|_{\spc{X}}
\intertext{and}
h(z)&=\tfrac1m\cdot\sum_j|y_j-z|_{\spc{X}}.
\end{align*}

Note that
\begin{align*}\tfrac1m\cdot\sum_j f(y_j)&=\tfrac1{m\cdot n}\cdot\sum_{i,j}|x_i-y_j|_{\spc{X}}=
\\
&=\tfrac1n\cdot\sum_i h(x_i);
\end{align*}
that is, the average value of $f(y_j)$ coincides with the average value of $h(x_i)$.
The latter contradicts~$({*})$.

\parit{Remark.}
The value $\ell$ is uniquely defined;
it is called the \index{rendezvous value}\emph{rendezvous value} of ${\spc{X}}$.
This is a result of Oliver Gross \cite{gross}.


%% file: uryson-sol.tex
\refstepcounter{chapter}
\setcounter{eqtn}{0}

\parbf{\ref{ex:compact-length}.}
By the Fréchet lemma (\ref{lem:frechet}) we can identify $\spc{K}$ with a compact subset in $\ell^\infty$.

Denote by $\spc{L}$ the \index{closed convex hull}\emph{closed convex hull} of $\spc{K}$;
that is, $\spc{L}$ is the minimal convex closed set in $\ell^\infty$ that contains $\spc{K}$.
(In other words, $\spc{L}$ is the minimal closed set containing $\spc{K}$ such that if $x,y\in \spc{L}$, then 
$t\cdot x+(1-t)\cdot y\in \spc{L}$ for any $t\in[0,1]$.)

Observe that $\spc{L}$ is a length space.
It remains to show that $\spc{L}$ is compact.

By construction, $\spc{L}$ is a closed subset of $\ell^\infty$; in particular, it is complete.
By \ref{totally-bounded}, it remains to show that $\spc{L}$ is totally bounded.

Recall that Minkowski sum $A + B$ of two sets $A$ and $B$ in a vector space is defined by
\[A + B 
\df
\set{a+b}{a\in A,\ b\in B}.\]
Observe that the Minkowski sum of two convex sets is convex.

Denote by $\bar B_\eps$ the closed $\eps$-ball in $\ell^\infty$ centered at the origin.
Choose a finite $\eps$-net $N$ in $\spc{K}$ for some $\eps>0$.
Note that $P=\Conv N$ is a convex polyhedron; in particular, $\Conv N$ is compact.

Observe that $N+\bar B_\eps$ is a closed $\eps$-neighborhood of $N$.
It follows that $N+\bar B_\eps\supset K$ and therefore $P+\bar B_\eps\supset \spc{L}$.
In particular, $P$ is a $2\cdot\eps$-net in $\spc{L}$;
since $P$ is compact and $\eps>0$ is arbitrary, $\spc{L}$ is totally bounded (see \ref{ex:compact-net}).

\parit{Remark.}
Alternatively, one may use that \textit{the injective envelope of a compact space is compact}; see \ref{ex:inj=complete-geodesic-contractible:geodesic}, \ref{ex:Inj(compact)}, and \ref{prop:InjX-is-injective}.

\parbf{\ref{ex:frechet}.}
Modify the proof of \ref{lem:frechet}.

\begin{wrapfigure}{r}{23mm}
\vskip-6mm
\centering
\includegraphics{mppics/pic-200}
\end{wrapfigure}

\parbf{\ref{ex:inf-extension}.}
Consider the metric tree $\spc{T}$ shown on the diagram;
it is a half-line $[0,\infty)$ with attached an interval of length $n+1$ to each integer~$n\ge 0$.
Denote by $o$ the origin of the half-line
and by $x_n$ the endpoint of $n^{\text{th}}$ interval.

Observe that if $m\ne n$, then
\[|x_m-x_n|_{\spc{T}}\ge |o-x_n|_{\spc{T}}+1.\]
Show and use that for any binary sequence $\eps_n$ there is an extension function $f$ such that 
\[f(x_n)=|o-x_n|_{\spc{T}}+\eps_n.\]

\parit{Remark.}
An if-and-only-if condition on $\spc{X}$ that have separable $\spc{X}^\infty$ was found by Julien Melleray \cite[2.8]{melleray}.
A similar condition was used by Herbert Federer to describe metric spaces where Besicovitch covering lemma holds \cite[2.8.9]{federer}.

\parbf{\ref{ex:geodesics-urysohn}.}
Choose a separable space $\spc{X}$ that has an infinite number of geodesics between a pair of points with the given distance between them;
say a square in $\RR^2$ with $\ell^\infty$-metric will do.
Apply to $\spc{X}$ universality of Urysohn space (\ref{prop:sep-in-urys}).

\parbf{\ref{ex:compact-extension}.} 
First let us prove the following claim:

\begin{itemize}
\item 
Suppose $f\: K\to\RR$ is an extension function defined on a compact subset $K$ of the Urysohn space $\spc{U}$.
Then there is a point $p\in \spc{U}$ such that 
$\dist{p}{x}{}=f(x)$ for any $x\in K$.
\end{itemize}

Without loss of generality, we may assume that $f>0$.
Since $K$ is compact, we may fix $\eps>0$ such that $f(x)>\eps$ for any $x\in K$.

Consider the sequence $\eps_n=\tfrac\eps{100\cdot 2^n}$.
Choose a sequence of $\eps_n$-nets $N_n\subset K$.
Applying the universality of $\spc{U}$ recursively, we may choose a point $p_n$ such that $\dist{p_n}{x}{}=f(x)$ for any $x\in N_n$ and $\dist{p_n}{p_{n-1}}{}\z=10\cdot\eps_{n-1}$.
Observe that the sequence $p_n$ is Cauchy and its limit $p$ meets 
$\dist{p}{x}{}=f(x)$ for any $x\in K$.

Now, choose a sequence $x_n$ of points that is dense in $\spc{S}$.
Applying the claim, we may extend the map from $K$ to $K\cup\{x_1\}$, further to $K\cup\{x_1,x_2\}$, and so on.
As a result, we extend the distance-preserving map $f$ to the whole sequence $x_n$.
It remains to extend it continuously to the whole space~$\spc{S}$.

\parbf{\ref{ex:sc-urysohn}.}
It is sufficient to show that any compact subspace $\spc{K}$ of the Urysohn space $\spc{U}$ can be contracted to a point.

Note that any compact space $\spc{K}$ can be extended to a contractible compact space $\spc{K}'$; for example, we may embed $\spc{K}$ into $\ell^\infty$ and pass to its convex hull, as it was done in \ref{ex:compact-length}.

By \ref{thm:compact-homogeneous}, there is an isometric embedding of $\spc{K}'$ that agrees with the inclusion $\spc{K}\hookrightarrow\spc{U}$.
Since $\spc{K}$ is contractible in $\spc{K}'$, it is contractible in $\spc{U}$.

\parit{A better way.}
One can contract the whole Urysohn space using the following construction.

Note that points in $\spc{X}_\infty$ constructed in the proof of \ref{prop:univeral-separable} can be multiplied by $t\in [0,1]$ --- simply multiply each function by $t$.
That defines a map 
\[\lambda_t\:\spc{X}_\infty\to \spc{X}_\infty\]
that rescales all distances by factor $t$.
The map $\lambda_t$ can be extended to the completion of $\spc{X}_\infty$, which is isometric to $\spc{U}_d$ (or $\spc{U}$).

Observe that 
the map $\lambda_1$ is the identity  and $\lambda_0$ maps the whole space to a single point, say $x_0$ --- this is the only point of $\spc{X}_0$.
Further, note that $(t,p)\mapsto \lambda_t(p)$ is a continuous map; in particular, $\spc{U}_d$ and $\spc{U}$ are contractible.

As a bonus, observe that for any point $p\in \spc{U}_d$ the curve $t\mapsto \lambda_t(p)$ is a geodesic path from $p$ to $x_0$.

\parit{Source:} \cite[$\text{(d)}$ on page 82]{gromov-2007}.

\parbf{\ref{ex:no-isom}.}
Consider two infinite metric trees as on the diagram. 

\begin{Figure}
\vskip-0mm
\centering
\includegraphics{mppics/pic-205}
\end{Figure}

\parit{Remark.}
A more sophisticated example: $\spc{X}\z=\ell^\infty$ and $\spc{Y}=L^\infty([0,1])$.
Try to prove that it qualifies; see also \cite{buehler}.




\parbf{\ref{ex:sphere-in-urysohn}}; \ref{SHORT.ex:sphere-in-urysohn:sphere} and \ref{SHORT.ex:sphere-in-urysohn:midpoint}.
Observe that $L$ and $M$ satisfy the definition of $d$-Urysohn space and apply the uniqueness (\ref{thm:urysohn-unique}).
Note that
\[\ell=\diam L=\min\{2\cdot r, d\}.\]

\parit{\ref{SHORT.ex:sphere-in-urysohn:homogeneous}.} 
Use \ref{SHORT.ex:sphere-in-urysohn:sphere}, maybe twice.

\parbf{\ref{ex:shere}.}
Let $p$ be the center of the sphere;
without loss of generality, we can assume that $\dist{p}{x}{}\le \dist{p}{y}{}$.

Consider function $f\:\{p,x,y\}\to\RR$ defined by $f(p)=1$, $f(x)=1+\dist{p}{x}{}$, and $f(y)=1+\dist{p}{y}{}-\eps$.
Suppose $\eps>0$ is sufficiently small;
show that $f$ is an extension function on $\{p,x,y\}$.

By the extension property, there is a point $z\in \spc{U}$ such that $\dist{p}{z}{}=f(p)$, $\dist{x}{z}{}=f(x)$, and $\dist{y}{z}{}=f(y)$.
Whence the statement follows.

\parit{Source:} This problem is taken from a survey of Julien Melleray
 \cite[Prop. 4.3]{melleray}, where it was attributed to Matatyahu Rubin.

\parbf{\ref{ex:ext(shere)}.} 
Observe that the complement $\spc{V}=\spc{U}\setminus B$ is complete.
Show that it $\spc{V}$ satisfies the extension property.
Conclude that $\spc{V}$ is an Urysohn space and apply \ref{thm:urysohn-unique}.

For the second part, observe that there is an isometry $\iota\:\spc{U}\to \spc{V}$.
Moreover, if $p$ is the center of $B$, then we can assume that $\iota$ has a fixed point $x$ such that $\dist{p}{x}{}>2$.

Consider the unit sphere $S$ centered at $x$.
The restriction of $\iota$ to $S$ is an isometry of $S$.
Use \ref{ex:shere} to show that it cannot be extended to an isometry of $\spc{U}$.

\parit{Source:} \cite[Sec. 4.4]{melleray}.

\parbf{\ref{ex:katetov}.}
Apply \ref{thm:urysohn-unique} and the construction in \ref{thm:urysohn-exists+}.

\parbf{\ref{ex:homogeneous}}; \ref{SHORT.ex:homogeneous:euclidean}.
The euclidean plane is homogeneous in every sense.

\parit{\ref{SHORT.ex:homogeneous:hilbert}.}
The Hilbert space $\ell^2$ is finite-set-homogeneous, but not compact-set-homogeneous, nor countable-set-homogeneous.

\parit{\ref{SHORT.ex:homogeneous:ell-infty}.}
$\ell^\infty$ is one-point-homogeneous, but not two-point-homogeneous.
Try to show that there is no isometry of $\ell^\infty$ such that
\begin{align*}
(0,0,0,\dots)&\mapsto (0,0,0,\dots),
\\
(1,1,1,\dots)&\mapsto (1,0,0,\dots).
\end{align*}

\parit{\ref{SHORT.ex:homogeneous:ell-1}.}
$\ell^1$ is one-point-homogeneous, but not two-point-homogeneous.
Try to show that there is no isometry of $\ell^\infty$ such that
\begin{align*}
(0,0,0,\dots)&\mapsto (0,0,0,\dots),
\\
(2,0,0\dots)&\mapsto (1,1,0,\dots).
\end{align*}

\parbf{\ref{ex:homogeneous-tree}.}
Let $\spc{T}$ be a one-point-homogeneous metric tree.
Note that all points in $\spc{T}$ have the same degree $d$;
that is, for any point $t\in \spc{T}$ the set of connected components of the complement $\spc{T}\setminus\{t\}$ has the same cardinality $d$.

Show that if $d=0$, then $\spc{T}$ is a one-point space;
there is no tree with $d=1$,
and if $d=2$, then $\spc{T}\iso\RR$.

Suppose $d\ge 3$.
Choose a geodesic $\gamma$ in $\spc{T}$.
Show that number of connected components of $\spc{T}\setminus\gamma$ has cardinality continuum.
Observe and use that one can choose a point $p_\alpha$ in each connected component such that $\dist{p_\alpha}{p_\beta}{\spc{T}}>1$ if $\alpha\ne\beta$.

\parbf{\ref{ex:horobry}.}
Assume $F_{\spc{X}}$ is not an embedding.
That is, there is a sequence of points $x_1,x_2,\dots$ 
and a point $x_\infty$ such that $f_{x_n}\to f_{x_\infty}$ in $C(\spc{X},\RR)$
as $n\to \infty$, 
while $|x_n-x_\infty|_{\spc{X}}>\eps$ 
for some fixed $\eps>0$ and all~$n$.

By \ref{prop:length+proper=>geodesic}, any pair of points $x,y\in \spc{X}$ can be connected by a minimizing geodesic $[xy]$.
Choose $\bar x_n$ on a geodesic $[x_\infty x_n]$ such that $|x_\infty-\bar x_n|=\eps$.
Note that 
\begin{align*}
f_{x_n}(x_\infty)-f_{x_n}(\bar x_n)&=\eps,
\\
f_{x_\infty}(x_\infty)-f_{x_n}(\bar x_n)&=-\eps
\end{align*}
for all $n$.

Since $\spc{X}$ is proper, we can pass to a subsequence of $x_n$ so that the sequence  $\bar x_n$ converges;
denote its limit by $\bar x_\infty$.
The above identities imply that
\begin{align*}
f_{x_n}(\bar x_\infty)&\not\to f_{x_\infty}(\bar x_\infty)
\quad
\text{or}
\\
f_{x_n}(x_\infty)&\not\to f_{x_\infty}( x_\infty)
\end{align*}
--- a contradiction.

For the second part, take $\spc{Y}$ to be the set of non-negative integers with the metric $\rho$ defined by
\[\rho(m,n)=m+n\] 
for $m\ne n$.

\medskip

\parit{Source:}
I learned this example from Linus Kramer and Alexander Lytchak;
it was also mentioned in the lectures of Anders Karlsson
and attributed to Uri Bader \cite[2.3]{karlsson}.

\parbf{\ref{ex:cut}.}
Suppose that our metric is $\sum a_S\cdot\delta_S$ with $a_S\ge 0$ for any $S\subset F$.
Enumerate all the subsets $S_1,\dots,S_{2^n}$;
set $S_i=F$ for all $i>2^n$. 
Consider the maps $x\mapsto (a_1,a_2,\dots)$ where $a_i=0$ if $x\in S_i$ and otherwise $a_i=1$.
Observe that it defines a distance-preserving map $F\to \ell^1$. 

The if part is proved.
For the only-if part, check the statement for subsets of the real line, and use it.

\parbf{\ref{ex:K23}.}
Show that for any proper subset $S$ in the vertex set there are three vertices $x,y,z$ such that $\dist{x}{y}{} +\dist{y}{z}{}=\dist{x}{z}{}$ and either 
$x,z\in S$ and $y\notin S$, or $x,z\notin S$ and $y\in S$.
Then apply \ref{ex:cut}.

\parbf{\ref{ex:RP-not}.}
For the first part, show and use that the quotient of $\RP^2$ by the isotropy group of one point is isometric to a line segment.

For the second part, choose three points on a closed geodesic at equal distances from each other.
Show and use that there is an isometric three-point set in $\RP^2$ that does not lie on a closed geodesic.

\parit{Source:} \cite[V \S 2]{busemann-1942}.

\parbf{\ref{ex:hom-cube}.}
Denote by $\dim(x_1,\dots,x_m)$ the dimension of the minimal face of the cube that contains all the points $x_1,\dots,x_m\in Q$.
Show and use that 
\[\dim(x_1,\dots,x_m)=\dim(x_1',\dots,x_m')\]
for any isometry $x\mapsto x'$ of $Q$.

\parit{Source:} \cite[prop. 6 and 7]{berestovskii-nikonorov}.

%% file: injective-sol.tex
\parbf{\ref{ex:conv-short}}; \textit{only-if part}.
To check convexity, assume that $B$ is a two-point subset.
For closeness, assume that $B$ is a countable set of $A$.

\parit{If part.}
Apply the Kirszbraun theorem together with the closest-point projection.

\refstepcounter{chapter}
\setcounter{eqtn}{0}

\parbf{\ref{ex:inj=complete-geodesic-contractible}.}
Choose an injective space $\spc{Y}$.

\textit{\ref{SHORT.ex:inj=complete-geodesic-contractible:complete}.}
Fix a Cauchy sequence $x_n$ in $\spc{Y}$;
we need to show that it has a limit $x_\infty\in \spc{Y}$.
Consider metric on $\spc{X}=\NN\cup\{\infty\}$ defined by 
\begin{align*}
\dist{m}{n}{\spc{X}}&\df\dist{x_m}{x_n}{\spc{Y}},
\\
\dist{m}{\infty}{\spc{X}}&\df\lim_{n\to\infty}\dist{x_m}{x_n}{\spc{Y}}.
\end{align*}
Since the sequence is Cauchy, so is the sequence $\ell_n=\dist{x_m}{x_n}{\spc{Y}}$ for any $m$.
Therefore, the last limit is defined.

By construction, the map $n\mapsto x_n$ is distance-preserving on $\NN\subset \spc{X}$.
Since $\spc{Y}$ is injective, this map can be extended to $\infty$ as a short map; set $\infty\mapsto x_\infty$.
Since $\dist{x_n}{x_\infty}{\spc{Y}}\le \dist{n}{\infty}{\spc{X}}$ 
and $\dist{n}{\infty}{\spc{X}}\to 0$, we get that
$x_n\to x_\infty$ as $n\to\infty$.

\textit{\ref{SHORT.ex:inj=complete-geodesic-contractible:geodesic}.}
Applying the definition of injective space, we get a midpoint for any pair of points in $\spc{Y}$.
By \ref{SHORT.ex:inj=complete-geodesic-contractible:complete},
$\spc{Y}$ is a complete space.
It remains to apply Menger's lemma (\ref{lem:mid>geod:geod}).

\textit{\ref{SHORT.ex:inj=complete-geodesic-contractible:contractible}.}
Let $k\:\spc{Y}\hookrightarrow \ell^\infty(\spc{Y})$ be the Kuratowski embedding (\ref{lem:kuratowski}).
Observe that $\ell^\infty(\spc{Y})$ is contractible;
in particular, there is a homotopy $k_t\:\spc{Y}\hookrightarrow \ell^\infty(\spc{Y})$ such that $k_0=k$ and $k_1$ is a constant map.
(In fact, one can take $k_t=(1-t)\cdot k$.)

Since $k$ is distance-preserving and $\spc{Y}$ is injective,
there is a short map $f\:\ell^\infty(\spc{Y})\to \spc{Y}$ such that the composition $f\circ k$ is the identity map on $\spc{Y}$.
The composition $f\circ k_t\:\spc{Y}\hookrightarrow \spc{Y}$ provides the needed homotopy. 

\parbf{\ref{ex:bicombing}.}
By \ref{lem:kuratowski}, the space $\spc{Y}$ can be considered as a subset in $\ell^\infty(\spc{Y})$.
Let $r\: \ell^\infty(\spc{Y})\to \spc{Y}$ be the short retraction provided by the definition of injective space.
Observe that $m(x,y)\df r(\tfrac{x+y}2)$ meets the condition.

\parit{Remark.} The same argument can be used to construct the so-called \index{geodesic bicombing}\emph{geodesic bicombing} on injective space --- a useful tool introduced by Urs Lang \cite[3.6]{lang-2013}.

\parbf{\ref{ex:injective-spaces}.}
Suppose that a short map $f\:A\to\spc{Y}$ is defined on a subset $A$ of a metric space $\spc{X}$.
We need to construct a short extension $F$ of $f$.
Without loss of generality, we may assume that $A\ne\emptyset$;
otherwise, map the whole $\spc{X}$ to a single point.
By Zorn's lemma, it is sufficient to enlarge $A$ by a single point $x\notin A$.

\parit{\ref{SHORT.ex:injective-spaces:R}.}
Suppose $\spc{Y}=\RR$.
Set 
\[F(x)=\inf\set{f(a)-\dist{a}{x}{}}{a\in A}.\] 
Observe that $F$ is short and $F(a)=f(a)$ for any $a\in A$.

\parit{\ref{SHORT.ex:injective-spaces:tree}.}
Suppose  $\spc{Y}$ is a complete metric tree.
Fix points $p\in \spc{X}$ and $q\in\spc{Y}$.
Given a point $a\in A$,
let $x_a\in\cBall[f(a),\dist{a}{p}{}]$ be the point closest to $f(x)$.
Note that $x_a\in[q\,f(a)]$ and either $x_a=q$ or $x_a$ lies on distance $\dist{a}{p}{}$ from $f(a)$.

Note that the geodesics $[q\,x_a]$ are nested;
that is, for any $a,b\in A$ we have either $[q\,x_a]\z\subset [q\,x_b]$ or $[q\,x_b]\z\subset [q\,x_a]$.
Moreover, in the first case, we have $\dist{x_b}{f(a)}{}\le \dist{p}{a}{}$ and in the second $\dist{x_a}{f(b)}{}\le \dist{p}{b}{}$.

It follows that the closure of the union of all geodesics $[q\,x_a]$ for $a\in\spc{A}$ is a geodesic.
Denote by $x$ its endpoint; it exists since $\spc{Y}$ is complete.
It remains to observe that $\dist{x}{f(a)}{}\le \dist{p}{a}{}$ for any $a\in\spc{A}$;
that is, one can take $f(p)=x$.

\parit{\ref{SHORT.ex:injective-spaces:ell-infty}.}
Show and use that \textit{any $\ell^\infty$-product of injective spaces is injective};
in particular, if $\spc{Y}$ and $\spc{Z}$ are injective, then so is the product $\spc{Y}\times\spc{Z}$ equipped with the metric 
\[\dist{(y,z)}{(y',z')}{\spc{Y}\times\spc{Z}}=\max\{\,\dist{y}{y'}{\spc{Y}},\dist{z}{z'}{\spc{Z}}\,\}.\]

\parbf{\ref{ex:extr-ball}}; \ref{SHORT.ex:extr-ball:one}.
Let $\spc{B}=\cBall[o,R]_{\spc{Y}}$.
Choose a metric space $\spc{X}$ with a subset $A$.
Given a short map $f\:A\to \spc{B}$ we need to find its short extension $\spc{X}\to \spc{B}$.

Since $\diam\spc{B}\le 2\cdot R$, we may assume that  $\diam \spc{X}\le 2\cdot R$;
if not pass to the metric defined by $\dist{x}{y}{}\df\max\{\,\dist{x}{y}{\spc{X}},2\cdot R\,\}$.

Let us add point $w$ to $\spc{X}$ such that $\dist{w}{x}{}=R$ for any $x\in\spc{X}$;
denote the obtained space $\spc{X}'$.
Let $f'\:A\cup\{w\}\to \spc{B}$ be an extension of $f$ by $w\mapsto o$; note that $f'$ is short.

Since $\spc{Y}$ is injective, there is a short extension $F\:\spc{X}'\to \spc{Y}$ of $f'$.
Show and use that $F(\spc{X}')\subset \spc{B}$.

\parit{\ref{SHORT.ex:extr-ball:many}.}
Try to modify the argument in \ref{SHORT.ex:extr-ball:one}.~Namely, let $\spc{B}=\bigcap_\alpha\cBall[o_\alpha,R_\alpha]_{\spc{Y}}$.
Note that one may assume that $\diam \spc{X}\z\le 2\cdot \inf_\alpha\{\,R_\alpha\,\}$.
Consider the space $\spc{X}'\z=\spc{X}\cup\{w_\alpha\}$ such that $\dist{w_\alpha}{x}{}\z=R_\alpha$ for any $x\in \spc{X}$ and $\dist{w_\alpha}{w_\beta}{}=R_\alpha+R_\beta$ if $\alpha\ne\beta$.
Further, consider an extension of $f$ by $w_\alpha\mapsto o_\alpha$.

\parbf{\ref{ex:extr-fixed}.}
Let $\diam \spc{Y}=2\cdot R$.
We can assume that $R>0$; otherwise, there is nothing to prove.
Denote by $\spc{Z}$ a minimal (with respect to inclusion) intersection of closed $R$-balls in $\spc{Y}$ such that $s(\spc{Z})\subset\spc{Z}$.

Consider 
the intersection 
\[\spc{Y}'=\spc{Z}\cap\left(\bigcap_{p\in \spc{Z}} \cBall[p,R]_{\spc{Y}}\right).\]
By \ref{ex:extr-ball:many}, $\spc{Y}'$ is injective.
Use that $\spc{Z}$ is minimal to show that $s(\spc{Y}')\subset \spc{Y}'$.
Show that $\diam \spc{Y}'\le \tfrac12\cdot\diam \spc{Y}$.

Consider a sequence of nested injective spaces $\spc{Y}=\spc{Y}_0\supset \spc{Y}_1\supset\dots$ such that $\spc{Y}_{n+1}\z=\spc{Y}_{n}'$.
Choose a point $y_n\in \spc{Y}_{n}$ for each $n$.
Show that the sequence $y_n$ is Cauchy.
By \ref{ex:inj=complete-geodesic-contractible:complete}, $y_n$ converges, say to $y_\infty$.
Observe that $y_\infty$ is a fixed point of $s$.

\parbf{\ref{ex:circle}}; \textit{only-if part}.
Suppose $r$ is extremal.
By \ref{lem:extremal-lipschitz}, $r$ is $1$-Lipschitz.
Since $\mathbb{S}^1$ is compact, \ref{lem:opposite-compact} implies that for any $p\in \mathbb{S}^1$ there is $q\in \mathbb{S}^1$ such that 
\[r(p)+r (q) = \dist{p}{q}{\mathbb{S}^1}.\]
Therefore
\begin{align*}
\pi&=\dist{p}{(-p)}{\mathbb{S}^1}\le 
\\
&\le 
r(p)+r(-p)=
\\
&=
r(p)+r(q) +r(-p) -r(q)\le
\\
&\le
\dist{p}{q}{\mathbb{S}^1}+\dist{q}{(-p)}{\mathbb{S}^1}=
\\
&=\pi.
\end{align*}
So, we have equality in both places, and the only-if part follows.

\parit{If part.}
Assume $r$ is a 1-Lipschitz function such that $r(p)+r(-p)=\pi$.
Then 
\begin{align*}
\dist pq{\mathbb{S}^1}&=
\dist{p}{(-p)}{\mathbb{S}^1}-\dist{q}{(-p)}{\mathbb{S}^1}\ge
\\
&\ge\pi -(r(-p)-r(q))=
\\
&=r(p)+r(q).
\end{align*}
Therefore $r$ is admissible.

Finally, if $r$ is not extremal, then there is an admissible function $s\le r$ such that $s(p)<r(p)$ for some $p$.
The latter contradicts the equality $r(p)+r(-p)=\pi$.

\parit{Source:} \cite[Proposition 2.7]{zuest}.

\parbf{\ref{ex:retraction}.}
Show and use that
$s^*(x)\ge \dist{x}{y}{}-s(y)$
for any $x,y\in \spc{X}$.

\parit{Remarks.}
It is easy to check that $q\:s\z\mapsto \tfrac12\cdot(s+s^*)$ is a short map on the space of admissible functions (with sup-norm).
Moreover, iterating $q$ and passing to the limit, we get a short retraction from the space of admissible functions to the space of extremal functions on $\spc{X}$ \cite[see 3.1 in][]{lang-2013}.
The existence of such a map also follows from \ref{thm:inj-envelope}.

\parbf{\ref{ex:one-point-gluing}.}
Apply \ref{thm:injective=hyperconvex:balls}.

\parit{Comment.}
Conditions under which gluings of injective spaces are injective were studied by Benjamin Miesch and Maël Pavón \cite{miesch,miesch-pavon}.

\parbf{\ref{ex:Rm-ell-infty}.}
Let $B=\cBall[0,1]$ and $P\supset B$ be a parallelepiped of minimal volume.
Choose coordinates so that $P$ is described by inequalities
$|x_i|\le 1$ for all $i$;
let $e_1,\dots,e_m$ be the standard basis of these coordinates.

Let $B_i=\cBall[(1+R)\cdot e_i,R]$ for some $R>0$.
Show that $ e_i\in B$ for any $i$; in particular $B\cap B_i\ne\emptyset$.
Show $P$ can be chosen so that $B_i\cap B_j\ne \emptyset$ for all $i$ and~$j$ and all large $R>0$.
Apply hyperconvexity to show that $e_1+\dots+ e_m\in B$.
The same way, show that $\pm e_1\pm \dots\pm e_m\in B$ for all choices of signs.
Conclude that $B=P$.

\parbf{\ref{ex:compact-hyperconvex}.}
Observe that closed balls are compact and
apply the finite intersection property.

\parbf{\ref{ex:urysohn-hyperconvex}.}
Denote by $\spc{U}_d$ the $d$-Urysohn space,
so $\spc{U}_\infty$ is the Urysohn space.

The extension property implies finite hyperconvexity.
It remains to show that $\spc{U}_d$ is not countably hyperconvex.

Suppose that $d<\infty$. 
Show that for any point $x\in\spc{U}_d$ there is a point $y\in\spc{U}_d$ such that $\dist{x}{y}{\spc{U}_d}=d=\diam\spc{U}_d$.
It follows that there is no point $z\in\spc{U}_d$ such that $\dist{z}{x}{\spc{U}_d}\le \tfrac d2$ for any $x\in\spc{U}_d$.
Whence $\spc{U}_d$ is not countably hyperconvex.

Use \ref{ex:sphere-in-urysohn:midpoint} to reduce the case $d=\infty$ to the case $d<\infty$.

\parbf{\ref{ex:almost-hyperconvex}.}
Choose $\eps_0>0$.
Let $p_0$ be a point provided by the definition of almost hyperconvexity;
that is $\dist{x_\alpha}{p_0}{}\le r_\alpha+\eps_0$ for a given $\eps_0>0$.
We may assume that $\delta_0=\sup\{\,\dist{x_\alpha}{p_0}{}- r_\alpha\,\}>0$; otherwise, the problem is solved.
Clearly, $\delta_0\le \eps_0$.

Applying hyperconvexity for $\eps_1<\tfrac1{10}\cdot\delta_0$,
we get a point 
$p_1$ such that $\dist{x_\alpha}{p_1}{}\le r_\alpha+\eps_1$ and $\dist{p_0}{p_1}{}\le \delta_0+\eps_1$.
Again, we may assume that $\delta_1=\sup\{\,\dist{x_\alpha}{p_1}{}- r_\alpha,\dist{p_0}{p_1}{}\,\}>0$, and we have $\delta_1\le \eps_1$.

Continuing this way, we get a sequence $p_0,p_1,\dots$ that either terminates and in this case the problem is solved, or it is an infinite Cauchy sequence.
In the latter case, its limit $p_\infty$ satisfies $\dist{x_\alpha}{p_\infty}{}\le r_\alpha$ for any $\alpha$.

\parit{Comment.}
This solution reminds the proof of \ref{prop:completion-univeral};
a more exact statement was proved by Benjamin Miesch and Maël Pavón \cite[2.2]{miesch-pavon2016};
namely, they show that almost $n$-hyperconvexity implies $(n-1)$-hyperconvexity.

\parbf{\ref{ex:Inj(compact)}.}
Show and use that the functions in $\Inj\spc{X}$ are 1-Lipschitz and uniformly bounded.

\parbf{\ref{ex:tripod+square}}; \ref{SHORT.ex:tripod+square:2}.
Use \ref{lem:opposite-compact} to show that if $f$ is extremal if and only if $f(v)=x$ and $f(w)=1-x$ for some $x\in [0,1]$.
Conclude that $\Inj\spc{X}$ is isometric to the unit interval $[0,1]$.

\parit{\ref{SHORT.ex:tripod+square:tripod}.}
Let $f$ be an extremal function.
By \ref{lem:opposite-compact}, at least two of the numbers $f(a)+f(b)$, $f(b)+f(c)$, and $f(c)+f(a)$ are $1$.
It follows that for some $x\in[0,\tfrac12]$, we have 
\begin{align*}
f(a)&=1\pm x,&
f(b)&=1\pm x,&
f(c)&=1\pm x,
\end{align*}
where we have one ``minus'' and two ``pluses'' in these three formulas.

Suppose that
\begin{align*}
g(a)&=1\pm y,& g(b)&=1\pm y,& g(c)&=1\pm y
\end{align*}
is another extremal function.
Then $|f-g|\z=|x-y|$ if $g$ has ``minus'' at the same place as $f$ and $|f-g|=|x+y|$ otherwise.

It follows that $\Inj\spc{X}$ is isometric to a \textit{tripod} --- three segments of length $\tfrac12$ glued at one end.

\begin{Figure}
\begin{minipage}{.48\textwidth}
\centering
\includegraphics{mppics/pic-3}
\end{minipage}\hfill
\begin{minipage}{.48\textwidth}
\centering
\includegraphics{mppics/pic-4}
\end{minipage}
\vskip-4mm
\end{Figure}

\parit{\ref{SHORT.ex:tripod+square:square}.}
Assume $f$ is an extremal function.
Use \ref{lem:opposite-compact} to show that
\begin{align*}
2&=f(x)+f(y)=
\\
&=f(p)+f(q);
\end{align*}
in particular, two values $a=f(x)-1$ and $b\z=f(p)-1$ completely describe the function $f$.
Since $f$ is extremal, we also have that 
\[(1\pm a)+(1\pm b)\ge 1\]
for all 4 choices of signs;
equivalently, 
\[|a|+|b|\le 1.\]

It follows that $\Inj\spc{X}$ is isometric to the rhombus $|a|+|b|\le 1$ in the $(a,b)$-plane with the metric induced by the $\ell^\infty$-norm.

\parit{Remarks.}
If $\spc{X}$ has $n$-points, then (evidently) $\Inj\spc{X}$ is a polyhedral complex in $(\RR^n,\ell^\infty) \z=\ell^\infty(\spc{X})$;
each face of the complex is defined by equalities and inequalities of the following type: $x_i+x_j\ge \const$ and  $x_i+x_j= \const$.
It is easy to see (and follows from \ref{ex:Rm-ell-infty}) that each face is isometric to a convex polyhedron in  $(\RR^k,\ell^\infty)$ for some $k\le n$;
in fact $k\le n/2$.
The structure of the complex can be encoded by certain graphs with the vertex set $\spc{X}$ \cite[see Section 4 in][]{lang-2013}.

\parbf{\ref{ex:kur-inj}.}
Recall that $x\mapsto \distfun_x$ gives an isometric embedding $\spc{X}\z\hookrightarrow\ell^\infty(\spc{X})$;
so we can identify $\spc{X}$ with a subset of $\ell^\infty(\spc{X})$.
Further, $\Inj\spc{X}$ is a subset of $\ell^\infty(\spc{X})$.
It is sufficient to show that $\Inj\spc{X}=G$.

Use \ref{lem:opposite-compact} to show that $\Inj\spc{X}\subset G$.

Given $g\in G$, show that $g(x)=\dist{g}{x}{\ell^\infty(\spc{X})}$.
Conclude that $g$ is admissible and apply \ref{lem:opposite-compact}.

\parit{Source:} Suggested by Rostislav Matveyev.

\parbf{\ref{ex:4-on-a-line}.}
Recall that 
\[\dist{f}{g}{\Inj\spc{X}}=\sup\set{|f(x)-g(x)|}{x\in\spc{X}}\]
and 
\[\dist{f}{p}{\Inj\spc{X}}=f(p)\]
for any $f,g\in \Inj\spc{X}$ and $p\in \spc{X}$.

Since $\spc{X}$ is compact we can find a point $p\in\spc{X}$ such that 
\begin{align*}
\dist{f}{g}{\Inj\spc{X}}&=|f(p)-g(p)|=
\\
&=\left|\dist{f}{p}{\Inj\spc{X}}-\dist{g}{p}{\Inj\spc{X}}\right|.
\end{align*}
Without loss of generality, we may assume that 
\[\dist{f}{p}{\Inj\spc{X}}
=
\dist{g}{p}{\Inj\spc{X}}
+
\dist{f}{g}{\Inj\spc{X}}.\]
Applying \ref{lem:opposite-compact}, we can find a point $q\in\spc{X}$ such that 
\[\dist{q}{p}{\Inj\spc{X}}
=
\dist{f}{p}{\Inj\spc{X}}
+
\dist{f}{q}{\Inj\spc{X}},\]
whence the result.

Since $\Inj\spc{X}$ is injective (\ref{prop:InjX-is-injective}), by \ref{ex:inj=complete-geodesic-contractible:geodesic} it has to be geodesic. It remains to note that the concatenation of geodesics $[pq]$, $[gf]$, and $[fq]$ is the required geodesic $[pq]$.

\parbf{\ref{ex:delta-hyp}.} The only-if part follows since $\spc{X}$ is isometric to a subset of $\Inj\spc{X}$.

The if part means that 
\[
\begin{aligned}
\dist{f}{g}{}+\dist{v}{w}{}\le
\max\{\,
&\dist{f}{v}{}+\dist{g}{w}{},\\
\dist{f}{w}{}+&\dist{g}{v}{}
\,\}+2\cdot\delta
\end{aligned}
\eqlbl{eq:fgvw-hyp}\]
for any $f,g,v,w\in \Inj\spc{X}$.

Suppose $\spc{X}$ is compact. 
Applying \ref{ex:4-on-a-line}, we can choose $p,q,x,y\in \spc{X}$  such that 
\[
\begin{aligned}
\dist{p}{f}{}+\dist{f}{g}{}+\dist{g}{q}{}&=\dist{p}{q}{}
\\
\dist{x}{v}{}+\dist{v}{w}{}+\dist{w}{y}{}&=\dist{x}{y}{}
\end{aligned}
\eqlbl{eq:pfgq+xvwy}
\]

Since $\spc{X}$ is $\delta$-hyperbolic, we have
\[\begin{aligned}
\dist{p}{q}{}+\dist{x}{y}{}\le
\max\{\,&\dist{p}{x}{}+\dist{q}{y}{},
\\
\dist{p}{y}{}+&\dist{q}{x}{}\,\}+2\cdot\delta.
\end{aligned}\]
Show that this inequality, together with the triangle inequality and \ref{eq:pfgq+xvwy} imply \ref{eq:fgvw-hyp}.

For the noncompact case, prove an approximate version of \ref{eq:pfgq+xvwy} and apply it the same way.

\parbf{\ref{ex:inj-envelope}.}
Show that there is a unique isometry of $\Inj\spc{X}$ that is identity on $\spc{X}$.
Use it together with \ref{thm:inj-envelope}.

\parbf{\ref{ex:d-p-inclusion}.}
Show that there is a pair of short maps 
$\Inj\spc{X}\to\Inj\spc{U}\to\Inj\spc{X}$ 
such that their composition is the identity on $\spc{X}$.
You may need to apply the Katětov extension (\ref{ex:extension-of-extension:a}).
Make a conclusion.

\parbf{\ref{ex:hemisphere-inj}.}
Apply \ref{lem:opposite-compact} to show that for any $u\in\mathbb{S}^2_+$ the restriction $f_u\z\df\distfun_u|_{\mathbb{S}^1}$ is an extremal function on $\mathbb{S}^1$.
Moreover, the function $f_u$ uniquely determines $u$. 
Make a conclusion.

\parbf{\ref{ex:3-4-hypreconvex}.}
Observe that coordinate functions are monotonic on any geodesic in $\ell^1$.
Use it to show that $\ell^1$ is a \emph{median space};
that is, for any three points $x,y,z$ there is a \textit{unique} point $m$ (it is called the \index{median}\emph{median} of $x$, $y$, and $z$) that lies on {}\emph{some} geodesics $[xy]$, $[xz]$ and $[yz]$.
Apply it to show that $\ell^1$ is 3-hyperconvex.

The 4-hyperconvexity fails for the unit balls centered at four even vertices of the cube $([0,1]^3,\ell^1)$.

\parbf{\ref{ex:ultrametric}.}
Choose three points $x,y,z\in\spc{X}$ and set $\spc{A}=\{x,z\}$.
Let $f\:\spc{A}\z\to \spc{A}$ be the identity map.
Then $F(y)=x$ or $F(y)=z$.
In both cases,
the strong triangle inequality follows.

\parbf{\ref{ex:ultrametric-converse}}; \textit{main part.}
Choose a maximal (with respect to inclusion) subset $A\z\supset \spc{K}$ that admits a short retraction $f\:A\to \spc{K}$;
it exists by Zorn's lemma.
If $A$ is the whole space, then the problem is solved.
Otherwise, choose $p\notin A$.

Choose a sequence of points $a_n\in A$ such that $\dist{a_n}{p}{}$ converge to the exact lower bound on the distances from points in $A$ to $p$.
Since $\spc{K}$ is compact, we can pass to a subsequence of $a_n$ such that $f(a_n)$ converges.
Let 
\[f(p)=\lim f(a_n).\]

It remains to check that 
\[\dist{f(a)}{f(p)}{}\le\dist{a}{p}{}\eqlbl{eq:short-retract}\]
for any $a\in A$.
Choose $\eps>0$; note that 
\begin{align*}
\dist{a_n}{p}{}&<\dist{a}{p}{}+\eps,
&
\dist{f(a_n)}{f(p)}{}&<\eps
\end{align*}
for all large~$n$.
Therefore, 
\begin{align*}
&\!\!\!\!\dist{f(a)}{f(p)}{}\le
\\
&\le\max\{\,\dist{f(a)}{f(a_n)}{},
\dist{f(a_n)}{f(p)}{}\,\}\le
\\
&\le \dist{f(a)}{f(a_n)}{}+\eps\le
\\
&\le \dist{a}{a_n}{} +\eps\le 
\\
&\le \max\{\,\dist{a}{p}{},\dist{a_n}{p}{}\,\}+\eps< 
\\
&< \dist{a}{p}{}+2\cdot\eps.
\end{align*}
Since $\eps>0$ is arbitrary, we get \ref{eq:short-retract}.

\parit{Example.}
Consider set of $\{\infty,1,2,\dots\}$ with metric defined by 
\[|m-n|
\df
1+\frac1{\min\{m,n\}}\]
for $m\ne n$.
Observe that the space is complete, the subset $\{1,2,\dots\}$ is closed, but it is not a short retract of the ambient space.

\parbf{\ref{ex:petrunin-stadler}.} Consider the space $\spc{K}^{\spc{X}}$ of all maps $\spc{X}\z\to \spc{K}$ equipped with the product topology.

Denote by $\mathfrak{S}_F$ the set of maps $h\in \spc{K}^\spc{X}$ such that the restriction $h|_F$  is short and agrees with $f$ in $F\cap A$.
Note that the sets $\mathfrak{S}_F\subset \spc{K}^\spc{X}$ are closed and any finite intersection of these sets is nonempty.

According to Tikhonov's theorem, $\spc{K}^{\spc{X}}$ is compact.
By the finite intersection property, the intersection $\bigcap_F\mathfrak{S}_F$ for all finite sets $F\subset X$ is nonempty.
Hence the statement follows.

\parit{Source:} \cite[7.1]{petrunin-stadler}.

%% file: converge-sol.tex
\refstepcounter{chapter}
\setcounter{eqtn}{0}

\parbf{\ref{ex:diam}.}
Suppose that $\dist{A}{B}{\Haus\spc{X}}<r$.
Choose a pair of points $a,a'\in A$ on maximal distance from each other.
Observe that there are points $b,b'\in B$ such that 
$\dist{a}{b}{\spc{X}},\dist{a'}{b'}{\spc{X}}<r$.
Whence 
\[\dist{a}{a'}{\spc{X}}-\dist{b}{b'}{\spc{X}}\le 2\cdot r\]
and therefore
\[\diam A-\diam B\le 2\cdot\dist{A}{B}{\Haus\spc{X}}.\]

Swap $A$ and $B$ and repeat the argument.

\parbf{\ref{ex:Hausdorff-bry}}; \ref{SHORT.ex:Hausdorff-bry:conv}.
Given a set $A\z\subset \RR^2$, denote by $A^r$ its closed $r$-neighborhood.
Show and use that 
\[(\Conv A)^r=\Conv(A^r).\]

\parit{\ref{SHORT.ex:Hausdorff-bry:bry}.}
The answer is ``no'' in both parts.

For the first part let $A$ be a unit disk and $B$ a finite $\eps$-net in $A$.
Evidently, $|A-B|_{\Haus\RR^2}\z<\eps$, 
but
$|\partial A-\partial B|_{\Haus\RR^2}\approx 1$.

For the second part take $A$ to be a unit disk and $B=\partial A$ to be its boundary circle.
Note that $\partial A=\partial B$; in particular, $\dist{\partial A}{\partial B}{\Haus\RR^2}=0$ while $\dist{ A}{B}{\Haus\RR^2}=1$.

\parit{Remark.}
A more interesting example for \ref{SHORT.ex:Hausdorff-bry:bry} is provided by the so-called \textit{lakes of Wada} --- an example of three (and more) disjoint open topological disks in the plane that have identical boundaries.

\parbf{\ref{ex:Haus-func}.}
Checking two functions $\distfun_A$ and $\distfun_B$ leads to 
\[\dist{A}{B}{}\le\sup_f\, \{\,\max_{a\in A}\{f(a)\}-\max_{b\in B}\{f(b)\,\}.\]
Use \ref{obs:Haus-nbhds} to prove the opposite inequality.

\parbf{\ref{ex:Haus-support}.} By \ref{ex:Hausdorff-bry:conv}, it is sufficient to show that  
\[\dist{A}{B}{\Haus\RR^n}= \sup_{|u|=1}\{|h_A(u)-h_B(u)|\}\]
for any nonempty compact convex sets $A,B\z\subset \RR^n$.

Prove the 1-dimensional case of this equality.
Further, denote by $A_\ell$ the orthogonal projection of $A$ to a line $\ell$.
Show and use that
\[\dist{A}{B}{\Haus\RR^n}=\sup_\ell\{\dist{A_\ell}{B_\ell}{\Haus\ell}\},\]
where the least upper bound is taken for all lines~$\ell$.

\parbf{\ref{ex:H-sections}};
\ref{SHORT.ex:H-sections:S}.
Given $t\in (0,1]$, consider the real interval $\tilde C_t=[\tfrac 1t+t, \tfrac 1t+1]$.
Denote by $C_t$ the image of $\tilde C_t$ under the covering map $\pi\:\RR\to \mathbb{S}^1\z=\RR/\ZZ$.

Set $C_0=\mathbb{S}^1$.
Note that the Hausdorff distance from $C_0$ to $C_t$ is $\tfrac t2$.
Therefore $\{C_t\}_{t\in[0,1]}$ is a family of compact subsets in $\mathbb{S}^1$ that is continuous in the sense of Hausdorff.

Assume there is a continuous section $c(t)\z\in C_t$, for $t\in [0,1]$.
Since $\pi$ is a covering map,
we can lift the path $c$ to a path $\tilde c\:[0,1]\to \RR$ such that $\tilde c(t)\in \tilde C_t$ for all $t$.
In particular, $\tilde c(t)\to\infty$ as $t\to0$,
a contradiction.

\parit{\ref{SHORT.ex:H-sections:R}.}
Consider path $c(t)\df\min C_t$.

\parit{Source:} Suggested by Stephan Stadler.

\parbf{\ref{ex:haus-contractible}.}
Show that $f\:(p,r)\mapsto \cBall[p,r]$ defines a continuous map $\spc{X}\z\times[0,\infty)\to \Haus\spc{X}$.
Observe that $\cBall[p,r]=\spc{X}$ for all $r\ge\diam\spc{X}$.
Composing $f$ with the retraction $r\:\Haus\spc{X}\to \spc{X}$ we get a homotopy from the identity map to a constant map on $\spc{X}$, hence the statement.

\parit{By the way, is there a good description of such spaces?}
Note that if $\spc{X}$ is injective or discrete, then a short retraction $\Haus\spc{X}\to \spc{X}$ exists.
On the other hand, the euclidean plane does not have such a retraction \cite{petrunin-2022-MO-center, przeslawski-yost}.
See also \cite{petrunin-435157}.

\parbf{\ref{ex:closure-union}.} 
Show that for any $\eps>0$ there is a positive integer $N$ such that $\bigcup_{n\le N} K_n$ is an $\eps$-net in the union $\bigcup_{n} K_n$.
Observe that $\bigcup_{n\le N} K_n$
is compact.
Apply \ref{ex:compact-net} and \ref{totally-bounded}.

\parbf{\ref{ex:Haus-length}}; \textit{if part.}
Choose two compact sets $A,B\z\subset \spc{X}$;
suppose that $\dist{A}{B}{\Haus\spc{X}}<r$.

Choose finite $\eps$-nets $\{a_1,\dots,a_m\}\subset A$ and $\{b_1,\dots b_n\}\subset B$.
For each pair $a_i,b_j$ construct a constant-speed path $\gamma_{i,j}$ from $a_i$ to $b_j$ such that 
\[\length \gamma_{i,j}<\dist{a_i}{b_j}{}+\eps.\]
Set 
\[C(t)=\set{\gamma_{i,j}(t)}{\dist{a_i}{b_j}{\spc{X}}<r+\eps}.\]
Observe that $C(t)$ is finite; in particular, it is compact.

Show and use that 
\begin{align*}
\dist{A}{C(t)}{\spc{X}}&<t\cdot r+10\cdot\eps,
\\
\dist{C(t)}{B}{\spc{X}}&<(1-t)\cdot r+10\cdot\eps.
\end{align*}
Apply \ref{ex:closure-union} and \ref{lem:mid>geod}.

\parit{Only-if part.}
Choose points $p,q\in\spc{X}$. 
Show that the existence of $\eps$-midpoints between $\{p\}$ and $\{q\}$ in $\Haus\spc{X}$ implies the existence of $\eps$-midpoints between $p$ and $q$ in $\spc{X}$.
Apply \ref{lem:mid>geod}.

\parbf{\ref{ex:Haus-G-delta}}; \ref{SHORT.ex:Haus-G-delta:closed}
Suppose that a sequence of compact subsets $K_n\subset \RR^2$ converges to $K_\infty$ in the sense of Hausdorff.
Assume $K_\infty$ is not connected; show that so is  $K_n$ for large~$n$.

\parit{\ref{SHORT.ex:Haus-G-delta:curves}.}
Choose a finite subset $F$ that is $\eps$-close to~$K$.
Show that one can obtain a tree $T$ by connecting some vertices of $F$ by line segments of length smaller than $2\cdot\eps$.
Let $\gamma$ be a curve that bounds a neighborhood of $T$.
Show that if the neighborhood is sufficiently small, then $\gamma$ is $2\cdot\eps$-close to~$K$.

\parit{Remarks.}
You might be surprised to learn that most connected compact sets in the pane are homeomorphic to each other --- they are homeomorphic to the so-called \index{pseudo-arc}\emph{pseudo-arc} \cite{bing-1948}; here the word \textit{most} understood in the sense of \ref{sec:G-delta}.
In particular, most of the compact connected plane sets are not path-connected.

\parbf{\ref{ex:Huas-perimeter-area}.}
Let $A$ be a compact convex set in the plane.
Denote by $A^r$ the closed $r$-neighborhood of $A$.
Recall that by Steiner's formula we have
\[\area A^r=\area A+r\cdot\perim A+\pi\cdot r^2.\]
Taking the derivative and applying the coarea formula, we get 
\[\perim A^r=\perim A+2\cdot\pi\cdot r.\]

Observe that if $A$ lies in a compact set $B$ bounded by a closed curve, then 
\[\perim A\le \perim B.\]
Indeed the closest-point projection $\RR^2\to A$ is short and it maps $\partial B$ onto $\partial A$.

It remains to use the following observation: if $A_n\to A_\infty$, then for any $r>0$ we have that the inclusions
\[A_\infty^r\supset A_n
\quad\text{and}\quad
A_\infty\subset A_n^r\]
hold for all large $n$.

\parbf{\ref{ex:round-disc}.}
Note that almost all points on $\partial D$ have a defined tangent line.
In particular, for almost all pairs of points $a,b\z\in\partial D$ the two angles $\alpha$ and $\beta$ between the chord $[ab]$ and $ \partial D$ are defined.

\begin{wrapfigure}{r}{27 mm}
\vskip-3mm
\centering
\includegraphics{mppics/pic-410}
\end{wrapfigure}

The convexity of $D'$ implies that $\alpha=\beta$;
here we measure the angles $\alpha$ and $\beta$ on one side from $[ab]$.
Show that if the identity $\alpha=\beta$ holds for almost all chords, then $D$ is a round disk.

\parbf{\ref{ex:generalized-selection}.}
Observe that all functions $\distfun_{A_n}$ are Lipschitz.
Suppose that for some (and therefore any) point $x$ the sequence $\distfun_{A_n}(x)$ is not bounded.
Then we can pass to a subsequence of $A_n$ so that $\distfun_{A_n}(x)\to\infty$ for any $x$;
in this case, $A_n$ converges to the empty set.

Assume the sequence $\distfun_{A_n}(x)$ is bounded for some (and therefore any) point $x$.
Then, passing to a subsequence of $A_n$, we may assume that the sequence $\distfun_{A_n}$ converges to some function $f$.

Set $A_\infty=f^{-1}\{0\}$.
It remains to show that $f=\distfun_{A_\infty}$.


\parbf{\ref{ex:d_GH-and-diam}};
\ref{SHORT.ex:d_GH-and-diam:point}.
Apply the definition (\ref{def:GH}) for space $\spc{W}$ obtained from $\spc{X}$ by adding one point that lies at distance $\tfrac12\cdot\diam \spc{X}$ from each point of $\spc{X}$.

\parit{\ref{SHORT.ex:d_GH-and-diam:scale}.}
Given a point $x\in\spc{X}$, denote by $a\cdot x$ and $b\cdot x$ the corresponding points in $a\cdot\spc{X}$ and $b\cdot \spc{X}$ respectively.
Show that there is a metric on $\spc{W}\z=a\cdot\spc{X}\sqcup b\cdot\spc{X}$ such that 
\[|a\cdot x-b\cdot x|_{\spc{W}}=\tfrac{|b-a|}2\cdot\diam\spc{X}\]
for any $x$ and the inclusions
$a\cdot\spc{X}\z\hookrightarrow\spc{W}$,
$b\cdot\spc{X}\z\hookrightarrow\spc{W}$ are distance-preserving.
Conclude that 
$|\spc{X}-\spc{O}|_{\GH}\le\tfrac12\cdot \diam \spc{X}.$
The opposite inequality follows from \ref{SHORT.ex:d_GH-and-diam:point}.

\parit{\ref{SHORT.ex:d_GH-and-diam:isometry}.}
Use \ref{SHORT.ex:d_GH-and-diam:point} and \ref{SHORT.ex:d_GH-and-diam:scale} to show that the isometry class of $\spc{O}$ is completely determined by the following property
\[|\spc{X}-\spc{Y}|_{\GH} \le \max\{\,|\spc{O}-\spc{X}|_{\GH},|\spc{O}-\spc{Y}|_{\GH}\,\}.\]
for any $\spc{X}$ and $\spc{Y}$.

\parit{Remark.}
In fact, \textit{the isometry group of space $\GH$ is trivial}.
The latter was proved by George Lowther \cite{lowther, ivanov-tuzhilin}.

\parbf{\ref{ex:GH<H}.}
Check a one-point set and the vertices of an equilateral triangle.
You may use \ref{ex:d_GH-and-diam:point}.

\parbf{\ref{ex:rectangle}.}
Suppose that
we can identify $\spc{A}_r$ and $\spc{B}_r$ with subspaces of a space $\spc{W}$
such that 
\[|\spc{A}_r-\spc{B}_r|_{\Haus \spc{W}}<\tfrac1{10}\]
for large $r$; see the definition of Gromov--Hausdorff metric (\ref{def:GH}).

Set $n=\lceil r \rceil$.
Note that there are $2\cdot n$ integer points in~$\spc{A}_r$: 
$a_1\z=(0,0)$, $a_2\z=(1,0),\ \dots,\ a_{2\cdot n}\z=(n,1)$.
Choose a point $b_i\in \spc{B}_r$ that lies at the minimal distance from $a_i$.
Note that $|b_i-b_j|>\tfrac 45$ if $i\ne j$.
It follows that $r>\tfrac 45\cdot (2\cdot n-1)$.
The latter contradicts $n=\lceil r \rceil$ for large~$r$.

\parit{Remark.}
Try to show that $|\spc{A}_r-\spc{B}_r|_{\GH}=\tfrac12$ for all large $r$.

\refstepcounter{chapter}
\setcounter{eqtn}{0}

\parbf{\ref{ex:GH-inj}.}
Suppose 
\[|\spc{X}-\spc{Y}|_{\Haus\spc{U}}<\eps.\]
Denote by $\hat{\spc{U}}$ the injective envelope of $\spc{U}$.
According to \ref{ex:d-p-inclusion}, the inclusions $\spc{X}\hookrightarrow\spc{U}$ and $\spc{Y}\hookrightarrow\spc{U}$ can be extended to distance-preserving inclusions $\hat{\spc{X}}\hookrightarrow\hat{\spc{U}}$ and $\hat{\spc{Y}}\hookrightarrow\hat{\spc{U}}$.
Therefore, we can and will consider  $\hat{\spc{X}}$ and $\hat{\spc{Y}}$ as subspaces of $\hat{\spc{U}}$.
It is sufficient to show that
\[|\hat{\spc{X}}-\hat{\spc{Y}}|_{\Haus\hat{\spc{U}}}<2\cdot \eps.\eqlbl{eq:hats|X-Y|<2eps}\]

Given $f\in \hat{\spc{U}}$,
let us find $g\in \hat{\spc{X}}$ such that 
\[|f(u)-g(u)|<2\cdot\eps\eqlbl{|g-f|}\]
for any $u\in\spc{U}$.
Note that the restriction $f|_{\spc{X}}$ is admissible on ${\spc{X}}$.
By \ref{obs:extremal:below}, there is $g\in \hat{\spc{X}}$ such that 
\[g(x)\le f(x)\eqlbl{g(x)=<f(x)}\]
for any $x\in\spc{X}$.

Recall that any extremal function is $1$-Lipschitz;
in particular, $f$ and $g$ are $1$-Lipschitz on $\spc{U}$.
Therefore, \ref{g(x)=<f(x)} and $|\spc{X}-\spc{Y}|_{\spc{U}}<\eps$ imply that
\[g(u)< f(u)+2\cdot \eps\]
for any $u\in\spc{U}$.
By \ref{lem:+-c}, we also have 
\[g(u)> f(u)-2\cdot \eps\]
for any $u\in\spc{U}$.
Whence \ref{|g-f|} follows.

It follows that $\hat{\spc{Y}}$ lies in a $2\cdot\eps$-neighborhood of $\hat{\spc{X}}$ in $\hat{\spc{U}}$.
The same way we show that $\hat{\spc{X}}$ lies in a $2\cdot\eps$-neighborhood of $\hat{\spc{Y}}$ in $\hat{\spc{U}}$.
Hence \ref{eq:hats|X-Y|<2eps} follows.

\parit{Remark.} 
This problem was discussed by Urs Lang, Maël Pavón, and Roger Züst \cite[3.1]{lang-pavon-zust}.
They also show that the constant 2 is optimal.
To see this, look at the injective envelopes of two four-point metric spaces shown on the diagram and observe that the Gromov--Hausdorff distance between the 4-point metric spaces is 1, while the distance between their injective envelopes approaches 2 as $s\to\infty$. 

\begin{Figure}
\vskip-0mm
\centering
\includegraphics{mppics/pic-505}
\end{Figure}

\parbf{\ref{ex:H-R}}; \textit{only-if part.}
Let us identify $\spc{X}$ and $\spc{Y}$ with subspaces of a metric space $\spc{W}$ such that 
\[|\spc{X}-\spc{Y}|_{\Haus \spc{W}}<\eps.\]

Set $x\approx y$ if and only if $\dist{x}{y}{\spc{W}}<\eps$.
It remains to check that $\approx$ is an $\eps$-approximation.

\parit{If part.}
Show that we can assume that 
\[R=\set{(x,y)\in\spc{X}\times\spc{Y}}{x\approx y}\] is a compact subset of $\spc{X}\times\spc{Y}$.
Conclude that
\[\bigl|\dist{x}{x'}{\spc{X}}-\dist{y}{y'}{\spc{Y}}\bigr|<2\cdot\eps'\]
for some $\eps'<\eps$.

Show that there is a metric on $\spc{W}=\spc{X}\sqcup\spc{Y}$ such that the inclusions $\spc{X}\hookrightarrow\spc{W}$ and
$\spc{Y}\hookrightarrow\spc{W}$ are distance-preserving and $\dist{x}{y}{\spc{W}}=\eps'$ if $x\approx y$.
Conclude that 
\[|\spc{X}-\spc{Y}|_{\Haus \spc{W}}\le\eps'<\eps.\]

\parbf{\ref{ex:eps-isom}};
\ref{SHORT.ex:eps-isom:GH>isom}.
Let $\approx$ be an $\eps$-approximation provided by \ref{ex:H-R}.
For any $x\in\spc{X}$ choose a point $f(x)\in\spc{Y}$ such that $x\approx f(x)$.
Show that $x\mapsto f(x)$ is an $2\cdot\eps$-isometry.

\parit{\ref{SHORT.ex:eps-isom:isom>GH}.}
Let $x\in\spc{X}$ and $y\in\spc{Y}$.
Set $x\approx y$ if $\dist{y}{f(x)}{\spc{Y}}<\eps$.
Show that $\approx$ is an $\eps$-approximation. 
Apply \ref{ex:H-R}.

\parbf{\ref{ex:XYZ}.}
Consider the product space $[0,\eps]\times \ZZ_n$ with the natural $\ell^\infty$-product metric. 
Make three variations of it by changing the sizes of some segments.

\parbf{\ref{ex:GH-SC}}; \ref{SHORT.ex:GH-SC:circle}.
Suppose $\spc{X}_n$ are simply-connected length metric space, $\spc{X}_n\GHto \spc{X}$, and there is a nontrivial covering map $f\:\tilde{\spc{X}}\to \spc{X}$.
We will arrive at a contradiction by showing that there is a nontrivial covering map $f_n\:\tilde{\spc{X}}_n\to \spc{X}_n$ for large $n$.

Choose a base point $p\in \spc{X}$ and its inverse image $\tilde p\in \tilde{\spc{X}}$.
Consider two paths $\alpha,\alpha'\:[0,1]\to \spc{X}$ that start at $p$;
denote by $\tilde \alpha,\tilde \alpha'\:[0,1]\z\to \tilde{\spc{X}}$ their liftings.
Show that there is $\eps>0$ such that if $\dist{\alpha(t)}{\alpha'(t)}{\spc{X}}<\eps$ for any $t$, then $\dist{\tilde \alpha(1)}{\tilde \alpha'(1)}{\tilde{\spc{X}}}<\eps$.

Now suppose $n$ is large.
Choose an $\tfrac\eps{10}$-approximation $\approx$ for $\spc{X}_n$ and $\spc{X}$.
Choose $q\in\spc{X}_n$ such that $q\approx p$.
Show that for any path $\beta\:[0,1]\to \spc{X}_n$ that starts at $q$ there is a path $\alpha\:[0,1]\to \spc{X}$ that starts at $p$ such that $\alpha(t)\z\approx\beta(t)$ for any $t$.
Observe that if $\alpha$ and $\alpha'$ are two choices of such paths, then $\dist{\alpha(t)}{\alpha'(t)}{\spc{X}}<\eps$.

Mimicking the standard construction of a covering map, we get the needed $f_n\:\tilde{\spc{X}}_n\to \spc{X}_n$.



\parit{\ref{SHORT.ex:GH-SC:nonsc-limit}.}
Let $\spc{V}$ be a cone over Hawaiian earrings.
Consider the \textit{doubled cone} $\spc{W}$ --- two copies of $\spc{V}$ with glued base points (see the diagram).

\begin{Figure}
\vskip-0mm
\centering
\includegraphics{mppics/pic-2}
\end{Figure}

The space $\spc{W}$ can be equipped with a length metric
(for example, the induced length metric from the shown embedding).

Show that $\spc{V}$ is simply-connected, but $\spc{W}$ is not; use the van Kampen theorem.

If we delete from the earrings all small circles and repeat the construction,
then the obtained double cone becomes simply-connected and remains close to $\spc{W}$.
That is, $\spc{W}$ is a Gromov--Hausdorff limit of simply-connected spaces.

\parit{Remark.}
Note that the limit space in \ref{SHORT.ex:GH-SC:nonsc-limit}, does not admit a nontrivial covering.

\parbf{\ref{ex:sphere-to-ball}};
\textit{\ref{SHORT.ex:sphere-to-ball:2}.}
Suppose that a metric on $\mathbb{S}^2$ is close to the unit disk $\DD^2$.
Show that $\mathbb{S}^2$ contains a circle $\gamma$ that is close to the boundary curve of $\DD^2$.
By the Jordan curve theorem, $\gamma$ cuts $\mathbb{S}^2$ into two disks, say $D_1$ and $D_2$.

By \ref{ex:GH-SC:circle}, the Gromov--Hausdorff limits of $D_1$ and $D_2$ have to contain the whole $\DD^2$; otherwise, the limit would admit a nontrivial covering.

Consider points $p_1\in D_1$ and $p_2\in D_2$ that are close to the center of $\DD^2$.
On one hand, the distance $\dist{p_1}{p_2}{\mathbb{S}^2}$ has to be small.
On the other hand, any curve from $p_1$ to $p_2$ must cross $\gamma$, so its length is about 2 (or larger) --- a contradiction.

\parit{\ref{SHORT.ex:sphere-to-ball:3}.}
Show that one can remove fine tunnels from the standard 3-ball in such a way that (1) the topology does not change, (2) the induced length metric is very close to the original one, and (3) the tunnels come sufficiently close to any point in the ball.

Consider the \index{doubling}\emph{doubling} of the obtained ball along its boundary;
that is, two copies of the ball with glued corresponding points on their boundaries.
The obtained space is homeomorphic to $\mathbb{S}^3$.
Observe that the obtained space is sufficiently close to the original ball.

\parit{Source:} \cite[Exercises 7.5.13 and 7.5.17]{burago-burago-ivanov}. 

\parbf{\ref{ex:utb+pack}.} Apply \ref{ex:pack-net}.

\parbf{\ref{pr:doubling}.}
Let $\mu$ be a $C$-doubling measure on a space $\spc{X}$ from $\bm{Q}(C,D)$.
Without loss of generality, we may assume that $\mu(\spc{X})\z=1$.

The doubling condition implies that 
\[\mu[\oBall(p,\tfrac D{2^n})]\ge\tfrac 1{C^n}\]
for any point $x\in \spc{X}$.
It follows that 
\[\pack_{\frac D{2^n}}\spc{X}\le C^n.\]

By \ref{ex:pack-net}, for any $\eps\ge\frac D{2^{n-1}}$, the space $\spc{X}$ admits an $\eps$-net with at most $C^n$ points.
Whence $\bm{Q}(C,D)$ is uniformly totally bounded.

\parbf{\ref{pr:under}}; \ref{SHORT.pr:under:if}.
Choose $\eps>0$.
Since $\spc{Y}$ is compact, we can choose a finite $\eps$-net $\{y_1,\dots,y_{n}\}$ in $\spc{Y}$.

Suppose $f\:\spc{X}\to \spc{Y}$ be a distance-noncontracting map.
Choose one point $x_i$ in each nonempty subset $B_i=f^{-1}[\oBall(y_i,\eps)]$.
Note that the subset $B_i$ has diameter at most $2\cdot \eps$ and 
\[\spc{X}=\bigcup_iB_i.\]
Therefore, the set of points $\{x_i\}$ is a $2\cdot\eps$-net in~$\spc{X}$.

\parit{\ref{SHORT.pr:under:only-if}.}
Let $\bm{Q}$ be a uniformly totally bounded family of spaces. 
Suppose that each space in $\bm{Q}$ has an $\tfrac1{2^n}$-net with at most $M_n$ points; we may assume that $M_0=1$.

Consider the space $\spc{Y}$ of all infinite integer sequences $m_0,m_1,\dots$ such that $1\le m_n\le M_n$ for any $n$.
Given two sequences $\bm{\ell}\z=(\ell_1,\ell_2,\dots)$, and $\bm{m}\z=(m_1,m_2,\dots)$ of points in $\spc{Y}$, set 
\[\dist{\bm{\ell}}{\bm{m}}{\spc{Y}}=\tfrac C{2^{n}},\]
where $n$ is the minimal index such that $\ell_n\ne m_n$ and $C$ is a positive constant.

Observe that $\spc{Y}$ is compact.
Indeed it is complete and the sequences with constant tails, starting from index $n$, form a finite $\tfrac C{2^{n}}$-net in~$\spc{Y}$.

Given a space $\spc{X}$ in $\bm{Q}$,
choose a sequence of $\tfrac1{2^n}$ nets 
$N_n\subset\spc{X}$ for each $n$.
We can assume that $|N_n|\le M_n$; let us label the points in $N_n$ by $\{1,\dots,M_n\}$.
Consider the map $f\:\spc{X}\to\spc{Y}$ defined by $f:x\z\mapsto (m_1(x),m_2(x),\dots)$ where $m_n(x)$ is the label of a point in $N_n$ that lies at the distance $<\tfrac1{2^n}$ from $x$.

If $\tfrac1{2^{n-2}}\ge \dist{x}{x'}{\spc{X}}>\tfrac1{2^{n-1}}$, then $m_n(x)\z\ne m_n(x')$.
It follows that $\dist{f(x)}{f(x')}{\spc{Y}}\ge \tfrac C{2^{n}}$.
In particular, if $C>10$, then 
\[\dist{f(x)}{f(x')}{\spc{Y}}\ge \dist{x}{x'}{\spc{X}}\]
for any $x,x'\in \spc{X}$.
That is, $f$ is a distance-noncontracting map $\spc{X}\to \spc{Y}$.

\parbf{\ref{ex:GH-G-delta}.}
Let $\spc{K}$ be a compact space.
Denote by $a(\spc{K})$ the largest diameter of connected component in a compact space $\spc{K}$.
Further, let 
\[b(\spc{K})=\max_{p\in \spc{K}}\min_{q\ne p}\{\dist{p}{q}{\spc{K}}\}.\]
Note that $b(\spc{K})=0$ if and only if $\spc{K}$ has no isolated points.

Show that if $a(\spc{K})=b(\spc{K})=0$, then $\spc{K}$ is homeomorphic to the Cantor set.

Further show that the sets
\begin{align*}
A_\eps&=\set{\spc{K}\in\GH}{a(\spc{K})<\eps}
\\
B_\eps&=\set{\spc{K}\in\GH}{b(\spc{K})<\eps}
\end{align*}
are open and dense in $\GH$.
Apply \ref{thm:baire}.

\parbf{\ref{ex-GH-length}}; \ref{SHORT.ex-GH-length:separable}.
Show that (the isometry classes of) finite metric spaces with only rational distances form a countable dense subset in $\GH$. 

\parit{\ref{SHORT.ex-GH-length:length}+\ref{SHORT.ex-GH-length:geodesic}.}
Choose two compact metric spaces $\spc{X}$ and $\spc{Y}$.
Let $\spc{W}\supset \spc{X}',\spc{Y}'$ be as in \ref{prop:GH=H};
so $\spc{X}'\iso \spc{X}$, $\spc{Y}'\iso \spc{Y}$, and
\[\ell=\dist{\spc{X}'}{\spc{Y}'}{\Haus\spc{W}}=\dist{\spc{X}}{\spc{Y}}{\GH}\]
for some $\ell\ge 0$.

We can assume that $\spc{W}=\spc{X}'\cup\spc{Y}'$, so $\spc{W}$ is compact.
Choose a compact geodesic extension $\spc{G}$ of $\spc{W}$;
it exists by \ref{ex:compact-length} (or \ref{ex:Inj(compact)}).

Given $t\in[0,\ell]$, consider the set
\[\spc{Z}_t=\set{w\in \spc{G}}{\distfun_{\spc{X}'}w\le t,\  \distfun_{\spc{Y}'}w\le \ell-t}.\]
Observe that $t\mapsto \spc{Z}_t$ is a geodesic in $\Haus(\spc{G})$ from $\spc{X}'$ to $\spc{Y}'$.
Conclude that $t\mapsto [\spc{Z}_t]$ is a geodesic in $\GH$ from $[\spc{X}]$ to $[\spc{Y}]$.

\parit{Source:} \cite{ivanov-nikolaeva-tuzhilin}.

\parbf{\ref{ex:GH-po}}; \ref{SHORT.ex:GH-po:a}.
To check that $\dist{*}{*}{\GH'}$ is a metric, it is sufficient to show that
\[\dist{\spc{X}}{\spc{Y}}{\GH'}=0 
\quad\Longrightarrow\quad
\spc{X}\iso\spc{Y};\]
the remaining conditions are trivial.

If $\dist{\spc{X}}{\spc{Y}}{\GH'}=0$, then there is a sequence of maps $f_n\:\spc{X}\to \spc{Y}$ such that 
\[\dist{f_n(x)}{f_n(x')}{\spc{Y}}\ge \dist{x}{x'}{\spc{X}}-\tfrac1n.\]

Choose a countable dense subset $S\subset \spc{X}$ and pass to a subsequence such that $f_n(x)$ converges for any $x\in S$; denote by $f_\infty\:S\z\to \spc{Y}$ the limit map.
Note that $f_\infty$ is distance-noncontracting, and it can be extended to a distance-noncontracting map $f_\infty\:\spc{X}\to \spc{Y}$.

The same way we can construct a distance-noncontracting map 
$g_\infty\:\spc{Y}\to \spc{X}$.

By \ref{ex:non-contracting-map}, the compositions $f_\infty\circ g_\infty\:\spc{Y}\to \spc{Y}$ and $g_\infty\z\circ f_\infty\:\spc{X}\to \spc{X}$ are isometries.
Therefore, $f_\infty$ and $g_\infty$ are isometries as well.

\parit{\ref{SHORT.ex:GH-po:b}.} The implication 
\[|\spc{X}_n-\spc{X}_\infty|_{\GH}\to 0 
\quad\Rightarrow\quad 
\dist{\spc{X}_n}{\spc{X}_\infty}{\GH'}\to 0\]
follows from \ref{ex:eps-isom:GH>isom}. 

Now suppose $\dist{\spc{X}_n}{\spc{X}_\infty}{\GH'}\to 0$.
Show that $\{\spc{X}_n\}$ is a uniformly totally bonded family.

If $\dist{\spc{X}_n}{\spc{X}_\infty}{\GH}\not\to 0$, then we can pass to a subsequence such that $\dist{\spc{X}_n}{\spc{X}_\infty}{\GH}\ge\eps$ for some $\eps>0$.
By the Gromov selection theorem, we can assume that $\spc{X}_n$ converges in the sense of Gromov--Hausdorff.
From the first implication, the limit $\spc{X}_\infty'$ has to be isometric to $\spc{X}_\infty$;
on the other hand, $\dist{\spc{X}_\infty'}{\spc{X}_\infty}{\GH}\ge \eps$ --- a contradiction.

\parbf{\ref{ex:GH-urysohn}.}
Apply \ref{thm:compact-homogeneous} and \ref{prop:GH-with-fixed-Z}.


%% file: ultralimit-sol.tex
\refstepcounter{chapter}
\setcounter{eqtn}{0}

\parbf{\ref{ex:ultrakatetov}.} 
Let $F=\set{n\in \NN}{f(n)=n}$; we need to show that $\omega(F)=1$.

Consider an oriented graph $\Gamma$ with vertex set $\NN\setminus F$ such that $m$ is connected to $n$ if $f(m)=n$.
Show that each connected component of $\Gamma$ has at most one cycle.
Use it to subdivide vertices of $\Gamma$ into three sets $S_1$, $S_2$, and $S_3$ such that $f(S_i)\cap S_i=\emptyset$ for each $i$.

Conclude that $\omega(S_1)=\omega(S_2)=\omega(S_3)=0$ and hence \[\omega(F)=\omega(\NN\setminus(S_1\cup S_2\cup S_3))=1.\]

\parit{Source:} 
The presented proof was given by Robert Solovay \cite{solovay}, but
the key statement is due to Miroslav Katětov \cite{katetov}.

\parbf{\ref{ex:linear}.}
Choose a nonprincipal ultrafilter $\omega$ and set $L(\bm{s})=s_\omega$.
It remains to observe that $L$ is linear.

\parit{Remark.} 
This construction identifies ultrafilters with vectors in $(\ell^\infty)^*$.
Recall that $\ell^\infty=(\ell^1)^*$ and $\ell^1\subsetneq(\ell^\infty)^*$.
A principle ultrafilter is a basis vector in $\ell^1$; 
nonprincipal ultrafilters lie in $(\ell^\infty)^*\setminus\ell^1$.
The set of ultrafilters is the closure of basis vectors in $\ell^1$ with respect to weak*-topology on $(\ell^\infty)^*$.

\parbf{\ref{ex:ultrakatetov+}.}
Apply \ref{ex:ultrakatetov}.

\parbf{\ref{ex:lim(tree)}.}
Let $\gamma$ be a path from $p$ to $q$ in a metric tree $\spc{T}$.
Assume that $\gamma$ passes thru a point $x$ on distance $\ell$ from $[pq]$.
Then 
\[\length\gamma\ge \dist{p}{q}{}+2\cdot \ell.
\eqlbl{eq:+ell}\]

Suppose that $\spc{T}_n$ is a sequence of metric trees that $\omega$-converges to $\spc{T}_\omega$.
By \ref{obs:ultralimit-is-geodesic}, the space $\spc{T}_\omega$ is geodesic.

The uniqueness of geodesics follows from \ref{eq:+ell}.
Indeed, if for a geodesic $[p_\omega q_\omega]$ there is another geodesic $\gamma_\omega$ connecting its ends, then it has to pass thru a point $x_\omega\notin [p_\omega q_\omega]$.
Choose sequences $p_n,q_n,x_n\in\spc{T}_n$ such that $p_n\to p_\omega$, $q_n\to q_\omega$, and $x_n\to x_\omega$ as $n\to\omega$.
Then 
\begin{align*}
\dist{p_\omega}{q_\omega}{}&=\length\gamma\ge 
\\
&\ge\lim_{n\to\omega}(\dist{p_n}{x_n}{}+\dist{q_n}{x_n}{})\ge
\\
&\ge \lim_{n\to\omega}(\dist{p_n}{q_n}{}+2\cdot\ell_n)=
\\
&=\dist{p_\omega}{q_\omega}{}+2\cdot\ell_\omega.
\end{align*}
Since $x_\omega\notin [p_\omega q_\omega]$, we have that $\ell_\omega>0$ --- a contradiction.

It remains to show that any geodesic triangle $\spc{T}_\omega$ is a tripod.
Consider the sequence of centers of tripods $m_n$ for given sequences of points $x_n,y_n,z_n\in \spc{T}_n$.
Observe that its ultralimit $m_\omega$ is the center of a tripod with ends at $x_\omega,y_\omega,z_\omega\in \spc{T}_\omega$.

\parbf{\ref{ex:ultracompact}.}
Construct $\bm{X}$ and distance-preserving embeddings $\spc{X}_n\hookrightarrow\bm{X}$ that satisfy \ref{propery:GH}.
Given $x_\infty\in \spc{X}_\infty$, choose a sequence $x_n\in \spc{X}_n$ such that $x_n\to x_\infty$ in $\bm{X}$.
Let $x_\omega$ be the $\omega$-limit of the sequence $x_n$ in $\bm{X}$.
Note that $x_\omega\in \spc{X}_\infty$.
Show that the map $x_\infty\mapsto x_\omega$ is defined; that is, it does not depend on the choice of the sequence $x_n$.
Further, show that the map $x_\infty\mapsto x_\omega$ is an isometry of $\spc{X}_\infty$.
Make a conclusion.

\parbf{\ref{ex:ultrapower}.}
Further, we consider $\spc{X}$ as a subset of $\spc{X}^\omega$.

\parit{\ref{SHORT.ex:ultrapower:a}.} Follows directly from the definitions.

\parit{\ref{SHORT.ex:ultrapower:compact}.}
Suppose $\spc{X}$ compact.
Given a sequence $x_1,x_2,\dots{}\in\spc{X}$, denote its $\omega$-limit in $\spc{X}^\omega$ by $x^\omega$ and its $\omega$-limit in $\spc{X}$ by $x_\omega$.

Observe that $x^\omega=\iota(x_\omega)$.
Therefore, $\iota$ is onto.

If $\spc{X}$ is not compact, we can choose a sequence $x_n$ such that $\dist{x_m}{x_n}{}>\eps$ for fixed $\eps>0$ and all $m\ne n$.
Observe that
\[\lim_{n\to\omega}\dist{x_n}{y}{\spc{X}}\ge \tfrac\eps2\]
for any $y\in\spc{X}$.
It follows that $x_\omega$ lies at the distance $\ge\tfrac\eps2$ from~$\spc{X}$.

\parit{\ref{SHORT.ex:ultrapower:proper}.}
A sequence of points $x_n$ in $\spc{X}$ will be called $\omega$-bounded if there is a real constant $C$ such that
\[\dist{p}{x_n}{\spc{X}}\le C\] 
for $\omega$-almost all $n$.

The same argument as in \ref{SHORT.ex:ultrapower:compact} shows that any $\omega$-bounded sequence has its $\omega$-limit in $\spc{X}$.
Further, if $(x_n)$ is not  $\omega$-bounded, then 
\[\lim_{n\to\omega}\dist{p}{x_n}{\spc{X}}=\infty;\]
that is, $x_\omega$ does not lie in the metric component of $p$ in $\spc{X}^\omega$.

\parbf{\ref{ex:isom-ultrapower}.}
Let us show that cardinality of $\spc{X}^\omega$ is at least continuum ---
it is sufficient to construct a continuum family $\mathcal{A}$ sequences of points on $\spc{X}$ such that for any two sequences $(a_n)$ and $(b_n)$ in $\mathcal{A}$ the equality $a_n=b_n$ holds only for finitely many $n$.

To do this, let us identify points in $\spc{X}$ with nonnegative integers.
Consider the set $\mathcal{A}$ of all sequences $a_n$ such that $a_0=0$ and $a_{n+1}\z=a_n+\eps_n\cdot 2^n$ where $\eps_n\in\{0,1\}$ for any $n$.
Observe that $\mathcal{A}$ has cardinality continuum and distinct sequences in $\mathcal{A}$ have distinct $\omega$-limits.

Show and use that the spaces $\spc{X}^\omega$ and $(\spc{X}^\omega)^\omega$ have discrete metrics and both have cardinality at most continuum.

\parit{A more conceptual construction of $\mathcal{A}$.}
Choose a compact metric space $\spc{K}$ with continuum points, say $\spc{K}=[0,1]$.
Identify $\spc{X}$ with a dense subset of $\spc{K}$.
For any point $k\in \spc{K}$, choose a sequence $a_n\in \spc{X}$ that converges to $k$.
Observe that the family of all these sequences meet the required condition.

\parbf{\ref{ex:ultrapower(ultrapower)}.}
Choose a bijection $\iota\:\NN\to \NN\times \NN$.
Given a set $S\subset \NN$, consider the sequence $S_1$, $S_2,\dots$ of subsets in $\NN$ defined by $m\in S_n$ if $(m,n)\z=\iota(k)$ for some $k\in S$.
Set $\omega_1(S)=1$ if and only if $\omega(S_n)=1$ for $\omega$-almost all $n$.
It remains to check that $\omega_1$ meets the conditions of the exercise.

\parit{Comment.}
It turns out that $\omega_1\ne \omega$ for any $\iota$;
see the post of Andreas Blass \cite{blass}.

\parbf{\ref{ex:two-geodesics-in-ultrapower}.}
Arguing as in \ref{obs:ultrapower-is-geodesic}, we get a pair of points $x$ and $y$ in $\spc{X}$ such that
\[\dist{p}{x}{}+\dist{x}{y}{}+\dist{y}{q}{}=\dist{p}{q}{}\]
and there is no midpoint between $x$ and $y$ in $\spc{X}$
(possibly $p=x$ and $q=y$).
Note that it is sufficient to show that there is a continuum of distinct midpoints in $\spc{X}^\omega$ between $x$ and $y$ in $\spc{X}$.

Since $\spc{X}$ is a length space, we can choose a $\tfrac1n$-midpoint $m_n\in\spc{X}$ between $x$ and $y$.
Note that the sequence $m_n$ contains no converging subsequence.
Conclude that we may pass to a subsequence of $m_n$ such that $\dist{m_i}{m_j}{}>\eps$ for a fixed $\eps>0$ and any $i\ne j$.

Argue as in \ref{ex:isom-ultrapower} to show that there is a continuum of distinct $\omega$-limits of subsequences of $m_n$;
each such limit is a midpoint between $x$ and $y$.

\parbf{\ref{ex:notproper-limit}.} Consider the infinite metric $\spc{T}$ tree with unit edges shown
on the diagram.
Observe that $\spc{T}$ is proper.

\begin{Figure}
\vskip-0mm
\centering
\includegraphics{mppics/pic-605}
\end{Figure}

Consider the vertex $v_\omega=\lim_{n\to\omega}v_n$ in the ultrapower $\spc{T}^\omega$.
Observe that $\omega$ has an infinite degree.
Conclude that $\spc{T}^\omega$ is not locally compact.

\parbf{\ref{ex:ultraT}.}
Consider a product of an infinite sequence of two-point spaces.

\parit{Remark.}
There are such examples of geodesic spaces with a cocompact isometric action of a finitely generated group \cite{thomas-velickovic}.

\parbf{\ref{ex:Asym(Lob)}.} Assume $\spc{L}$ is the Lobachevsky plane.

\parit{\ref{SHORT.ex:Asym(Lob):metric-tree}.}
Show that there is $\delta>0$ such that sides of any geodesic triangle in $\spc{L}$ intersect a disk of radius $\delta$.
Conclude that any geodesic triangle in $\Asym\spc{L}$ is a tripod.

\parit{\ref{SHORT.ex:Asym(Lob):homogeneous}.} Observe that $\spc{L}$ is one-point-homogeneous and use it.

\parit{\ref{SHORT.ex:Asym(Lob):continuum}.} 
By \ref{SHORT.ex:Asym(Lob):homogeneous}, it is sufficient to show that $p_\omega$ has a continuum degree.

Choose distinct geodesics $\gamma_1,\gamma_2\:[0,\infty)\z\to L$ that start at a point $p$.
Show that the limits of $\gamma_1$ and $\gamma_2$ run in the different connected components of $(\Asym\spc{L})\setminus \{p_\omega\}$.
Since there is a continuum of distinct geodesics starting at $p$,
we get that the degree of $p_\omega$ is at least continuum.

On the other hand, the set of sequences of points in $\spc{L}$  has cardinality continuum.
In particular, the set of points in $\Asym\spc{L}$ has cardinality at most continuum.
It follows that the degree of any vertex is at most continuum.

The proof for the Lobachevsky space goes along the same lines.

For the infinite three-regular tree, part \ref{SHORT.ex:Asym(Lob):metric-tree} follows from \ref{ex:lim(tree)}.
The three-regular tree is only vertex-homogeneous; the latter is sufficient to prove \ref{SHORT.ex:Asym(Lob):homogeneous}.
No changes are needed in~\ref{SHORT.ex:Asym(Lob):continuum}.

\parit{Remark.}
The properties \ref{SHORT.ex:Asym(Lob):homogeneous} and \ref{SHORT.ex:Asym(Lob):continuum} describe the tree $\spc{T}$ up to isometry \cite{dyubina-polterovich}.
In particular, the asymptotic space of the Lobachevsky plane does not depend on the choice of the ultrafilter and the sequence $\lambda_n\z\to \infty$.

\parbf{\ref{ex:T(Sx[0,1]/Sx0)}.}
Denote by $o_\omega$ the point in $\T^\omega_o\spc{X}$ that corresponds to $o$.
Argue as in \ref{ex:Asym(Lob):continuum} to show that $\T^\omega_o\spc{X}\setminus \{o_\omega\}$ has a continuum of connected components.
Further, show that each connected component $\spc{W}_\alpha$ is isometric to $\RR\times (0,\infty)$ with the metric described by
\begin{align*}
&\dist{(x_1,t_1)}{(x_2,t_2)}{}=
\\
&\qquad=\min\{\,\dist{(x_1,t_1)}{(x_2,t_2)}{\RR^2},t_1+t_2\,\}.
\end{align*}

Conclude that the space $\T^\omega_o\spc{X}$ can be described as follows.
Prepare continuum copies $\spc{W}_\alpha$ as above;
denote by $(x,t)_\alpha$ the point in $\spc{W}_\alpha$ with coordinates $(x,t)$.
The tangent space is the disjoint union of single point $o_\omega$ and all $\spc{W}_\alpha$ with metric
such that $\dist{(x_1,t_1)_\alpha}{(x_2,t_2)_\alpha}{}$ is the same as in $\spc{W}_\alpha$ and for the remaining pairs, we have $\dist{o_\omega}{(x,t)_\alpha}{}=t$ and $\dist{(x_1,t_1)_\alpha}{(x_2,t_2)_\beta}{}=t_1+t_2$
if $\alpha\ne\beta$.

%% file: library.bib
@book {alexander-kapovitch-petrunin-2019,
    AUTHOR = {Alexander, S. and Kapovitch, V. and Petrunin, A.},
     TITLE = {An invitation to {A}lexandrov geometry: CAT(0) spaces},
      YEAR = {2019},
       URL = {https://doi.org/10.1007/978-3-030-05312-3},
}

@article {aronszajn-panitchpakdi,
    AUTHOR = {Aronszajn, N. and Panitchpakdi, P.},
     TITLE = {Extension of uniformly continuous transformations and hyperconvex metric spaces},
   JOURNAL = {Pacific J. Math.},
  FJOURNAL = {Pacific Journal of Mathematics},
    VOLUME = {6},
      YEAR = {1956},
     PAGES = {405--439},
      ISSN = {0030-8730},
   MRCLASS = {54.0X},
  MRNUMBER = {84762},
MRREVIEWER = {J. L. Kelley},
       URL = {http://projecteuclid.org/euclid.pjm/1103043960},
}

@article{bing-1948,
  title={A homogeneous indecomposable plane continuum},
  author={Bing, R. H.},
  year={1948},
  JOURNAL = {Duke Math. J.},
  VOLUME ={15},
  PAGES = {729--742}
}

@article {birkhoff,
    AUTHOR = {Birkhoff, G.},
     TITLE = {Metric foundations of geometry. {I}},
   JOURNAL = {Trans. Amer. Math. Soc.},
  FJOURNAL = {Transactions of the American Mathematical Society},
    VOLUME = {55},
      YEAR = {1944},
     PAGES = {465--492},
      ISSN = {0002-9947},
   MRCLASS = {48.0X},
  MRNUMBER = {10393},
MRREVIEWER = {L. M. Blumenthal},
       DOI = {10.2307/1990304},
       URL = {https://doi.org/10.2307/1990304},
}

@book {blaschke,
    AUTHOR = {Blaschke, W.},
     TITLE = {Kreis und {K}ugel},
      YEAR = {1916},
      addendum={Русский перевод: В. Бляшке, \emph{Круг и шар} (1967).}
}

@book {bridson-haefliger,
    AUTHOR = {Bridson, M. and Haefliger, A.},
     TITLE = {Metric spaces of non-positive curvature},
    SERIES = {Grundlehren der Mathematischen Wissenschaften},
    VOLUME = {319},
 %PUBLISHER = {Springer-Verlag, Berlin},
      YEAR = {1999},
 %    PAGES = {xxii+643},
      ISBN = {3-540-64324-9},
   MRCLASS = {53C23 (20F65 53C70 57M07)},
  MRNUMBER = {1744486},
MRREVIEWER = {Athanase Papadopoulos},
       DOI = {10.1007/978-3-662-12494-9},
       URL = {https://doi.org/10.1007/978-3-662-12494-9},
}

@MISC {buehler,
    TITLE = {Does there exist an isometry between $L^p$ and $\ell^p$?},
    AUTHOR = {Buehler, T.},
    HOWPUBLISHED = {MathOverflow},
   % NOTE = {URL:https://mathoverflow.net/q/112776 (version: 2017-04-13)},
    EPRINT = {https://mathoverflow.net/q/112776},
    URL = {https://mathoverflow.net/q/112776}
}

@book {burago-burago-ivanov,
    AUTHOR = {Burago, D. and Burago, Yu. and Ivanov, S.},
     TITLE = {A course in metric geometry},
    SERIES = {Graduate Studies in Mathematics},
    VOLUME = {33},
 %PUBLISHER = {American Mathematical Society, Providence, RI},
      YEAR = {2001},
    % PAGES = {xiv+415},
      ISBN = {0-8218-2129-6},
   MRCLASS = {53C23},
  MRNUMBER = {1835418},
MRREVIEWER = {Mario Bonk},
       DOI = {10.1090/gsm/033},
       URL = {https://doi.org/10.1090/gsm/033},
}

@book {busemann-1942,
    AUTHOR = {Busemann, H.},
     TITLE = {Metric methods in {F}insler spaces and in the foundations of geometry},
    SERIES = {Annals of Mathematics Studies, No. 8},
 %PUBLISHER = {Princeton University Press, Princeton, N. J.},
      YEAR = {1942},
   %  PAGES = {viii+243},
   MRCLASS = {48.0X},
  MRNUMBER = {0007251},
MRREVIEWER = {S. M. Ulam},
}

@incollection {cameron,
    AUTHOR = {Cameron, P.},
     TITLE = {The random graph},
 BOOKTITLE = {The mathematics of {P}aul {E}rd\H{o}s, {II}},
    SERIES = {Algorithms Combin.},
    VOLUME = {14},
     PAGES = {333--351},
 %PUBLISHER = {Springer, Berlin},
      YEAR = {1997},
   MRCLASS = {05C80},
  MRNUMBER = {1425227},
MRREVIEWER = {Tomasz J. \L uczak},
       DOI = {10.1007/978-3-642-60406-5_32},
       URL = {https://doi.org/10.1007/978-3-642-60406-5_32},
}

@book {deza-laurent,
    AUTHOR = {Deza, M. and Laurent, M.},
     TITLE = {Geometry of cuts and metrics},
    SERIES = {Algorithms and Combinatorics},
    VOLUME = {15},
 %PUBLISHER = {Springer-Verlag, Berlin},
      YEAR = {1997},
   %  PAGES = {xii+587},
      ISBN = {3-540-61611-X},
   MRCLASS = {52-02 (05B30 05C12 51K05 52B12 52C07 68R05 90C28)},
  MRNUMBER = {1460488},
MRREVIEWER = {Alexander I. Barvinok},
       DOI = {10.1007/978-3-642-04295-9},
       URL = {https://doi.org/10.1007/978-3-642-04295-9},
addendum={Русский перевод: М. Деза, M. Лоран, \emph{Геометрия разрезов и метрик} (2001).}
}

@article {dyubina-polterovich,
    AUTHOR = {Dyubina, A. and Polterovich, I.},
     TITLE = {Explicit constructions of universal {$\mathbb{R}$}-trees and
              asymptotic geometry of hyperbolic spaces},
   JOURNAL = {Bull. London Math. Soc.},
  FJOURNAL = {The Bulletin of the London Mathematical Society},
    VOLUME = {33},
      YEAR = {2001},
    NUMBER = {6},
     PAGES = {727--734},
      ISSN = {0024-6093},
   MRCLASS = {57M07 (20F67 53C23 54F50)},
  MRNUMBER = {1853785},
       DOI = {10.1112/S002460930100844X},
       URL = {https://doi.org/10.1112/S002460930100844X},
}

@book {federer,
    AUTHOR = {Federer, H.},
     TITLE = {Geometric measure theory},
    SERIES = {Die Grundlehren der mathematischen Wissenschaften, Band 153},
 %PUBLISHER = {Springer-Verlag New York Inc., New York},
      YEAR = {1969},
    % PAGES = {xiv+676},
   MRCLASS = {28.80 (26.00)},
  MRNUMBER = {0257325},
MRREVIEWER = {J. E. Brothers},
addendum={Русский перевод: Г. Федерер, \emph{Геометрическая теория меры} (1987).}
}

@article {frolik,
    AUTHOR = {Frolík, Z.},
     TITLE = {Concerning topological convergence of sets},
   JOURNAL = {Czechoslovak Math. J},
  FJOURNAL = {Czechoslovak Mathematical Journal},
    VOLUME = {10(85)},
      YEAR = {1960},
     PAGES = {168--180},
      ISSN = {0011-4642},
   MRCLASS = {54.00},
  MRNUMBER = {0116303},
MRREVIEWER = {Ky Fan},
}

@article{frechet,
  author={Fr{\'e}chet, M.},
  title={Sur quelques points du calcul fonctionnel},
  journal={Rendiconti del Circolo Matematico di Palermo (1884-1940)},
  volume={22},
  number={1},
  pages={1--72},
  year={1906},
  publisher={Springer Milan}
}

@incollection {gromov-1981-hyperbolic,
    AUTHOR = {M. {Gromov}},
     TITLE = {Hyperbolic manifolds, groups and actions},
 BOOKTITLE = {Riemann surfaces and related topics: {P}roceedings of the 1978
              {S}tony {B}rook {C}onference},
    SERIES = {Ann. of Math. Stud.},
    VOLUME = {97},
     PAGES = {183--213},
 %PUBLISHER = {Princeton Univ. Press, Princeton, N.J.},
      YEAR = {1981},
   MRCLASS = {53C15 (53C45 58F17)},
  MRNUMBER = {624814 (82m:53035)},
MRREVIEWER = {M. Rees},
}

@book {gromov-2007,
    AUTHOR = {Gromov, M.},
     TITLE = {Metric structures for {R}iemannian and non-{R}iemannian spaces},
    SERIES = {Modern Birkh\"{a}user Classics},
 %  EDITION = {English},
     % NOTE = {Based on the 1981 French original,  With appendices by M. Katz, P. Pansu and S. Semmes, Translated from the French by Sean Michael Bates},
% PUBLISHER = {Birkh\"{a}user Boston, Inc., Boston, MA},
      YEAR = {2007},
  %   PAGES = {xx+585},
   %   ISBN = {978-0-8176-4582-3; 0-8176-4582-9},
  % MRCLASS = {53C23 (53-02)},
 % MRNUMBER = {2307192},
}

@incollection {gross,
    AUTHOR = {Gross, O.},
     TITLE = {The rendezvous value of metric space},
 BOOKTITLE = {Advances in game theory},
     PAGES = {49--53},
 PUBLISHER = {Princeton Univ. Press, Princeton, N.J.},
      YEAR = {1964},
   MRCLASS = {90.70},
  MRNUMBER = {0162643},
}

@book {hausdorff,
    AUTHOR = {Hausdorff, F.},
     TITLE = {Grundz\"{u}ge der {M}engenlehre},
      YEAR = {1914},
      addendum={Русский перевод: Ф. Хаусдорф \emph{Теория множеств} (1937);
      English translation F. Hausdorff \emph{Set Theory} (1957).}
}

@Article{hilbert,
    Author = {D. {Hilbert}},
    Title = {{Ueber die gerade Linie als k\"urzeste Verbindung zweier Punkte.}},
    FJournal = {{Mathematische Annalen}},
    Journal = {{Math. Ann.}},
    ISSN = {0025-5831; 1432-1807/e},
    Volume = {46},
    Pages = {91--96},
    Year = {1895},
    Publisher = {Springer, Berlin/Heidelberg},
    DOI = {10.1007/BF02096204},
    MSC2010 = {51A05},
    Zbl = {26.0415.01},
    language={german},
hyphenation={german}
}

@article {ivanov-nikolaeva-tuzhilin,
    AUTHOR = {Ivanov, A. O. and Nikolaeva, N. K. and Tuzhilin, A. A.},
     TITLE = {The {G}romov--{H}ausdorff metric on the space of compact metric
              spaces is strictly intrinsic},
   JOURNAL = {Mat. Zametki},
  FJOURNAL = {Matematicheskie Zametki},
    VOLUME = {100},
      YEAR = {2016},
    NUMBER = {6},
     PAGES = {947--950},
      ISSN = {0025-567X},
   MRCLASS = {53C23 (58D17)},
  MRNUMBER = {3588919},
MRREVIEWER = {Pavel D. Andreev},
       DOI = {10.4213/mzm11411},
       URL = {https://doi.org/10.4213/mzm11411},
}

@article {isbell,
    AUTHOR = {Isbell, J. R.},
     TITLE = {Six theorems about injective metric spaces},
   JOURNAL = {Comment. Math. Helv.},
  FJOURNAL = {Commentarii Mathematici Helvetici},
    VOLUME = {39},
      YEAR = {1964},
     PAGES = {65--76},
      ISSN = {0010-2571},
   MRCLASS = {54.35},
  MRNUMBER = {182949},
       DOI = {10.1007/BF02566944},
       URL = {https://doi.org/10.1007/BF02566944},
}

@article {isbell2,
    AUTHOR = {Isbell, J. R.},
     TITLE = {Injective envelopes of {B}anach spaces are rigidly attached},
   JOURNAL = {Bull. Amer. Math. Soc.},
  FJOURNAL = {Bulletin of the American Mathematical Society},
    VOLUME = {70},
      YEAR = {1964},
     PAGES = {727--729},
      ISSN = {0002-9904},
   MRCLASS = {46.10},
  MRNUMBER = {184061},
MRREVIEWER = {A. P. Robertson},
       DOI = {10.1090/S0002-9904-1964-11192-7},
       URL = {https://doi.org/10.1090/S0002-9904-1964-11192-7},
}

@online{karlsson,
author={Karlsson, A.},
title={Ergodic theorems for noncommuting random products},
URL = {http://www.unige.ch/math/folks/karlsson/}
}

@misc{karlsson-2023,
      title={A metric fixed point theorem and some of its applications}, 
      author={Anders Karlsson},
      year={2023},
      eprint={2207.00963},
      archivePrefix={arXiv},
      primaryClass={math.FA}
}

@incollection {katetov,
    AUTHOR = {Katětov, M.},
     TITLE = {On universal metric spaces},
 BOOKTITLE = {General topology and its relations to modern analysis and
              algebra, {VI} ({P}rague, 1986)},
    SERIES = {Res. Exp. Math.},
    VOLUME = {16},
     PAGES = {323--330},
 PUBLISHER = {Heldermann, Berlin},
      YEAR = {1988},
   MRCLASS = {54E35 (54A25 54E40)},
  MRNUMBER = {952617},
MRREVIEWER = {Stoyan Nedev},
}

@article {KSS,
    AUTHOR = {Kirchheim, B. and Spadaro, E. and Sz{é}kelyhidi, L.},
     TITLE = {Equidimensional isometric maps},
   JOURNAL = {Comment. Math. Helv.},
  FJOURNAL = {Commentarii Mathematici Helvetici. A Journal of the Swiss
              Mathematical Society},
    VOLUME = {90},
      YEAR = {2015},
    NUMBER = {4},
     PAGES = {761--798},
      ISSN = {0010-2571},
   MRCLASS = {Preliminary Data},
  MRNUMBER = {3433279},
       DOI = {10.4171/CMH/370},
       URL = {http://dx.doi.org/10.4171/CMH/370},
}

@article {kleiner-leeb,
    AUTHOR = {Kleiner, B. and Leeb, B.},
     TITLE = {Rigidity of quasi-isometries for symmetric spaces and
              {E}uclidean buildings},
   JOURNAL = {Inst. Hautes \'{E}tudes Sci. Publ. Math.},
  FJOURNAL = {Institut des Hautes \'{E}tudes Scientifiques. Publications
              Math\'{e}matiques},
    NUMBER = {86},
      YEAR = {1997},
     PAGES = {115--197 (1998)},
      ISSN = {0073-8301},
   MRCLASS = {53C35 (20E42 20F32 53C23 57M07)},
  MRNUMBER = {1608566},
MRREVIEWER = {Lee Mosher},
       URL = {http://www.numdam.org/item?id=PMIHES_1997__86__115_0},
}

@article{kuratowski,
  title={Quelques probl{\`e}mes concernant les espaces m{\'e}triques non-s{\'e}parables},
  author={Kuratowski, C.},
  journal={Fundamenta Mathematicae},
  volume={25},
  number={1},
  pages={534--545},
  year={1935}
}

@article {lang-2013,
    AUTHOR = {Lang, U.},
     TITLE = {Injective hulls of certain discrete metric spaces and groups},
   JOURNAL = {J. Topol. Anal.},
  FJOURNAL = {Journal of Topology and Analysis},
    VOLUME = {5},
      YEAR = {2013},
    NUMBER = {3},
     PAGES = {297--331},
      ISSN = {1793-5253},
   MRCLASS = {20F65 (30L05 52C10)},
  MRNUMBER = {3096307},
       DOI = {10.1142/S1793525313500118},
       URL = {https://doi.org/10.1142/S1793525313500118},
}

@article {lang-pavon-zust,
    AUTHOR = {Lang, U. and Pav\'{o}n, M. and Z\"{u}st, R.},
     TITLE = {Metric stability of trees and tight spans},
   JOURNAL = {Arch. Math. (Basel)},
  FJOURNAL = {Archiv der Mathematik},
    VOLUME = {101},
      YEAR = {2013},
    NUMBER = {1},
     PAGES = {91--100},
      ISSN = {0003-889X},
   MRCLASS = {54E35 (51K05 54D35)},
  MRNUMBER = {3073668},
MRREVIEWER = {Yasunao Hattori},
       DOI = {10.1007/s00013-013-0535-y},
       URL = {https://doi.org/10.1007/s00013-013-0535-y},
}

@article {menger,
    AUTHOR = {Menger, K.},
     TITLE = {Untersuchungen \"{u}ber allgemeine {M}etrik},
   JOURNAL = {Math. Ann.},
  FJOURNAL = {Mathematische Annalen},
    VOLUME = {100},
      YEAR = {1928},
    NUMBER = {1},
     PAGES = {75--163},
      ISSN = {0025-5831},
   MRCLASS = {DML},
  MRNUMBER = {1512479},
       DOI = {10.1007/BF01448840},
       URL = {https://doi.org/10.1007/BF01448840},
}

@article {melleray,
    AUTHOR = {Melleray, J.},
     TITLE = {Some geometric and dynamical properties of the {U}rysohn
              space},
   JOURNAL = {Topology Appl.},
  FJOURNAL = {Topology and its Applications},
    VOLUME = {155},
      YEAR = {2008},
    NUMBER = {14},
     PAGES = {1531--1560},
      ISSN = {0166-8641},
   MRCLASS = {54E50 (54E45 54H11)},
  MRNUMBER = {2435148},
MRREVIEWER = {Oscar Valero},
       DOI = {10.1016/j.topol.2007.04.029},
       URL = {https://doi.org/10.1016/j.topol.2007.04.029},
}

@article {miesch,
    AUTHOR = {Miesch, B.},
     TITLE = {Gluing hyperconvex metric spaces},
   JOURNAL = {Anal. Geom. Metr. Spaces},
  FJOURNAL = {Analysis and Geometry in Metric Spaces},
    VOLUME = {3},
      YEAR = {2015},
    NUMBER = {1},
     PAGES = {102--110},
   MRCLASS = {54E35 (30L05)},
  MRNUMBER = {3349339},
MRREVIEWER = {Oleksiy A. Dovgoshey},
       DOI = {10.1515/agms-2015-0007},
       URL = {https://doi.org/10.1515/agms-2015-0007},
}

@article {miesch-pavon,
    AUTHOR = {Miesch, B. and Pav\'{o}n, M.},
     TITLE = {Weakly externally hyperconvex subsets and hyperconvex gluings},
   JOURNAL = {J. Topol. Anal.},
  FJOURNAL = {Journal of Topology and Analysis},
    VOLUME = {9},
      YEAR = {2017},
    NUMBER = {3},
     PAGES = {379--407},
      ISSN = {1793-5253},
   MRCLASS = {54E35 (46B07)},
  MRNUMBER = {3661648},
       DOI = {10.1142/S1793525317500145},
       URL = {https://doi.org/10.1142/S1793525317500145},
}

@MISC{nazarov,
    TITLE = {Intrinsic metric with no geodesics},
    AUTHOR = {Nazarov, F.},
    HOWPUBLISHED = {MathOverflow},
   % NOTE = {(version: 2010-02-18)},
    EPRINT = {http://mathoverflow.net/q/15720},
    URL = {http://mathoverflow.net/q/15720}
}

@book{petrunin-2022-PIGTIKAL,
      title={{PIGTIKAL} (puzzles in geometry that I know and love)}, 
      author={A. Petrunin},
      SERIES = {AMR Research Monographs, Volume 2},
      year={2022},
}

@misc{petrunin2020mnfld,
  doi = {10.48550/ARXIV.2010.10040},   
  author = {Petrunin, A.},
  title = {Metric geometry on manifolds: two lectures},  
  copyright = {Creative Commons Attribution 4.0 International},
      year={2020},
    eprint={2010.10040},
    archivePrefix={arXiv},
    primaryClass={math.DG}
}

@article {petrunin-stadler,
    AUTHOR = {Petrunin, A. and Stadler, S.},
     TITLE = {Metric-minimizing surfaces revisited},
   JOURNAL = {Geom. Topol.},
  FJOURNAL = {Geometry \& Topology},
    VOLUME = {23},
      YEAR = {2019},
    NUMBER = {6},
     PAGES = {3111--3139},
      ISSN = {1465-3060},
   MRCLASS = {53C23 (30L05 53C43 53C45)},
  MRNUMBER = {4039186},
       DOI = {10.2140/gt.2019.23.3111},
       URL = {https://doi.org/10.2140/gt.2019.23.3111},
}

@article {petrunin-yashinski,
    AUTHOR = {Petrunin, A. and Yashinski, A.},
     TITLE = {Piecewise isometric mappings},
   JOURNAL = {St. Petersburg Math. J.},
    VOLUME = {27},
      YEAR = {2016},
    NUMBER = {1},
     PAGES = {155–-175},
      ISSN = {0234-0852},
   MRCLASS = {53C45 (26A16)},
  MRNUMBER = {3443270},
MRREVIEWER = {V. Oproiu},
       DOI = {10.1090/spmj/1381},
       URL = {https://doi.org/10.1090/spmj/1381},
}

@book {pogorelov,
    AUTHOR = {Погорелов, А. В.},
      TITLE ={Четвертая проблема {Г}ильберта},
            YEAR = {1974},
    addendum={English translation: A. Pogorelov, \emph{Hilbert's fourth problem} (1979).}
}

@article {przeslawski-yost,
    AUTHOR = {Przes\l awski, K. and Yost, D.},
     TITLE = {Continuity properties of selectors and {M}ichael's theorem},
   JOURNAL = {Michigan Math. J.},
  FJOURNAL = {Michigan Mathematical Journal},
    VOLUME = {36},
      YEAR = {1989},
    NUMBER = {1},
     PAGES = {113--134},
      ISSN = {0026-2285},
   MRCLASS = {49A50 (28B20 47H99 54C65)},
  MRNUMBER = {989940},
MRREVIEWER = {Ewa Bednarczuk},
       DOI = {10.1307/mmj/1029003885},
       URL = {https://doi.org/10.1307/mmj/1029003885},
}

@article {hu-kirk,
    AUTHOR = {Hu, T. and Kirk, W. A.},
     TITLE = {Local contractions in metric spaces},
   JOURNAL = {Proc. Amer. Math. Soc.},
  FJOURNAL = {Proceedings of the American Mathematical Society},
    VOLUME = {68},
      YEAR = {1978},
    NUMBER = {1},
     PAGES = {121--124},
      ISSN = {0002-9939},
   MRCLASS = {54E40 (54H25)},
  MRNUMBER = {464180},
MRREVIEWER = {Harold Bell},
       DOI = {10.2307/2040922},
       URL = {https://doi.org/10.2307/2040922},
}

@MISC {lowther,
    TITLE = {On the global structure of the {G}romov--{H}ausdorff metric space},
    AUTHOR = {G. Lowther},
    HOWPUBLISHED = {MathOverflow},
    %NOTE = {(version: 2015-08-21)},
    EPRINT = {https://mathoverflow.net/q/212608},
   % URL = {https://mathoverflow.net/q/212608}
}

@misc{petrunin-zamorabarrera,
title={What is differential geometry: curves and surfaces}, 
author={Petrunin, A. and Zamora Barrera, S.},
year={2020},
      eprint={2012.11814},
      archivePrefix={arXiv},
      primaryClass={math.HO}
}

@article {ivanov-tuzhilin,
    AUTHOR = {Ivanov, A. O. and Tuzhilin, A. A.},
     TITLE = {Isometry group of {G}romov--{H}ausdorff space},
   JOURNAL = {Mat. Vesnik},
  FJOURNAL = {Matematichki Vesnik},
    VOLUME = {71},
      YEAR = {2019},
    NUMBER = {1-2},
     PAGES = {123--154},
      ISSN = {0025-5165},
   MRCLASS = {53C23 (54E45)},
  MRNUMBER = {3895911},
MRREVIEWER = {Nan Li},
}

@article {rudin,
    AUTHOR = {Rudin, W.},
     TITLE = {Homogeneity problems in the theory of \v{C}ech compactifications},
   JOURNAL = {Duke Math. J.},
  FJOURNAL = {Duke Mathematical Journal},
    VOLUME = {23},
      YEAR = {1956},
     PAGES = {409--419},
      ISSN = {0012-7094},
   MRCLASS = {54.0X},
  MRNUMBER = {0080902},
MRREVIEWER = {L. Gillman},
       URL = {http://projecteuclid.org/euclid.dmj/1077466953},
}

@online{solovay,
  title={Maps preserving measures},
  author={Solovay, R. M.},
  url={https://math.berkeley.edu/~solovay/Preprints/Rudin_Keisler.pdf},
  year={2011}
}

@article {thomas-velickovic,
    AUTHOR = {Thomas, S. and Velickovic, B.},
     TITLE = {Asymptotic cones of finitely generated groups},
   JOURNAL = {Bull. London Math. Soc.},
  FJOURNAL = {The Bulletin of the London Mathematical Society},
    VOLUME = {32},
      YEAR = {2000},
    NUMBER = {2},
     PAGES = {203--208},
      ISSN = {0024-6093},
   MRCLASS = {20F65},
  MRNUMBER = {1734187},
MRREVIEWER = {Thomas Delzant},
       DOI = {10.1112/S0024609399006621},
       URL = {https://doi.org/10.1112/S0024609399006621},
}

@article {tits,
    AUTHOR = {Tits, J.},
     TITLE = {Sur certaines classes d'espaces homog\`enes de groupes de {L}ie},
   JOURNAL = {Acad. Roy. Belg. Cl. Sci. M\'{e}m. Coll. in 8$^\circ$},
  FJOURNAL = {Acad\'{e}mie Royale de Belgique. Classe des Sciences. M\'{e}moires.
              Collection in-8$^\circ$. Koninklijke Belgische Academie.
              Klasse der Wetenschappen. Verhandelingen. Verzameling
              in-8$^\circ$},
    VOLUME = {29},
      YEAR = {1955},
    NUMBER = {3},
     PAGES = {268},
      ISSN = {0365-0936},
   MRCLASS = {17.0X},
  MRNUMBER = {76286},
MRREVIEWER = {L. Auslander},
}

@article{urysohn,
  title={Sur un espace m{\'e}trique universel},
  author={Urysohn, P.},
  journal={Bull. Sci. Math},
  volume={51},
  number={2},
  pages={43--64},
  year={1927},
  addendum ={Русский перевод в П. С. Урысон \textit{Труды по топологиии другим областям математики}, Том II, (1951) 747---777.}
}

@article {uspenskij,
    AUTHOR = {Uspenskij, V.},
     TITLE = {The {U}rysohn universal metric space is homeomorphic to a {H}ilbert space},
   JOURNAL = {Topology Appl.},
  FJOURNAL = {Topology and its Applications},
    VOLUME = {139},
      YEAR = {2004},
    NUMBER = {1-3},
     PAGES = {145--149},
      ISSN = {0166-8641},
   MRCLASS = {54F65 (54E50)},
  MRNUMBER = {2051102},
MRREVIEWER = {K. D. Magill, Jr.},
       DOI = {10.1016/j.topol.2003.09.008},
       URL = {https://doi.org/10.1016/j.topol.2003.09.008},
}

@article {vaisala,
    AUTHOR = {V\"{a}is\"{a}l\"{a}, J.},
     TITLE = {A proof of the {M}azur--{U}lam theorem},
   JOURNAL = {Amer. Math. Monthly},
  FJOURNAL = {American Mathematical Monthly},
    VOLUME = {110},
      YEAR = {2003},
    NUMBER = {7},
     PAGES = {633--635},
      ISSN = {0002-9890},
   MRCLASS = {46B20 (39B52)},
  MRNUMBER = {2001155},
MRREVIEWER = {Szymon W\polhk asowicz},
       DOI = {10.2307/3647749},
       URL = {https://doi.org/10.2307/3647749},
}

@article {vershik,
    AUTHOR = {Vershik, A. M.},
     TITLE = {Random metric spaces and universality},
   JOURNAL = {Uspekhi Mat. Nauk},
  FJOURNAL = {Uspekhi Matematicheskikh Nauk},
    VOLUME = {59},
      YEAR = {2004},
    NUMBER = {2(356)},
     PAGES = {65--104},
      ISSN = {0042-1316},
   MRCLASS = {60B99 (54E35)},
  MRNUMBER = {2086637},
MRREVIEWER = {O. Lipovan},
       DOI = {10.1070/RM2004v059n02ABEH000718},
       URL = {https://doi.org/10.1070/RM2004v059n02ABEH000718},
}

@article {wijsman,
    AUTHOR = {Wijsman, R. A.},
     TITLE = {Convergence of sequences of convex sets, cones and functions. {II}},
   JOURNAL = {Trans. Amer. Math. Soc.},
  FJOURNAL = {Transactions of the American Mathematical Society},
    VOLUME = {123},
      YEAR = {1966},
     PAGES = {32--45},
      ISSN = {0002-9947},
   MRCLASS = {52.30},
  MRNUMBER = {196599},
MRREVIEWER = {Z. Frol\'{\i}k},
       DOI = {10.2307/1994611},
       URL = {https://doi.org/10.2307/1994611},
}

@MISC{blass,
    TITLE = {Selective ultrafilter and bijective mapping},
    AUTHOR = {A. Blass},
    HOWPUBLISHED = {MathOverflow},
    EPRINT = {https://mathoverflow.net/q/324261},
    URL = {https://mathoverflow.net/q/324261}
}

@MISC {petrunin-2022-MO-center,
    TITLE = {Center of convex figure},
    AUTHOR = {A. Petrunin},
    HOWPUBLISHED = {MathOverflow},
    %NOTE = {URL:https://mathoverflow.net/q/432694 (version: 2022-10-18)},
    EPRINT = {https://mathoverflow.net/q/432694},
    URL = {https://mathoverflow.net/q/432694}
}

@MISC {petrunin-435157,
    TITLE = {Short selection in the space of subsets},
    AUTHOR = {A. Petrunin},
    HOWPUBLISHED = {MathOverflow},
    %NOTE = {URL:https://mathoverflow.net/q/435157 (version: 2022-11-24)},
    EPRINT = {https://mathoverflow.net/q/435157},
    URL = {https://mathoverflow.net/q/435157}
}

@book{alexander-kapovitch-petrunin-2025,
    title={Alexandrov geometry: foundations},
    author={Alexander, S. and  Kapovitch, V. and Petrunin, A.},
    year={2022},
    eprint={1903.08539},
    archivePrefix={arXiv},
    primaryClass={math.DG}
}

@misc{berestovskii-nikonorov,
  doi = {10.48550/ARXIV.2206.13096},  
  url = {https://arxiv.org/abs/2206.13096},  
  author = {Berestovskii, V. N. and Nikonorov, Y. G.},
  title = {On $m$-point homogeneous polytopes in Euclidean spaces},  
  publisher = {arXiv},  
  year = {2022},  
  copyright = {arXiv.org perpetual, non-exclusive license},
    eprint={2206.13096},
    archivePrefix={arXiv},
    primaryClass={math.MG}
}

@misc{zuest,
      title={The Riemannian hemisphere is almost calibrated in the injective hull of its boundary}, 
      author={Züst, R.},
      year={2021},
      eprint={2104.04498},
      archivePrefix={arXiv},
      primaryClass={math.DG}
}

@misc{miesch-pavon2016,
  title = {Ball intersection properties in metric spaces},
  author = {Miesch, B. and Pavón, M.},  
  year = {2016},
  eprint={1610.03307},
  archivePrefix={arXiv},
  primaryClass={math.MG},
  addendum ={to appear in J. Topol. Anal.}
}

@misc{lebedeva-petrunin2211.09671,
  title = {All-set-homogeneous spaces},
  author = {Lebedeva, N. and Petrunin, A.},
  year = {2022}, 
  eprint={2211.09671},
  archivePrefix={arXiv},
  primaryClass={math.MG},
  copyright = {Creative Commons Zero v1.0 Universal},
  addendum ={to appear in St. Petersburg Math. J.}
}

@misc{nabutovsky,
    title={Linear bounds for constants in Gromov's systolic inequality and related results},
    author={A. Nabutovsky},
    year={2019},
    eprint={1909.12225},
    archivePrefix={arXiv},
    primaryClass={math.MG}
}

@MISC {bilokopytov,
    TITLE = {Is it possible to connect every compact set?},
    AUTHOR = {Bilokopytov, E.},
    HOWPUBLISHED = {MathOverflow},
    %NOTE = {URL:https://mathoverflow.net/q/359390 (version: 2020-05-05)},
    EPRINT = {https://mathoverflow.net/q/359390},
    URL = {https://mathoverflow.net/q/359390}
}

@MISC {petrunin-431426,
    TITLE = {$m$-point-homogeneous, but not $(m+1)$-point-homogeneous},
    AUTHOR = {A. Petrunin},
    HOWPUBLISHED = {MathOverflow},
    %NOTE = {URL:https://mathoverflow.net/q/431426 (version: 2022-11-12)},
    EPRINT = {https://mathoverflow.net/q/431426},
    URL = {https://mathoverflow.net/q/431426}
}
